%% file: FESTUNG-Pt3.tex
\documentclass[%
3p,%
preprint,
fleqn, 
times
]{elsarticle}

\makeindex

\usepackage{amsmath}
\usepackage{amssymb}
\usepackage[english=american]{csquotes} 
\usepackage{graphicx}
\usepackage[english]{babel} 
\usepackage{color}
\usepackage[table]{xcolor}
\usepackage{textcomp} 
\usepackage{listings}
  \definecolor{keywordcolor}{rgb}{0, 0.25, 0.5}
  \definecolor{commentcolor}{rgb}{0.2, 0.5, 0.2}
  \definecolor{stringcolor}{rgb}{0.5, 0.5, 0.2}
\newcommand*{\code}[1]{\text{\lstinline[basicstyle=\ttfamily\small]{#1}}}

\usepackage{caption}
\usepackage{subcaption}

\usepackage{microtype}
\usepackage{hyperref}
\usepackage{tabularx}             
  \newcolumntype{R}{>{\raggedleft\arraybackslash}X}  \newcolumntype{L}{>{\raggedright\arraybackslash}X}  \newcolumntype{C}{>{\centering\arraybackslash}X}
\usepackage{booktabs} 

\usepackage{pgfplots} 
\usepackage{pgfplotstable} 
\pgfplotsset{compat=newest}

\usepackage[xcolor=dvipdf,authormarkup=none]{changes} 
  \definechangesauthor[color=red,name=Balthasar]{BR}
  \definechangesauthor[color=blue,name=Alex]{AJ}
  \definechangesauthor[color=green,name=Vadym]{VA}

\usepackage[normalem]{ulem}
\usepackage{dsfont} 

\input{mycommands.tex}
\input{tableStyles.tex} 

\renewcommand{\vec}[1]{\mathbf{#1}}
\newtheorem{remark}{Remark}

\title{FESTUNG: A~\MatOct~toolbox for the discontinuous Galerkin method. Part III: Hybridized discontinuous Galerkin (HDG) formulation}

\author[Hasselt]{Alexander Jaust}
  \ead{alexander.jaust@uhasselt.be}
\author[FAU]{Balthasar Reuter}
  \ead{reuter@math.fau.de}
\author[AWI,FAU]{Vadym Aizinger\corref{cor}}
  \ead{vadym.aizinger@awi.de}
\author[Hasselt]{Jochen Sch\"utz}
  \ead{jochen.schuetz@uhasselt.be}
\author[FAU]{Peter Knabner}
  \ead{knabner@math.fau.de}

\address[Hasselt]{UHasselt -- Hasselt University, Faculty of Sciences, Agoralaan Gebouw D, 3590 Diepenbeek, Belgium}
\address[FAU]{Friedrich--Alexander University of Erlangen--N\"urnberg, Department of Mathematics, 
Cauerstra{\ss}e~11, 91058~Erlangen, Germany}
\address[AWI]{Alfred Wegener Institute, Helmholtz Centre for Polar and Marine Research, Am Handelshafen 12, 27570 Bremerhaven, Germany}

\cortext[cor]{Corresponding author}

\begin{document}

\begin{abstract}
The third paper in our~series on open source \MatOct~implementation of the~discontinuous Galerkin~(DG) method(s) focuses on a~hybridized formulation. 
The main aim of this ongoing work is to develop rapid prototyping techniques covering a~range of standard DG methodologies and suitable for small to medium sized applications. 
Our FESTUNG package relies on fully vectorized matrix\,/\,vector operations throughout, and all details of the implementation are fully documented. 
Once again, great care is taken to maintain a~direct mapping between discretization terms and code routines as well as to ensure full compatibility to \Octave{}.
The current work formulates a~hybridized DG scheme for a~linear advection problem, describes hybrid approximation spaces on the mesh skeleton, and compares the performance of this discretization to the standard (element-based) DG method for different polynomial orders.
\end{abstract}


\maketitle

\section{Introduction}
\input{sections/introduction}

\section{Model problem}
\label{sec:model}
\input{sections/model}

\section{Discretization}
\label{sec:discretization}
\input{sections/discretization}

\section{Implementation}
\label{sec:implementation}
\input{sections/implementationNew}

\section{Numerical results}
\label{sec:results}
\input{sections/results}

\section{Register of routines}
\label{sec:registerroutines}
\input{sections/registerofroutines}

\section{Conclusion and Outlook}\label{sec:conclusion}

The third work in our series introduces a~hybridized DG formulation and expands our FESTUNG framework to include edge-based degrees of freedom as well as implicit time stepping schemes. 
A~comparison of the results obtained using HDG to those from the standard DG method indicates little difference in solution quality (except in the lowest order approximation $p=0$ case, where the unhybridized DG held a~slight edge) and confirms the well-known computational cost advantage of the HDG method for higher order polynomial spaces.
Since HDG implementations spend a~significant portion of computing time in element-local solves, this advantage over standard DG discretizations becomes even more pronounced for parallel implementations using the distributed memory programming paradigm. 
Our future work plans include more complex systems of PDEs and coupled multi-physics applications.

\section*{Acknowledgments}
The stay of A.~Jaust at the University of Erlangen-N\"urnberg was supported by the Research Foundation - Flanders (FWO) with a grant for a short study visit abroad and by the Special Research Fund (BOF) of Hasselt University.

\section*{Index of notation}
\noindent
\begin{small}
\begin{tabularx}{\linewidth}{@{}lX@{}}\toprule
\textbf{Symbol}    & \textbf{Definition}\\\midrule
$\diag(\vecc{A},\vecc{B})$ & $\coloneqq \begin{bmatrix}\vecc{A}&\quad\\\quad&\vecc{B}\end{bmatrix}$, block-diagonal matrix with blocks~$\vecc{A}$, $\vecc{B}$.\\
$\vec{a}\cdot\vec{b}$& $\coloneqq \sum_{m=1}^2a_mb_m$, Euclidean scalar product in~$\IR^2$.\\
$\grad$             & $\coloneqq \transpose{[\partial_{x^1}, \partial_{x^2}]}$, spatial gradient in the physical domain~$\Omega$.\\
$\circ$             & Composition of functions or Hadamard product.\\
$\otimes$,\; $\otimes_\mathrm{V}$& Kronecker product, `vectorial' Kronecker product.\\
$\alpha$            & Stabilization parameter.\\
$\betaMap$          & Mapping from $[0,1]$ to $[0,1]$ that adapts the edge orientation.\\
$c$                 & Concentration (scalar-valued unknown).\\
$c_{\partial\Omega}$& Concentration on the boundary $\partial\Omega$.\\
$c^0$,\; $c_\mathrm{D}$& Concentration prescribed at initial time~$t=0$ and on the inflow Dirichlet boundary.\\
$\vec{C}$           & $\in\IR^{KN}$, representation vector of~$c_h\in\IP_p(\setT_h)$ with respect to $\{\varphi_{kj}\}$.\\
$\delta_\mathrm{[condition]}$ & $\coloneqq \{1~\text{if condition is true, 0~otherwise}\}$, Kronecker delta.\\
$E_{kn}$, $\hat{E}_n$ & $n$th edge of the physical triangle~$T_k$, $n$th edge of the reference triangle~$\hat{T}$.\\
$E_{\kEdge}$        & $\kEdge$th physical edge of the mesh.\\
$\setE$,\; $\setT$  & Sets of edges and triangles.\\
$\setE_\mathrm{int}$,\; $\setE_\mathrm{ext}$,\; $\setE_\mathrm{in}$,\; $\setE_\mathrm{out}$ & Sets of interior, boundary, inflow, and outflow edges.\\
$\vec{f}$,\; $\hat{\vec{f}}$ & Advective flux function, numerical flux.\\
$\vec{F}_k$         & Affine mapping from $\hat{T}$ to $T_k$.\\
$\Gamma$            & Trace of the mesh.\\
$\ghat_n$           & Mapping from $[0,1]$ to $\hat{E}_n$.\\
$h$                 & Mesh fineness of $\setT_h$.\\
$h_T$               & $\coloneqq \diam(T)$, diameter of triangle~$T\in\setT_h$.\\
$\source(t,\vec{x})$      & Source\,/\,sink (scalar-valued coefficient function).\\ 
$J$                 & $\coloneqq (t^0,t_\mathrm{end})$, open time interval.\\
$K$,\; $\overline{K}$ & Number of triangles, number of edges.\\
$\kappa(\kEdge, l)$ & Index mapping from physical edge to neighboring physical triangles $T_{k^-}, T_{k^+}$.\\
$\lambda_h$         & Scalar-valued unknown on the trace of the mesh.\\
$\boldsymbol\Lambda$& $\in\IR^{\overline{K}\overline{N}}$, representation vector of~$\lambda_h\in\IP_p(\setE_h)$ with respect to $\{\mu_{\kEdge j}\}$.\\
$\mu_{\kEdge i}$,\; $\hat{\mu}_i$ & $i$th hierarchical basis function on $E_{\kEdge}$, $i$th hierarchical basis function on $[0,1]$.\\
$\boldsymbol{\nu}_{T}$,\; $\boldsymbol{\nu}_{kn}$,\; $\boldsymbol{\nu}_{\kEdge}$& Unit normal on~$\partial T$ pointing outward of~$T$, unit normal on~$E_{kn}$, unit normal on~$E_{\kEdge}$.\\
$N$                 & $\coloneqq (p+1)(p+2)/2$, number of local degrees of freedom of~$\IP_p(T)$.\\
$\overline{N}$      & $\coloneqq p+1$, number of local degrees of freedom of~$\IP_p(E)$.\\
$\omega_r$          & Quadrature weight associated with~$\hat{\vec{q}}_r$.\\
$\Omega$,\; $\partial\Omega$,\; $\partial\Omega_\mathrm{in}$,\; $\partial\Omega_\mathrm{out}$         & spatial domain in two dimensions, boundary of $\Omega$, inflow- and outflow boundaries,  $\partial\Omega = \partial\Omega_\mathrm{in}\cup\partial\Omega_\mathrm{out}$.\\
$p$                 & $= (\sqrt{8N+1}-3)/2$, polynomial degree.\\
$\vphi_{ki}$,\;  $\hat{\vphi}_i$      & $i$th hierarchical basis function on~$T_k$, $i$th hierarchical basis function on~$\hat{T}$.\\
$\IP_p(T)$          & Space of polynomials on~$T\in\setT_h$ of degree at most~$p$.\\
$\IP_p(\setT_h)$    & $\coloneqq \{ w_h:\overline{\Omega}\rightarrow \IR\,;\forall T\in\setT_h,\,   {w_h}|_T\in\IP_p(T)\}$.\\
$\IP_p(\setE_h)$    & $\coloneqq \{ w_h:\Gamma\rightarrow \IR\,;\forall E\in\setE_h,\,   {w_h}|_E\in\IP_p(E)\}$.\\
$R$,\; $\hat{\vec{q}}_r$   & Number of quadrature points, $r$th quadrature point in~$\hat{T}$.\\
$\rho(k,n)$         & Index mapping from edge $E_{kn}$ to physical edge $E_{\kEdge}$.\\
$t$                 & Time variable.\\
$t^n$               & $n$th time level.\\
$t_\mathrm{end}$    & End time.\\
$\mapEE_{n^-n^+}$   & Mapping from $\hat{E}_{n^-}$ to $\hat{E}_{n^+}$.\\
$\Delta t^n$        & $\coloneqq t^{n+1} - t^n$, time step size.\\
%
$T_k$,\; $\partial T_k$               & $k$th physical triangle, boundary of~$T_k$.\\
$\hat{T}$           & Bi-unit reference triangle.\\
$\vec{u}$           & Velocity (vector-valued coefficient function).\\
$\vec{x}$           & $=\transpose{[x^1,x^2]}$, space variable in the physical domain~$\Omega$.\\
$\hat{\vec{x}}$     & $=\transpose{[\hat{x}^1, \hat{x}^2]}$, space variable in the reference triangle~$\hat{T}$.\\
\bottomrule
\end{tabularx}
\end{small}

\bibliographystyle{elsarticle-num}
\bibliography{literature}

\end{document}

%% file: mycommands.tex
\newcommand*{\divergence}{\nabla \cdot}

\newcommand*{\vecu}{\vec{u}}

\newcommand{\overbar}[1]{\mkern 1.5mu\overline{\mkern-1.5mu#1\mkern-1.5mu}\mkern 1.5mu}

\newcommand*{\lambdah}{\lambda_{h}}

\newcommand*{\varphih}{\varphi_{h}}
\newcommand*{\muh}{\mu_{h}}

\newcommand*{\dt}{\Delta t}

\newcommand*{\parT}{\partial_{t}}
\newcommand*{\parXm}{\partial_{x^{m}}}
\newcommand*{\parX}[1]{\partial_{x^{#1}}}
\newcommand*{\parXhat}[1]{\partial_{\hat{x}^{#1}}}
\newcommand*{\parXmHat}{\parXhat{m}}
\newcommand*{\vx}{\vec{x}}

\newcommand*{\vxhat}{\vec{\hat{x}}}

\newcommand*{\ctx}{c\,(t, \vx)}

\newcommand*{\vu}{\vec{u}}
\newcommand*{\sumEl}{\sum_{j=1}^{N}}
\newcommand*{\sumFa}{\sum_{j=1}^{\Nedge}}

\newcommand*{\Ckj}{C_{kj}}
\newcommand*{\Lkj}{\Lambda_{kj}}

\newcommand*{\sumElLkj}{\sumFa \Lkj}
\newcommand*{\sumLkj}{\sumFa \Lkj}

\newcommand*{\sumElCkj}{\sumEl \Ckj}
\newcommand*{\sumElCkjt}{\sumElCkj(t)}

\newcommand*{\pkj}{\varphi_{kj}}
\newcommand*{\pki}{\varphi_{ki}}

\newcommand*{\muElementView}[1]{\mu_{kn#1}}
\newcommand*{\muiElementView}{\muElementView{i}}
\newcommand*{\mujElementView}{\muElementView{j}}

\newcommand*{\muFaceView}[1]{\mu_{\kEdge#1}}

\newcommand*{\mujFaceView}{\muFaceView{j}}

\newcommand*{\mkj}{\mu_{\kEdge{}j}}
\newcommand*{\mki}{\mu_{\kEdge{}i}}

\newcommand*{\Tk}{{T_{k}}}
\newcommand*{\TkBoundary}{\partial{\Tk}}
\newcommand*{\intTk}{\int_\Tk^{}}
\newcommand*{\sintTk}{\int_{\partial \Tk}^{}}
\newcommand*{\sintTkNB}{\int_{\partial \Tk\setminus \partial \Omega}^{}}
\newcommand*{\sintTkB}{\int_{\TkBoundary \cap \partial \Omega}^{}}

\newcommand*{\sintErefk}{\int_{\Eref_k}^{}}

\newcommand*{\dx}{\text{d}\vec{x}}
\newcommand*{\dxhat}{\text{d}\vec{\hat{x}}}
\newcommand*{\ds}{\text{d}s}

\newcommand*{\ch}{c_h}
\newcommand*{\cmh}{c^{-}_h}
\newcommand*{\cph}{c^{+}_h}

\newcommand*{\phih}{\varphi_{h}}

\newcommand*{\summ}{\sum_{m=1}^{2}}

\newcommand*{\mh}{\mu_h}
\newcommand*{\vNormal}{\boldsymbol{\nu}}

\newcommand*{\fluxNoDep}{\vec{f}}
\newcommand*{\flux}[1]{\fluxNoDep\,\left(#1\right)}

\newcommand*{\numFlux}{\hat{\vec{f}}}
\newcommand*{\numFluxDG}{\hat{\vec{f}}_{DG}}

\newcommand*{\Nedge}{\overbar{N}}
\newcommand*{\Kedge}{\overbar{K}}
\newcommand*{\kEdge}{\bar{k}}

\newcommand*{\R}{\mathbb{R}}

\newcommand*{\textInflow}{\text{on } \partial \Omega_{\text{in}}}
\newcommand*{\textOutflow}{\text{on } \partial \Omega_{\text{out}}}

\newcommand*{\Ekn}{E_{kn}}

\newcommand*{\Ekbar}{E_{\kEdge}}

\newcommand*{\Eref}{\hat{E}}
\newcommand*{\Erefn}{\Eref_{n}}
\newcommand*{\Ehatn}{\Erefn}

\newcommand*{\inE}{\setE_{\text{int}}}
\newcommand*{\bcE}{\setE_{\text{ext}}}

\newcommand*{\bcEin}{\setE_{\text{in}}}
\newcommand*{\bcEout}{\setE_{\text{out}}}

\newcommand*{\sumInE}{\sum_{\Ekn \in \TkBoundary \cap \inE }}

\newcommand*{\sintEkn}{\int_{\Ekn}}

\newcommand*{\oneDBasisHat}[1]{\hat{\mu}_{#1}(s)}

\newcommand*{\Ci}{\boldsymbol{C}^{(i)}}
\newcommand*{\Li}{\boldsymbol{\Lambda}^{(i)}}

\newcommand*{\matB}{\mat{B}}

\newcommand*{\matM}{\mat{M}}
\newcommand*{\matMbar}{\mat{\overbar{M}}}
\newcommand*{\matMhat}{\mat{\hat{M}}}

\newcommand*{\matK}{\mat{K}}

\newcommand*{\vecK}{\vec{K}}

\newcommand*{\matGm}{\mat{G}^{m}}
\newcommand*{\matG}{\mat{G}}
\newcommand*{\matGhat}{\mat{\hat{G}}}

\newcommand*{\matR}{\mat{R}}
\newcommand*{\matRhat}{\mat{\hat{R}}}
\newcommand*{\matRmu}{\matR_{\mu}}
\newcommand*{\matRmuHat}{\matRhat_{\mu}}
\newcommand*{\matRphi}{\matR_{\varphi}}

\newcommand*{\vecFphiD}{\vec{F}_{\varphi,\mathrm{in}}}

\newcommand*{\vecU}{\vec{U}}

\newcommand*{\matT}{\mat{T}}

\newcommand*{\matS}{\mat{S}}
\newcommand*{\matShat}{\mat{\hat{S}}}

\newcommand*{\matSout}{\mat{S}_{\text{out}}}

\newcommand*{\matMphi}{\matM_{\varphi}}

\newcommand*{\matMmuBar}{\matMbar_{\mu}}
\newcommand*{\matMmuTilde}{\tilde{\matM}_{\mu}}

\newcommand*{\vecKmuD}{\vecK_{\mu,\mathrm{in}}}
\newcommand*{\matKmuOut}{\matK_{\mu,\mathrm{out}}}

\newcommand*{\vecQ}{\vec{Q}}
\newcommand*{\vecR}{\vec{R}}

\newcommand*{\matP}{\mat{P}}
\newcommand*{\matN}{\mat{N}}

\newcommand*{\matL}{\mat{L}}
\newcommand*{\matLbar}{\mat{\overbar{ L}}}

\newcommand*{\vecBphi}{\vec{B}_{\varphi}}
\newcommand*{\vecBmu}{\vec{B}_{\mu}}

\newcommand*{\vecC}{\vec{C}}
\newcommand*{\vecLambda}{\vec{\Lambda}}

\newcommand*{\phat}{\hat{\varphi}}

\newcommand*{\omr}{\omega_{r}}

\newcommand*{\sumR}{\sum_{r=1}^{R}}

\newcommand*{\aii}{a_{ii}}

\newcommand*{\vphat}[1]{\hat{\varphi}_{#1}}

\newcommand*{\mhat}[1]{\hat{\mu}_{#1}}

\newcommand*{\nkn}{\vNormal_{kn}}
\newcommand*{\refEn}{\Eref_{n}}

\newcommand*{\sintUnit}{\int_{0}^{1}}

\newcommand*{\ghat}{\boldsymbol{\hat{\gamma}}}
\newcommand*{\ghatn}{\ghat_{n}}
\newcommand*{\ghatns}{\ghatn(s)}

\newcommand*{\Fk}{\boldsymbol{F}_{k}}

\newcommand*{\betaMap}{\hat{\beta}_{kn}}

\newcommand*{\OmegaIn}{\Omega_{\mathrm{in}}}
\newcommand*{\OmegaOut}{\Omega_{\mathrm{out}}}

\newcommand*{\refIntPtTri}{\vec{\hat{q}}_{r}}

\newcommand*{\permutMat}{\vecc{\Delta}_{n}}

\newcommand*{\um}{u^m}
\newcommand*{\uOne}{u^{1}}
\newcommand*{\uTwo}{u^{2}}

\newcommand*{\setT}{\ensuremath{\mathcal{T}}}
\newcommand*{\setE}{\ensuremath{\mathcal{E}}}

\newcommand*{\vphi}{\varphi}                                     
\newcommand*{\mapEE}{\hat{\vec{\vartheta}}}                      
\renewcommand*{\vec}[1]{{\boldsymbol{#1}}}                       
\DeclareMathAlphabet{\mathbfsf}{\encodingdefault}{\sfdefault}{bx}{n}
\newcommand*{\vecc}[1]{\mathbfsf{#1}}                            
\newcommand*{\grad}{\vec{\nabla}}                                
\newcommand*{\dd}{\mathrm{d}}                                    
\newcommand*{\abs}[1]{\ensuremath{|#1|}}                         
\newcommand*{\transpose}[1]{{#1}^\mathrm{T}}                     
\newcommand*{\invtrans}[1]{{#1}^\mathrm{-T}}                     
\DeclareMathOperator*{\diag}{diag}                               
\DeclareMathOperator*{\diam}{diam}                               

\newcommand*{\IP}{\ensuremath\mathds{P}}
\newcommand*{\IR}{\ensuremath\mathds{R}}
\newcommand*{\IN}{\ensuremath\mathds{N}}

\newcommand*{\Matlab}{\mbox{MATLAB}}
\newcommand*{\Octave}{GNU~\mbox{Octave}}
\newcommand*{\MatOct}{\Matlab\,/\,\allowbreak\Octave} 

\newcommand*{\setEh}{\setE_h}

\newcommand*{\mat}[1]{\vecc{#1}}
\newcommand{\I}{I} 
\newcommand{\II}{I\!I} 
\newcommand{\III}{I\!I\!I} 
\newcommand{\IV}{I\!V} 
\newcommand{\V}{V} 
\newcommand{\VI}{V\!I} 
\newcommand{\VII}{V\!I\!I} 
\newcommand{\VIII}{V\!I\!I\!I} 
\newcommand{\IX}{I\!X} 
\newcommand{\X}{X} 
\newcommand{\XI}{X\!I}

\newcommand{\source}{\xi}
\newcommand*{\vecH}{\vec{H}}

\newcommand{\dimT}{d_{\setT_h}}
\newcommand{\dimE}{d_{\setE_h}}

%% file: tableStyles.tex
\pgfplotsset{table/steadyTable/.style=
	{	
		every head row/.style={
			before row={\toprule}, 
			after row={\midrule} 
		},
		every last row/.style={after row=\bottomrule}, 
		font=\footnotesize, 
		column type = {@{}lr@{\;\;}rr@{\;\;}rr@{\;\;}rr@{\;\;}rr@{\;\;}r@{}},
		create on use/orderP0/.style={create col/dyadic refinement rate={errP0}},
		create on use/orderP1/.style={create col/dyadic refinement rate={errP1}},
		create on use/orderP2/.style={create col/dyadic refinement rate={errP2}},
		create on use/orderP3/.style={create col/dyadic refinement rate={errP3}},
		create on use/orderP4/.style={create col/dyadic refinement rate={errP4}},
		columns={lvl,errP0,orderP0,errP1,orderP1,errP2,orderP2,errP3,orderP3,errP4,orderP4},
		columns/Nel/.style={
			column name={$K$}, /pgf/number format/.cd, sci,precision=2,sci e, sci zerofill=true,
		},		
		columns/lvl/.style={
			column name={$j$}
		},
		columns/errP0/.style={column name={$p=0$},/pgf/number format/.cd,sci,precision=3,sci e, sci zerofill=true},
		columns/orderP0/.style={column name={$\alpha$},/pgf/number format/.cd,fixed,precision=2, zerofill=true},
		columns/errP1/.style={column name={$p=1$},/pgf/number format/.cd,sci,precision=3,sci e, sci zerofill=true},
		columns/orderP1/.style={column name={$\alpha$},/pgf/number format/.cd,fixed,precision=2, zerofill=true},
		columns/errP2/.style={column name={$p=2$},/pgf/number format/.cd,sci,precision=3,sci e, sci zerofill=true},
		columns/orderP2/.style={column name={$\alpha$},/pgf/number format/.cd,fixed,precision=2, zerofill=true},
		columns/errP3/.style={column name={$p=3$},/pgf/number format/.cd,sci,precision=3,sci e, sci zerofill=true},
		columns/orderP3/.style={column name={$\alpha$},/pgf/number format/.cd,fixed,precision=2, zerofill=true},
		columns/errP4/.style={column name={$p=4$},/pgf/number format/.cd,sci,precision=3,sci e, sci zerofill=true},
		columns/orderP4/.style={column name={$\alpha$},/pgf/number format/.cd,fixed,precision=2, zerofill=true},
	}
}

\pgfplotsset{table/unsteadyTable/.style=
{
	every head row/.style={
		before row={\toprule}, 
		after row={\midrule} 
	},
	every last row/.style={after row=\bottomrule}, 
	font=\footnotesize, 
	column type = {@{}lr@{\;\;}rr@{\;\;}rr@{\;\;}rr@{\;\;}r@{}},
	create on use/orderP0/.style={create col/dyadic refinement rate={errP0}},
	create on use/orderP1/.style={create col/dyadic refinement rate={errP1}},
	create on use/orderP2/.style={create col/dyadic refinement rate={errP2}},
	create on use/orderP3/.style={create col/dyadic refinement rate={errP3}},
	columns={Nel,errP0,orderP0,errP1,orderP1,errP2,orderP2,errP3,orderP3},
	columns/Nel/.style={
		column name={$K$}, /pgf/number format/.cd, sci,precision=2,sci e, sci zerofill=true,
	},
	columns/errP0/.style={column name={$p=0$},/pgf/number format/.cd,sci,precision=3,sci e, sci zerofill=true},
	columns/orderP0/.style={column name={order},/pgf/number format/.cd,fixed,precision=2, zerofill=true},
	columns/errP1/.style={column name={$p=1$},/pgf/number format/.cd,sci,precision=3,sci e, sci zerofill=true},
	columns/orderP1/.style={column name={order},/pgf/number format/.cd,fixed,precision=2, zerofill=true},
	columns/errP2/.style={column name={$p=2$},/pgf/number format/.cd,sci,precision=3,sci e, sci zerofill=true},
	columns/orderP2/.style={column name={order},/pgf/number format/.cd,fixed,precision=2, zerofill=true},
	columns/errP3/.style={column name={$p=3$},/pgf/number format/.cd,sci,precision=3,sci e, sci zerofill=true},
	columns/orderP3/.style={column name={order},/pgf/number format/.cd,fixed,precision=2, zerofill=true}, }
}

%% file: sections/introduction.tex

The discontinuous Galerkin (DG) method first introduced in the early 70s in \cite{ReedHill1973} went on to have an~illustrious career as one of the most popular numerical methods especially (but not exclusively) for fluid simulation and engendered a~whole family of numerical schemes (see, e.g. \cite{ABCM,ShuOverviewDG} and the references therein). The reasons for this success are many  \cite{BaBo2011,Wang201553}: an~extremely flexible framework easily lending itself to many different types of equations, stability and conservation properties comparable to those of the finite volume method, natural support for high order discretizations and various types of adaptivity ($h$-, $p$-, $r$-) as well as for complex domain geometries, etc. Due to a~favorable computation-to-communication ratio, the method also fits extremely well into the popular parallel and hybrid computational paradigms~\cite{ortwein201401}.

However, one aspect of the DG methodology places it at a~clear disadvantage compared to the classical finite element and finite volume methods: a~large number of degrees of freedom with corresponding memory requirements. This drawback becomes even more restrictive in time-implicit or stationary numerical solvers that rely on matrix assembly, where increases in the lengths of vectors of unknowns become quadratically compounded in the sizes of corresponding matrix blocks. One idea proposed to speed up the solution of linear systems resulting from DG discretizations exploits the specifics of DG approximation spaces. Since DG basis functions are usually element-local (and often also hierarchical), the linear system can be trivially split into parts corresponding to different polynomial orders resulting in a scheme somewhat inspired by multigrid solvers, where different polynomial approximations on a~fixed mesh play the role of fine and coarse mesh solutions of the classical multigrid. This approach can be carried out for each polynomial degree as in the $p$-multigrid method \cite{Fidkowski2005, Luo2006, Bassi2009} or using a~two scale technique as in the hierarchical scale separation (HSS) method or variations thereof \cite{AizingerKK2015,JaustSA2016,SchuetzAizinger2017,Thiele2017}.
Another common way to deal with this issue---and in many cases even to speed up the time-to-solution \cite{CGorHDG,Yakovlev2016,WBMS13}---is to use \emph{hybridization}. 
Roughly speaking, this means that the discretized PDE is augmented with an unknown---let us call it $\lambda_h$---supported on the skeleton of the mesh consisting of element edges in 2D and element faces in 3D.
Using static condensation on the algebraic system level produces a~significantly smaller system than that obtained for an~unhybridized DG method at the price of additional cell-wise (small and uncoupled) linear systems that have to be solved alongside the globally coupled system on the mesh skeleton. Since all local solves are element-local and fully decoupled, the parallel communication cost of this algorithm part is zero.

The idea of using hybridization can be traced back to the 60s \cite{Veubeke}, it has subsequently been used in the context of mixed methods \cite{AB85,BDM85}. 
In those works, $\lambda_h$ was not only considered an~implementation feature but also as a~way to obtain a~more accurate solution via postprocessing. 
Based on the work in \cite{CoGo04}, Cockburn and coworkers introduced the hybridized discontinuous Galerkin method in a unifying framework in \cite{COGOLA}. 
Subsequently, the method has been extended to various types of equations such as the Stokes \cite{CockburnHDGStokes,EgWal13} and Darcy-Stokes equations \cite{Egger2013}, the incompressible and compressible Navier-Stokes equations \cite{NguyenNavierStokes,NgPe12,NgPeCo11,SchMa11}, the Maxwell equations \cite{NguyenMaxwell}, and---particularly relevant for this work---the convection(-diffusion) equation \cite{NPC09,NPC09L,EgSch09}. 
For an~interesting unification framework, we would also like to draw reader's attention to a~recent publication~\cite{BuiThanh2015}. 

The current work applies and extends to the hybridized schemes the framework and design principles introduced in~\cite{FrankRAK2015,ReuterAWFK2016} for unhybridized DG formulations and epitomized in our \MatOct~toolbox \textsl{FESTUNG} (\textsl{F}inite \textsl{E}lement \textsl{S}imulation \textsl{T}oolbox for \textsl{UN}structured \textsl{G}rids) available at~\cite{FESTUNG,FESTUNGGithub}. 
Citing from \cite{FrankRAK2015}, we aim to
\begin{enumerate}
\item Design a~general-purpose software package using the DG~method for a~range of standard applications 
and provide this toolbox as a~research and learning tool in the open source format.
\item Supply a~well-documented, intuitive user-interface to 
ease adoption by a~wider community of application and engineering professionals.
\item Relying on the vectorization capabilities of \MatOct, optimize the computational performance 
of the toolbox components and demonstrate these software development strategies.
\item Maintain throughout full compatibility with \Octave~to support users of open source software.
\end{enumerate}
For details about basic data structures and a~general overview of the solver structure, we refer the interested reader to our first publication~\cite{FrankRAK2015}, which applies the local discontinuous Galerkin (LDG) method to the diffusion operator.
A~DG discretization of the same model problem as in this work combined with higher-order explicit time stepping schemes and arbitrary order vertex-based slope limiters \cite{Kuzmin2010, Aizinger2011} is presented in the second paper in series~\cite{ReuterAWFK2016}.
The implementation presented in the current publication makes use of a~new solver structure tailored towards improved readability and maintainability of the code and designed to ease coupling of different solvers.
A~detailed description of this new structure with an~outline of the coupling capabilities is in preparation~\cite{ReuterRAK2017}.

The rest of the paper is structured as follows:
The model problem is introduced in the next section accompanied by a~detailed description of the space and time discretization in Section~\ref{sec:discretization}. 
Important aspects of the implementation including local mappings, numerical quadrature, and the assembly of the system matrix are presented in Section~\ref{sec:implementation}.
In Section~\ref{sec:results}, the code is verified using analytical convergence tests, and the numerical results are compared to those of the unhybridized DG implementation of the model problem from our previous publication~\cite{ReuterAWFK2016}.
Section~\ref{sec:registerroutines} lists the routines mentioned in this article, and Section~\ref{sec:conclusion} contains some conclusions and a~brief outlook of future tasks.

%% file: sections/model.tex
%
Let $J\coloneqq (0,t_\mathrm{end})$ be a~finite time interval, and let $\Omega \subset \IR^2$ be a~polygonally bounded Lipschitz domain with boundary~$\partial\Omega$.
We solve the \emph{linear advection equation}
\begin{equation}
  \parT \, c\,(t, \vx) + \divergence \flux{t,\vx,\ctx} = \source\,(t,\vx) \quad \text{in} \; J\times\Omega\,,
  \label{eq:model}
\end{equation}
where the scalar quantity~$c:J\times\Omega \rightarrow \IR$ is unknown. The flux function $\fluxNoDep{}:J\times\Omega\times\IR \rightarrow \IR^2$ determines the type of the problem and may depend on time $t$ and space coordinate $\vx$.
Within the context of this work, we assume
\begin{equation}\label{eq:flux}
  \flux{t,\vx,\ctx} \coloneqq \transpose{\begin{bmatrix} \uOne(t,\vx), & \uTwo(t,\vx) \end{bmatrix}} \, \ctx\,
\end{equation}
with given velocity~$\vu \coloneqq \transpose{[\uOne, \uTwo]}$ and source function~$\source:J\times\Omega\rightarrow\IR$ independent of~$c$.
Additionally, Eq.~\eqref{eq:model} is equipped with some initial 
\begin{equation}\label{eq:initial}
c   \;=\; c^0 \quad\text{on}~\{0\}\times\Omega 
\end{equation}\label{eq:bc}
and Dirichlet boundary condition
\begin{equation}
c  \;=\; c_\mathrm{D} \quad \text{on}~J\times{\partial\Omega}_\mathrm{in}(t)
\end{equation}
on inflow boundaries~$\partial\OmegaIn(t)\coloneqq \{ \vec{x} \in \partial\Omega \,|\, \vu(t, \vx) \cdot \vNormal(\vx) < 0 \}$ with $\vNormal(\vx)$ denoting the outward pointing unit normal vector. 
Since \eqref{eq:model} specifies a~first-order hyperbolic equation, no boundary conditions need to be prescribed on the outflow boundary~$\partial\OmegaOut(t)\coloneqq\partial\Omega\setminus\OmegaIn(t)$. 
Functions~$c^0:\Omega \mapsto \IR$ and~$c_\mathrm{D}:J\times\partial\OmegaIn(t)\mapsto\IR$ are given. 
The model problem is the same as the one considered in the second paper of our series~\cite{ReuterAWFK2016}; however, a~slightly different notation is used here.

%% file: sections/discretization.tex

\newcommand{\cjs}[1]{{{#1}}}
\newcommand{\va}[1]{{\color{red}{#1}}}

\subsection{Notation}
\label{sec:notation}

Let $\setT_h=\{T\}$ be a~regular family of non-overlapping partitions of a~polygonally bounded domain~$\Omega$ into~$K$ closed triangles~$T$.
The hybridized discontinuous Galerkin scheme employed in this work uses unknowns located on element boundaries, i.e., on the edges of elements (also referred to as the trace of the mesh). 
We introduce the set of all edges~$\setEh = \inE\cup\bcE$ consisting of sets of interior edges~$\inE$ and boundary edges~$\bcE$.
The latter is split into inflow~$\bcEin$ and outflow edges~$\bcEout$.
Furthermore, on each edge~$E_{\kEdge}\in\setEh$, a~unit normal~$\vNormal_{\kEdge}$ is defined such that~$\vNormal_{\kEdge}$ is exterior to element~$T_{\kEdge^-}${, where $T_{\kEdge^-}$ is an~arbitrarily chosen but fixed element adjoining $E_{\kEdge}$}. 
Obviously, boundary edges have only one adjacent element. For interior edges, the element opposite to $T_{\kEdge^-}$ is called~$T_{\kEdge^+}$. 
\cjs{The minus sign in $T_{\kEdge^-}$ is suppressed when no confusion is possible.}
The total number of edges is denoted by 
{
\[\Kedge \coloneqq |\setEh|, \qquad \mbox{and} \qquad \Gamma = \bigcup_{\kEdge=1}^{\Kedge} \Ekbar 
\]}
is the trace of the mesh.
We refer to elements using \cjs{indices~$k$, and} add bars~$\bar\cdot$ (e.g., $\kEdge$ or $\Kedge$) whenever we refer to edges to distinguish one from another.
Additionally, we make use of the notation from the previous publications~\cite{FrankRAK2015,ReuterAWFK2016} where we referred to edges of an element~$T_k$ as~$E_{kn},n\in\{1,2,3\}$.
Note the difference between this element-local edge numbering (\enquote{$E_{kn}$ is the $n$th edge of the $k$th element}) and global edge numbering (\enquote{$E_{\kEdge}$ is the $\kEdge$th edge in the mesh}).

\begin{figure}[t!]\centering
\includegraphics{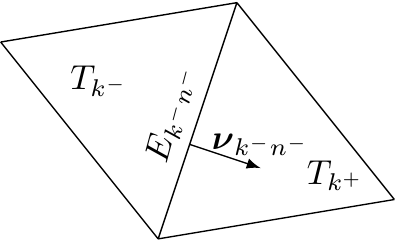}
\captionsetup{justification=raggedright,singlelinecheck=false}
\caption{Two triangles adjacent to edge~${E_{\kEdge}}$. It holds~${E_{\kEdge} = }E_{k^-n^-} = E_{k^+n^+}$ and $\vec{\nu}_{k^-n^-} = -\vec{\nu}_{k^+n^+}$.}
\label{fig:T1T2}
\end{figure}

For the description of the method, we need mappings allowing us to switch back and forth between element-local and global indices.
For that reason, we introduce a~mapping~$\rho(k,n)$ that relates the $n$th edge~$E_{kn}$ of element~$T_k$ to its global index in the set of edges~$\setEh$
\begin{equation}
  \rho: \{1, \ldots, K\} \times \{1,2,3\} {\rightarrow \{1, \ldots, \Kedge\}, \quad }(k,n) \mapsto \kEdge \,.
\label{eq:mapkntokEdge}
\end{equation}
This mapping is not injective since for each interior edge there exists a~pair of index tuples $(k^{-},n^{-})$ and $(k^{+},n^{+})$ that map to the same edge index $\kEdge$, i.\,e.,~$E_{\kEdge} = E_{k^-n^-} = E_{k^+n^+}$ (cf. Fig.~\ref{fig:T1T2}).
We define a~second mapping~$\kappa(\kEdge,l)$ to identify elements $T_{k^-}, T_{k^+}$ adjacent to an edge~$E_{\kEdge}$.
In this, argument $l\in\{1,2\}$ denotes the edge-local index of the adjacent elements; it is constructed so that~$l=1$ refers to the \enquote{inner} element~$T_{k^{-}}$ (that always exists) while~$l=2$ refers to the \enquote{outer} element $T_{\kEdge^{+}}$ that exists only for interior edges~$\setE_{\mathrm{int}}$,
\begin{equation}
  \kappa: \{1,\ldots,\Kedge\} \times \{1,2\} {\rightarrow \{0, \ldots, K\}, \quad }(\kEdge,l) \mapsto \left\{ \begin{matrix}
    k^- \in \{1,\ldots,K\}\,, & \text{if } l=1 \\
    k^+ \in \{0,\ldots,K\}\,, & \text{if } l=2 \end{matrix} \right\} \,.
  \label{eq:mapKappa}
\end{equation}
As element indices start counting from 1, we set $\kappa(\kEdge, 2) = 0 \,, \forall \Ekbar\in\bcE$ to mark the absence of the \enquote{outer} element.

\subsection{Semi-discrete form and hybridization}
\label{section:semidiscrete}
{The particular feature of the HDG method that distinguishes it from \cjs{conventional DG methods} is that solutions are not only sought on elements, but also on the skeleton of the mesh. To this end, one has to define approximation spaces on both $\setT_h$ and $\setE_h$.}
Given~$\IP_p(T)$ and~$\IP_p(E)$, the spaces of complete polynomials of degree at most~$p$ on an element~$T\in\setT_h$ \cjs{and} edge~$E\in\setE_h$, {respectively, } we denote \cjs{global} broken polynomial spaces~{$\IP_p(\setT_h), \IP_p(\setE_h)$} on the triangulation~$\setT_h$ and its set of edges~$\setE_h$ by
\begin{equation}\label{eq:spaces}
  {\IP_p(\setT_h)} \coloneqq \lbrace \varphih : \bar{\Omega} \rightarrow \R: \; \forall T \in \mathcal{T}_h, \; \varphih|_{T} \in \mathbb{P}_p(T) \rbrace\,, \qquad
  {\IP_p(\setE_h)} \coloneqq \lbrace \muh : \Gamma \rightarrow \R: \; \forall E \in \setEh, \; \muh|_{E} \in \mathbb{P}_p(E) \rbrace\,.
\end{equation}
{%
To support the tutorial character of this paper, we briefly discuss the derivation of the HDG method for the advection equation and refer the reader to~\cite{NPC09L} for more details. 
The standard DG discretization of \eqref{eq:model} is obtained by multiplying \eqref{eq:model} with a~test function~$\varphi_h \in \IP_p(\setT_h)$, applying integration by parts on element $\Tk \in \mathcal{T}_h$, and replacing the flux~$\fluxNoDep$ (see \eqref{eq:flux}) on element boundaries~$\partial\Tk$ by a numerical flux function $\numFluxDG(\cmh, \cph)$ well-known from finite volume discretizations~\cite{GR1,GR2,Kroener} \va{that depends on the values of the solution $\ch$ from both sides} of the boundary of element $\Tk \in \setT_h$. Omitting the time argument of $c_h$ in the flux function for the sake of readability, this yields
\begin{align}\label{eq:semi-discrete-dg}
  \intTk \parT c_h(t) \, \varphi_{h} \, \dx 
  -
  \intTk \flux{\ch} \cdot \nabla \phih \, \dx
  +
  \sintTk \numFluxDG(\cmh,\cph) \cdot \vNormal_{\Tk} \,\va{\phih} \,\ds
  &\,=\,    \intTk \source\cjs{(t)} \, \varphi_{h} \, \dx \,\quad&& \forall \phih \in \IP_p(\setT_h)\,.
\end{align}
At the domain boundary $\partial\Omega$, $\cph$ has of course to be replaced by a suitably chosen boundary value $c_{\partial \Omega}(\cmh)$ as in
\begin{align}
 \label{eq:outflow_operator}
 c_{\partial \Omega}(\cmh) = \begin{cases} \cmh, \quad \text{on outflow boundary ${\partial\Omega}_\mathrm{out}$}\,, \\ 
c_D, \quad \text{on inflow boundary ${\partial\Omega}_\mathrm{in}$}\,. \end{cases}
\end{align}
In this work, the numerical flux is chosen to be a~slightly modified local Lax-Friedrichs/Rusanov flux~\cite{LF,Rusanov}
\begin{equation}
\numFluxDG(\cmh,\cph) \; \coloneqq \;
\flux{\frac{\cmh+\cph}{2}} \; + \; \frac\alpha 2 (\cmh - \cph ) \,\vNormal_{\Tk}\,; 
\label{eq:numfluxDG}
\end{equation}
however, other choices are possible \cite{BuiThanh2015}. The stabilization parameter~$\alpha$ must satisfy $\alpha \ge \left\lbrace \max\{|\flux{\ch^-} \cdot \vNormal_{\cjs{\Tk}} |,  |\flux{\ch^+} \cdot \vNormal_{\cjs{\Tk}} | \} \right\rbrace$. In this work $\alpha$ is a user-defined constant that is the same for every edge. At the boundary, we simply replace the flux by $\flux{c_{\partial \Omega}(\cmh)}$\cjs{, also see \eqref{eq:outflow_operator}.}
\cjs{\eqref{eq:numfluxDG} can be formally rewritten to depend on the mean value of $\cmh$ and $\cph$ and either on $\cmh$ or on $\cph$ instead of $\cmh$ and $\cph$.}
Upon defining a quantity $\lambdah(t)$ on an edge $E_{\bar k}$ as 
\begin{align}
 \label{eq:lambdastrong}
 \lambdah(t) := \begin{cases} \frac{1}{2}(\cmh(t) + \cph(t))\,, &\qquad E_{\bar k} \in \inE \,,
\\ c_{\partial \Omega}(\cmh(t))\,, &\qquad E_{\bar k} \in \bcE\,, \end{cases}
\end{align}
the flux $\numFluxDG$ in \eqref{eq:numfluxDG} can also be written as
\begin{equation}
\numFlux(\lambdah, \cmh) \coloneqq
\begin{cases} \flux{\lambdah} - \alpha (\lambdah - \cmh) \,\vNormal_{\Tk}\,, &\qquad E_{\bar k} \in \inE \,,\\
\flux{\lambdah}\,, &\qquad E_{\bar k} \in \bcE \,.\end{cases}
\label{eq:numfluxDG2}
\end{equation}
Note that $\lambdah$ on an interior edge is a polynomial of degree $p$. 
Discretizing the strong formulation of $\lambdah$ given in ~\eqref{eq:lambdastrong} in a~DG framework is straightforward. Together with the slightly modified (only the flux function differs!) version of \eqref{eq:semi-discrete-dg} this yields the HDG discretization: \cjs{Seek $(c_h(t),\lambdah(t))\in \IP_p(\setT_h)\times \IP_p(\setE_h)$ such that for all $t \in J$ there holds}
\begin{subequations}\label{eq:semi-discrete}
\begin{align}\label{eq:semi-discrete1}
  &\intTk \parT c_h(t) \, \varphi_{h} \, \dx 
  -
  \intTk \flux{\ch} \cdot \nabla \phih \, \dx
  +
  \sintTk \numFlux(\lambdah, \cmh) \cdot \vNormal_{\Tk} \,\va{\phih} \,\ds
  \,=\,    \intTk \source \, \varphi_{h} \, \dx \,,
  \quad&& \forall \phih \in \IP_p(\setT_h)\,,
  \\
\label{eq:semi-discrete2}
  &\cjs{\int_{\Ekbar} \mh \left\{\begin{array}{cl}
    \alpha \left( 2 \lambdah(t) - \cmh(t) - \cph(t) \right)  & \text{on } \Ekbar\in\inE \\
    \lambdah(t) - c_{\partial \Omega}(\cmh(t)) & \text{on } \Ekbar\in\bcE
  \end{array}\right\} \, \ds  \,=\, 0\,,}
  \quad&& \forall \mh \in \IP_p(\setE_h)\,.
\end{align}
\end{subequations}
Note that $\alpha$ in \eqref{eq:semi-discrete2} \cjs{has been introduced} to make the formulation compatible to the one used for viscous problems \cite{NPC09L}. 

It can immediately be seen that the above formulation is nothing else than our \cjs{point-of-departure DG method, see \eqref{eq:semi-discrete-dg}, written in another, yet equivalent, way.}
}%
At first glance, this scheme does not seem to have any advantages compared to other DG schemes: on the contrary, an~additional equation and an~additional unknown~$\lambdah$ are apparent. 
However, the choice of the numerical flux on element boundaries leads to an~inter-element coupling solely by the function~$\lambdah$.
In Sec.~\ref{sec:staticcondensation}, we explain in detail how this structure of the system can be exploited to reduce the number of globally coupled unknowns, which turns out to be especially attractive for discretizations that rely on implicit solution techniques such as time-implicit schemes or stationary problems. 
{
To give a first glimpse of this procedure consider the stationary part of \eqref{eq:semi-discrete} (i.e. $\parT c_h = 0$). This results in a~linear system of equations of form 
\begin{align*}
 \left( \begin{matrix} \matLbar & \matMbar \\ \vecc N & \vecc P \end{matrix} \right) \left( \begin{matrix} \vecc C \\ \vecc \Lambda \end{matrix} \right) = \left( \begin{matrix} \vecBphi \\ \vecBmu \end{matrix} \right)
\end{align*}
with matrices \cjs{$\matLbar \in \R^{\dimT\times\dimT}$, $\matMbar\in \R^{\dimT\times\dimE}$, $\vecc N\in \R^{\dimE\times\dimT}$ and $\vecc P\in\R^{\dimE\times\dimE}$, where $\dimT = \dim(\IP_h(\setT_h))$ and $\dimE = \dim(\IP_h(\setE_h))$.} 
$\vecc C$ and $\vecc \Lambda$ are basis coefficients corresponding to $\ch$ and $\lambdah$, respectively, see \eqref{eq:basis_coefficients}. It is important to note that due to the inter-element coupling in $\lambda_h$ only, the matrix $\matLbar$ (possibly reordered) is block-diagonal and can therefore be easily inverted. Thus, $\vecc C$ can be computed by 
\begin{align*}
 \vecc C = \matLbar^{-1} \left(\vecBphi - \matMbar \vecc \Lambda\right), 
\end{align*}
which can be plugged into the second equation to obtain 
\begin{align*}
 \vecc N \matLbar^{-1} \left(\vecBphi - \matMbar \vecc \Lambda\right) + \vecc P \vecc\Lambda = \vecBmu.
\end{align*}
This is an algebraic system in $\vecc\Lambda$ only, that for higher order polynomial spaces tends to be much smaller than the corresponding linear system for an~unhybridized DG scheme of the same order. This has a tremendous influence on the efficiency of iterative solvers and is one of the reasons for the success of the HDG method. 

In the following sections, we make these introductory remarks more specific and explain their implementation in the software. 
}

\subsubsection{Local basis representation}
\label{sec:localbasis}

The DG function spaces defined in~\eqref{eq:spaces} do not have any continuity constraints across element or edge boundaries.
Consequently, a~two-dimensional basis function~$\varphi_{ki}: \bar{\Omega} \rightarrow \IR$ for~${\IP_p(\setT_h)}$ is only supported on~$T_k \in \setT_h$ and must fulfill
\begin{equation*}
\forall k \in \{1,\ldots,K\}\,, \qquad 
\IP_p(T_k) = \mathrm{span} \big\{ \varphi_{ki} \big\}_{i\in\{1,\ldots,N\}}  \qquad
\mathrm{ with } \quad
N \coloneqq \frac{(p+1)(p+2)}{2} \,.
\end{equation*}
In the same way, a~one-dimensional basis function~$\mki: \Gamma \rightarrow\IR$ for~${\IP_p(\setE_h)}$ is only supported on the edge~$\Ekbar \subset \Gamma$, ensuring
\begin{equation*}
\forall \kEdge \in \{ 1,\ldots,\Kedge \}\,, \qquad 
\IP_p(\Ekbar) = \mathrm{span} \big\{ \mki \big\}_{i\in\{1,\ldots,\Nedge\}} \qquad
\mathrm{ with } \quad
\Nedge \coloneqq p + 1\,.
\end{equation*}
We denote the numbers of local degrees of freedom on a~triangle or an~edge by~$N$ or~$\Nedge$, respectively.
In this work, we choose orthonormal, hierarchical Legendre polynomials as basis functions. 
For details and closed-form expressions of the two-dimensional basis functions~$\varphi_{ki}$ we refer to our first paper~\cite{FrankRAK2015} {(those can be easily obtained from an~arbitrary polynomila basis via e.g. Gram-Schmidt orthogonalization procedure using a~symbolic algebra software)} .
The one-dimensional basis functions up to order four given on the unit interval~$[0,1]$ are defined as
\begin{equation}
\begin{aligned}
&\oneDBasisHat{1} \coloneqq 1 \,, \qquad
\oneDBasisHat{2} \coloneqq \sqrt{3} (1 - 2s) \,,\qquad
\oneDBasisHat{3} \coloneqq \sqrt{5} (1 - 6s + 6s^2) \,,\\
&\oneDBasisHat{4} \coloneqq \sqrt{7} (-1 + 12s -30s^2 + 20s^3) \,,\qquad
\oneDBasisHat{5} \coloneqq \sqrt{9} (1 - 20s + 90s^2 -140s^3 + 70s^4) \,.
\end{aligned}
\label{eq:phi1D}
\end{equation}
The local discrete solutions $\ch$ on $T_k\in\setT_h$ and $\lambdah$ on $\Ekbar\in\setE_h$ are represented using local bases on elements and edges
\begin{equation}
  \label{eq:basis_coefficients}
  \ch(t,\vx)\big|_{T_k} = \sum_{j=1}^{N} \Ckj(t) \, \varphi_{kj}(\vx) \,,\qquad
  \lambdah(t, \vx)\big|_{\Ekbar}  = \sum_{j=1}^{\Nedge} \Lambda_{\kEdge j}(t)\, \mujFaceView(\vx)\,.
\end{equation}
%

\subsubsection{System of equations}
Testing~\eqref{eq:semi-discrete1} with $\varphi_h = \pki, i=1, \ldots, N$ yields a~time-dependent system of equations with the contribution from~$T_k$ given by
\begin{subequations}\label{eq:soe}
\begin{equation}
  \begin{aligned}
    &
    \underbrace{\parT \sumElCkjt \intTk \pkj\,\pki \,\dx}_{\I~(\mat\matMphi)}
    - 
    \underbrace{\sumElCkjt \summ \intTk \um(t, \vx)\, \pkj \,\parXm \pki\, \dx}_{\II~(-\matG^1~\text{and}~-\matG^2)}
    +\underbrace{\sumElLkj{(t)} \sintTkNB \left(\vecu(t, \vx) \cdot \vNormal\right)\, \mkj\, \pki \,  \ds}_{\III~(\matS)}
    \\
    &\quad 
    - 
    \underbrace{\alpha \sumLkj{(t)} \sintTkNB \mkj\,\pki\, \ds}_{\IV~(\alpha \matRmu)}
    + 
    \underbrace{\alpha \sumElCkj{(t)} \sintTkNB \pkj\, \pki\, \ds}_{\V~(\alpha \matRphi)}
    \\
    &\quad  
    +\underbrace{
    \sintTkB \pki 
    \left(\vecu(t, \vx) \cdot \vNormal\right)
    \left\lbrace
    \begin{matrix}
      \sumLkj{(t)} \, \mkj   & \textOutflow \\
      c_\mathrm{D}(t,\vx) &  \textInflow
    \end{matrix}
    \right\rbrace
    \ds
    }_{\VI~(\matSout~\text{and}~\vecFphiD)}
    =
    \underbrace{\intTk \source\cjs{(t)} \, \pki\, \dx}_{\VII~(\vecH)}\, .
  \end{aligned}
  \label{eq:soe:eq1}
\end{equation}
In \eqref{eq:flux} $\fluxNoDep{}$ was defined as a~linear function. 
\cjs{Upon testing with $\mh=\mki$, t}he semi-discrete form of \eqref{eq:semi-discrete2} is then given by
\begin{equation}
  \begin{aligned}
  &
  \underbrace{\alpha \sumLkj{(t)} \sintTkNB \mkj \, \mki \, \ds}_{\VIII~(\alpha \matMmuBar)}
  - \underbrace{\alpha \sumElCkj{(t)} \sintTkNB \pkj \, \mki \, \ds}_{\IX~(- \alpha \matT)}
  \\
  &\quad+
  \underbrace{\sumLkj{(t)} \sintTkB \mkj \, \mki \, \ds}_{\X~(\matMmuTilde)}
  -
  \underbrace{
  \sintTkB \mki
  \left\lbrace
  \begin{matrix}
  \sumElCkj{(t)} \, \pkj & \textOutflow \\
  c_\mathrm{D}(t, \vx)  &  \textInflow
  \end{matrix}
  \right\rbrace
  \ds }_{\XI~(- \matKmuOut \text{ and } \vecKmuD)}
  =0
  \end{aligned}
  \label{eq:soe:eq2}
\end{equation}
\end{subequations}
We can rewrite system~\eqref{eq:soe} in matrix form
\begin{subequations}
  \label{eq:systemAsMatrices}
\begin{align}
  \matMphi \parT \vecC{(t)}
  +
  \underbrace{\left(- \matG^1{(t)} - \matG^2{(t)} + \alpha \matRphi \right)}_{\eqqcolon\matLbar(t)} \vecC{(t)}
  + 
  \underbrace{\left( \matS{(t)} + \matSout{(t)}
    - \alpha\matRmu \right)}_{\eqqcolon\matMbar(t)} \vecLambda{(t)}
  &= \underbrace{\vecH{(t)} - \vecFphiD{(t)}  }_{\eqqcolon\vec{B}_{\varphi}(t)} \,,
  \label{eq:systemAsMatrices:eq1}\\
  \underbrace{(- \alpha \matT - \matKmuOut{(t)})}_{\eqqcolon\vecc{N}{(t)}} \vecC{(t)} + \underbrace{\left(\alpha \matMmuBar + \matMmuTilde \right)}_{\eqqcolon\vecc{P}} \vecLambda{(t)} &= \underbrace{ - \vecKmuD{(t)} }_{\eqqcolon\vecBmu(t)}
  \label{eq:systemAsMatrices:eq2}
\end{align}
\end{subequations}
with representation vectors
\begin{equation*}
\vecC(t) \coloneqq \transpose{ \Big[ C_{11}(t) \cdots C_{1N}(t) \cdots\, \cdots C_{K1}(t) \cdots C_{KN}(t) \Big]} \,, \qquad
\vecLambda \coloneqq \transpose{ \Big[ \Lambda_{11}(t) \cdots \Lambda_{1\Nedge}(t) \cdots\, \cdots \Lambda_{\Kedge1}(t) \cdots \Lambda_{\Kedge\Nedge}(t) \Big]} \,.
\end{equation*}
Recall that the inter-element coupling of the solution~$\ch$ in system~\eqref{eq:systemAsMatrices} acts only via the edge function~$\lambdah$ via~$\matMbar(t) \in \IR^{KN\times \Kedge\Nedge}$ and~$\matN \in \IR^{\Kedge\Nedge\times KN}$ -- block matrices generally combining rectangular blocks into rectangular block structures. 
{The rectangular~\cjs{$N\times\Nedge$ (or transposed)} block shapes arise from the difference in the number of local degrees of freedom on elements~($N$) and edges~($\Nedge$) for $p>0$ (i.\,e., $N>\Nedge$), whereas the rectangular block-matrix structure is the result of the difference between the number of edges~($\Kedge$) and the number of elements~($K$) in the triangulation. \cjs{Thus,~$K\times\Kedge$ such rectangular~$N\times\Nedge$ blocks build up the global matrix, resulting in a~global $KN\times\Kedge\Nedge$ matrix~$\matMbar$ (or~$\matN$ with all dimensions transposed).}}

\subsubsection[Contributions from volume terms I, II, VII]{Contributions from volume terms~\I,~\II~and~\VII}
All matrices in system~\eqref{eq:systemAsMatrices} have sparse block structure, and we define the non-zero entries in the remainder of this section.
For all remaining entries we tacitly assume zero fill-in.

The integral of term~\I~gives the standard mass matrix $\matMphi \in \R^{KN \times KN}$ with components defined as
\begin{subequations}\label{eq:matMphi}
\begin{equation}
  [\matMphi]_{(k-1)N+i,(k-1)N+j} \coloneqq \intTk \pkj\,\pki\, \dx\,.
\end{equation}
This leads to a block diagonal matrix because basis and test functions~$\pki$, $i \in\{1,\ldots,N\}$ are only supported on element $\Tk$. 
Therefore 
\begin{equation}
  \matMphi =
  \begin{bmatrix}
  \matM_{\varphi,T_1} & & \\  & \ddots & \\ & & \matM_{\varphi,T_N}\\
  \end{bmatrix}
  \quad\text{with \emph{local mass matrix}}\quad
  \matM_{\varphi,\Tk}  \coloneqq 
  \intTk
  \begin{bmatrix}
    \varphi_{k1}\varphi_{k1} & \ldots & \varphi_{k1}\varphi_{kN} \\
    \vdots           & \ddots & \vdots \\
    \varphi_{kN}\varphi_{k1} & \ldots & \varphi_{kN}\varphi_{kN} \\
  \end{bmatrix} \dx \in \R^{N \times N}\,.
  \label{eq:massmatrix}
\end{equation}
\end{subequations}
We abbreviate~$\matMphi = \mathrm{diag}( \matM_{\varphi,T_1}, \ldots, \matM_{\varphi,T_N} )$.

The definition of block matrices~$\matGm \in \R^{KN \times KN}, \ m \in \{1,2\}$ differs slightly from the form in our previous publication~\cite{ReuterAWFK2016} because we do not use the projected advection velocity but rather evaluate the velocity function\,/\,flux at each quadrature point.
The component-wise entries are given by
\begin{subequations}\label{eq:matG}
\begin{equation}
  [\matGm]_{(k-1)N+i,(k-1)N+j} = \intTk \um \,\pkj\, \parXm \pki \,\dx
\end{equation}
again leading to a block-diagonal matrix~$\matGm = \mathrm{diag}(\matGm_{T_1}, \dots \matGm_{T_{K}})$, where each block is given by
\begin{equation}
  \matGm_{T_{K}} \coloneqq
  \intTk
  \um
  \begin{bmatrix}
    \varphi_{k1} \parXm\, \varphi_{k1} & \dots & \varphi_{kN}\, \parXm \varphi_{k1}  \\
    \vdots & \ddots & \vdots  \\
    \varphi_{k1} \,\parXm \varphi_{kN} & \ldots & \varphi_{kN}\, \parXm \varphi_{kN}  
  \end{bmatrix}
  \dx\,.
\end{equation}
\end{subequations}
The source term enters the discretization as an additional term $\vecH\in\R^{KN}$ on the right-hand-side of the equation. The entries are given as
\begin{subequations}\label{eq:vecH}
\begin{equation}
  [\vecH]_{(k-1)N+i} = \intTk \source\,\pki\,\dx
  \label{eq:vecHentry}
\end{equation}
such that the full vector is easily assembled as
\begin{equation}
\vecH = \begin{bmatrix} \vecH_{T_1} \\ \vdots \\ \vecH_{T_K} \end{bmatrix} 
\qquad\text{with}\quad
  \vecH_{\Tk} = 
  \intTk 
  \source
  \begin{bmatrix}
  \varphi_{k1} \\
  \vdots \\
  \varphi_{kN}
  \end{bmatrix} 
  \dx.
\end{equation}
\end{subequations}

\subsubsection[Contribution of edge terms III, IV, V, VI -- first equation]{Contribution of edge terms \III, \IV, \V, \VI --- first equation}
Compared to the DG discretization used in our previous works~\cite{FrankRAK2015,ReuterAWFK2016}, the number of edge integrals has significantly increased. 
This is caused by the following factors:
\begin{enumerate}
  \item Edge integrals are split into integrals over interior edges and over edges on the domain boundaries.
  \item The numerical flux function~\eqref{eq:numfluxDG2} contains three terms compared to only one for the upwind flux used in~\cite{ReuterAWFK2016}.
  \item An~additional unknown is introduced that is only defined on edges resulting in an~additional equation~\eqref{eq:semi-discrete2}.
\end{enumerate}
To improve readability, we split the presentation of edge integrals into two sections:
first, edge terms in the original equation for~$c_h$ and then the edge terms in the hybrid equation.
Throughout the assembly description and within the implementation we use the \emph{element-based view}, i.e., all edge terms are presented in a~form allowing for the assembly as nested loops over elements~$\Tk$, $k\in\{1,\dots,K\}$ and then edges~$\Ekn$, $n\in\{1,2,3\}$ of each element.
This is different from the \emph{edge based view} which would allow to assemble the terms in a~single loop over all edges~$\Ekbar$, $\kEdge\in\{1,\dots,\Kedge\}$.
We made this choice since the data structures in FESTUNG favor the element-based view.
For that reason, from now on, we always consider~$\mujElementView = \mkj$ using the mapping $\rho: (k,n) \mapsto \kEdge$ specified in~\eqref{eq:mapkntokEdge}.

Term~\III~contributes to matrix $\matS\in\IR^{KN\times\Kedge\Nedge}$ as 
\begin{subequations}\label{eq:matS}
\begin{equation}
  [\matS]_{(k-1)N+i,(\kEdge-1)\Nedge+j} =\sumInE  \sintEkn (\vecu\cdot \vNormal_{kn}) \,\mujElementView\, \pki\,  \ds\,.
  \label{eq:SbarEntry}
\end{equation}
The entries are structured into $N\times\Nedge$-blocks contributed by each edge on every element given as
\begin{equation}
\matS_{\Ekn} = \sintEkn (\vecu \cdot \vNormal_{kn})
\begin{bmatrix}
\muElementView{1} \, \varphi_{k1} & \ldots & \muElementView{\Nedge}  \, \varphi_{k1} \\
\vdots           & \ddots & \vdots \\
\muElementView{1}  \, \varphi_{kN} & \ldots & \muElementView{\Nedge}  \, \varphi_{kN} 
\end{bmatrix}
\ds\,.
\label{eq:Sbar}
\end{equation}
\end{subequations}
Note that this is the contribution of a single interior edge~$\Ekn \in \partial T_k \cap \setE_\mathrm{int}$ of triangle~$\Tk$.
As this integral is also evaluated on the neighboring element, this block will produce \emph{two} contributions to matrix $\matS$.

The matrix from term \IV{} is a~block matrix $\matRmu \in \R^{KN \times \Kedge \Nedge}$ with blocks of size $ N \times \Nedge$ given by
\begin{subequations}\label{eq:RbarMuDef}
\begin{equation}
  [\matRmu]_{(k-1)N+i,(\kEdge-1)\Nedge+j} = \sumInE \sintEkn  \mujElementView \, \pki \, \ds
  \label{eq:RmuEntry}
\end{equation}
with $i\in\{1, \ldots, N\}$ and $j\in\{ 1, \dots, \Nedge \}$.
The local matrix of a single edge reads 
\begin{equation}
  \matR_{\mu,\Ekn} = \sintEkn 
  \begin{bmatrix}
    \muElementView{1}\,\varphi_{k1} & \ldots & \muElementView{\Nedge}\,\varphi_{k1} \\
    \vdots           & \ddots & \vdots \\
    \muElementView{1}\,\varphi_{kN} & \ldots & \muElementView{\Nedge}\,\varphi_{kN} 
  \end{bmatrix}
  \ds\,.
  \label{eq:RbarMu}
\end{equation}
\end{subequations}
Term \V{} gives another block diagonal contribution $\matRphi \in \R^{KN \times KN}$, where each entry is given by 
\begin{subequations}\label{eq:matRphi}
\begin{equation}
  [\matRphi]_{(k-1)N+i,(k-1)N+j}  = \sum_{\Ekn\in\partial\Tk\cap\inE} \sintEkn \pkj \, \pki \, \ds.
  \label{eq:RphiEntry}
\end{equation}
The element-local matrix is then 
\begin{equation}
  \matR_{\varphi,\Tk} = 
  \sum_{\Ekn\in\partial\Tk\cap\inE}
  \sintEkn
  \begin{bmatrix}
    \varphi_{k1} \, \varphi_{k1} & \ldots & \varphi_{kN} \, \varphi_{k1} \\
                   \vdots & \ddots & \vdots \\
    \varphi_{k1} \, \varphi_{kN} & \ldots & \varphi_{kN} \, \varphi_{kN}
  \end{bmatrix}
  \ds\,,
\end{equation}
\end{subequations}
where each edge $\Ekn$ of each triangle $\Tk$ is visited exactly once, and~$\matRphi = \diag(\matR_{\vphi,T_1},\ldots,\matR_{\vphi,T_K})$.
This is slightly different from the previous work \cite{ReuterAWFK2016}, where the test functions $\varphi^{-}$ and $\varphi^{+}$ from two neighboring elements may have been multiplied because elements would be coupled directly with each other.

Term \VI{} incorporates the boundary conditions on boundary edges $\Ekbar \in \bcE$. In the case of an inflow boundary condition, this contributes to the right-hand-side. 
Each entry of vector~$\vecFphiD\in\IR^{KN}$ is given as
\begin{equation}
  [\vecFphiD]_{(k-1)N+i} = 
  \sum_{\Ekn\in\partial\Tk\cap\bcEin} \sintEkn (\vecu \cdot \vNormal_{kn}) \, c_\mathrm{D}\, \pki \, \ds\,.
  \label{eq:FphiDentry}
\end{equation}
Outflow boundary conditions depend on $\lambdah$ and therefore on the solution. 
This gives us an additional contribution~$\matSout\in\IR^{KN\times\Kedge\Nedge}$ to the left hand side, where each entry is given by
\begin{equation}
  [\matSout]_{(k-1)N+i,(\kEdge-1)\Nedge+j} 
  = \sum_{\Ekn\in\partial\Tk\cap\bcEout} \sintEkn (\vecu\cdot\vNormal_{kn}) \, \mujElementView \, \pki \,\ds\,.
  \label{eq:matSout}
\end{equation}
This is almost identical to~\eqref{eq:SbarEntry} with the only difference being the set of edges considered, and thus the sub-blocks take the same form as in equation~\eqref{eq:Sbar}.
In the implementation,~$\matS$ and~$\matSout$ are assembled together.

\subsubsection[Contribution of edge terms  VIII, IX, X, and XI -- hybrid equation]{Contribution of edge terms  \VIII{}, \IX{}, \X{}, and \XI{} --- hybrid equation}

The first term --- term \VIII{} --- is very similar to an~edge mass matrix with the only differences being \cjs{the} parameter $\alpha$ that has to be respected and the fact that every edge is visited twice because it is an integral over interior edges $\Ekbar \in \inE$.
$\matMmuBar\in\IR^{\Kedge\Nedge\times\Kedge\Nedge}$ is given by
%
\begin{equation}
\label{eq:matMmuBar1}
\lbrack\matMmuBar\rbrack_{(\kEdge-1) \Nedge+i,(\kEdge-1)\Nedge+j} \coloneqq
\sum_{\Ekn\in\partial \Tk \cap \inE} \int_{\Ekn} \mujElementView\,\muiElementView\, \ds
\end{equation}
for $i,j\in \{1,2,\ldots,\Nedge\}$ and $\kEdge \in \{1,2,\ldots,\Kedge\}$.
This leads to a block diagonal matrix because the ansatz and test functions $\muiElementView = \mki, i \in\{1,2,\ldots,\Nedge\}$ have support only on the corresponding edge $\Ekn = \Ekbar$.
Term \IX{} is very similar to term \IV, where each entry of the resulting matrix~$\matT\in\IR^{\Kedge\Nedge\times KN}$ is given as 
\begin{subequations}\label{eq:matT}
\begin{equation}
[\matT]_{(\kEdge-1)\Nedge+i,(k-1)N+j} = \sumInE \sintEkn \pkj \, \muiElementView \, \ds
\label{eq:matTentry}
\end{equation}
with $i\in \{ 1, \dots, \Nedge \}$, $j\in \{1, \ldots, N\}$, and $\kEdge$ given by the mapping in \eqref{eq:mapkntokEdge}.
The contribution of a single edge is
\begin{equation}
\matT_{\Ekn} =\sintEkn 
\begin{bmatrix}
\varphi_{k1} \, \muElementView{1} & \ldots & \varphi_{kN} \, \muElementView{1}  \\
\vdots           & \ddots & \vdots \\
\varphi_{k1} \, \muElementView{1}    & \ldots &  \varphi_{kN} \, \muElementView{\Nedge} 
\end{bmatrix}
\ds.
\label{eq:matU}
\end{equation}
\end{subequations}
In fact, we have $\matT = \matRmu^\mathrm{T}$ from \eqref{eq:RbarMu}.
Term \X{} gives us the edge mass matrix on boundary edges
\begin{equation}
[\matMmuTilde]_{(\kEdge-1) \Nedge+i,(\kEdge-1)\Nedge+j} \coloneqq \sum_{\Ekn\in\partial\Tk\cap\bcE} \int_{\Ekn} \mujElementView \, \muiElementView \, \ds 
  \label{eq:matMmuTildeEntry}
\end{equation}
for $i,j\in \{1,2,\ldots,\Nedge\}$, $\kEdge\in \{1,2,\ldots,\Kedge\}$ and the matrix entries given in \eqref{eq:matMmuBar1}.
These integrals are over edges on the domain boundary, so that each integral is only evaluated once.

The last term --- term \XI{} --- incorporates boundary data into the hybrid equation. 
For the inflow boundary edges, we obtain a~contribution to the right-hand-side~$\vecKmuD\in\IR^{\Kedge\Nedge}$
\begin{subequations}\label{eq:matK}
\begin{equation}
[\vecKmuD]_{(\kEdge-1)\Nedge+i} = \sum_{\Ekn\in\partial\Tk\cap\bcEin} \,
\int_{\Ekn} c_\mathrm{D} \, \muiElementView\, \ds\,,
\end{equation}
and outflow boundaries add a contribution to the matrix~$\matKmuOut \in\IR^{\Kedge\Nedge\times KN}$ 
\begin{equation}
  [\matKmuOut]_{(\kEdge-1)\Nedge+i,(k-1)N+j} = \sum_{\Ekn\in\partial\Tk\cap\bcEout} \int_{\Ekn} \pkj \, \muiElementView \, \ds
  \label{eq:matKmuOutEntry}
\end{equation}
\end{subequations}
meaning that we obtain a~block matrix similar to~$\matT$ (cf.~\eqref{eq:matU}) from every outflow edge.

\subsection{Time discretization}
\label{sec:timeDiscretization}
The system of equations in~\eqref{eq:systemAsMatrices} can be reformulated in matrix notation as
\begin{equation*}
\begin{aligned}
\mat\matMphi \parT \vecC(t) &\;=\; \vecBphi(t) - \matLbar(t) \, \vecC(t) - \matMbar(t) \, \vecLambda(t) \;\eqqcolon\; \vecR_\mathrm{RK}(t, \vecC(t), \vecLambda(t)) \,,
\\
\vec{0} &\;=\; \vecBmu(t) - \matN{(t)} \vecC(t) - \matP \vecLambda(t) \,.
\end{aligned}
\end{equation*}
This is a~first order differential algebraic equation \cite{NgPe12}, and we use implicit time stepping schemes to discretize it in a~stable manner.
Here, we employ diagonally implicit Runge-Kutta (DIRK) schemes of orders 1 to 4 \cite{Alex77,HaiWanII}.

The time interval $[t^0, t^\text{end}]$ is discretized into not necessarily equidistant points $t^n$ with $t^0 = 0 < t^1 < t^2 < \ldots < t^{\text{end}}$. 
The time step size of the $n$th time step is given by $\dt^{n} \coloneqq t^{n+1} - t^{n}$, and we abbreviate coefficient vectors and matrices on the $n$th time level as~$\vecC^{n}$, etc.
Time step adaptation can be easily achieved in DIRK schemes with embedded error estimators, e.g., as carried out for HDG methods in~\cite{JS13}.

For a~{stiffly accurate} DIRK scheme with $s$ stages, the update at $t^{n+1}$ is obtained by solving
\begin{equation*}
  \begin{aligned}
  \matM_{\varphi} \vecC^{(i)}
  &= \matM_{\varphi} \vecC^{n} 
  + \dt \sum_{j=1}^{i} a_{ij} \, \vecR_{\text{RK}}\left(t^{(j)}, \vecC^{(j)}, \vecLambda^{(j)}\right)
  \\
  \vec0 &= \vecBmu^{(i)} - \matN^{{(i)}} \vecC^{(i)} - \matP \vecLambda^{(i)}\\
  \end{aligned}
  \begin{aligned}
    \qquad i=1\ldots,s,
  \end{aligned}
\end{equation*}
and setting
\begin{equation*}
  \vecC^{n+1} = \vecC^{(s)}
\end{equation*}
with $t^{(i)} = t^{n} + c_{i} \dt^{n}$ and coefficients~$a_{ij}, b_j, c_i$ defined 
in the routine~\code{rungeKuttaImplicit} (see Butcher tableau in Table~\ref{tab:butcher}).
Due to the implicit character of the time iteration, each stage of the DIRK method requires 
solution of a~linear equation system.
This is where the hybridization comes in handy: it reduces the size of the system that has to be solved.
More details on this are given in Sec.~\ref{sec:staticcondensation}.
%
%
\begin{table}
  \centering
  \begin{tabular}{c|ccc}
    $c_1$ & $a_{11}$ &        &  \\ 
    $\vdots$& $\vdots$ & $\ddots$   &  \\ 
    $c_s$ & $a_{s1}$ & $\cdots$   & $a_{ss}$ \\ \hline
    & $b_1$ & $\cdots$ & $b_s$ \\ 
  \end{tabular} 
  \caption{Butcher tableau of a diagonally implicit Runge-Kutta method. Coefficients of the upper triangular part are zeros.}
  \label{tab:butcher}
\end{table}
\begin{remark}
  All of the employed DIRK schemes are $A$- and $L$-stable \cite{HaiWanII}. 
\end{remark}
\begin{remark}
  The DIRK schemes are stiffly accurate, i.e.
\begin{equation*}
  b_{j} = a_{sj}, \quad j=1,\ldots,s\,.
\end{equation*}
  This means in particular that the last update $\vecC^{(s)}$ is actually the updated solution at the new time $t^{n+1}$. 
\end{remark}

\subsection{Static condensation}
\label{sec:staticcondensation}
In the $i$th stage of a~DIRK scheme, one has to solve the following system of equations:
\begin{equation*}
\begin{aligned}
\overbrace{ 
  \left( \matMphi 
  + \aii \, \dt \, \matLbar^{(i)} \right) 
}^{\eqqcolon \matL} \vecC^{(i)} 
+ \overbrace{\phantom{\Big(} \aii \, \dt \, \matMbar^{(i)} }^{\eqqcolon \matM} \vecLambda^{(i)}
&= 
\overbrace{ 
  \matMphi \vecC^{n} 
  + \aii  \, \dt \,\vecBphi^{(i)}
  + \dt \sum_{j=1}^{i-1} a_{ij} \, \vecR_{\text{RK}}(t^{(j)}, \vecC^{(j)}, \vecLambda^{(j)})
}^{=:\vecQ}\,,
\\
 \matN^{{(i)}} \vecC^{(i)} + \matP \vecLambda^{(i)} &= \vecBmu^{(i)}.
\end{aligned}
\end{equation*}
We can write this compactly as
\begin{subequations}\label{eq:staticcond}
\begin{align}
\matL \Ci + \matM \Li &= \vecQ, 
\label{eq:staticcondeq1}
\\
\matN \Ci + \matP \Li &= \vecBmu^{(i)}.
\label{eq:staticcondeq2}
\end{align}
\end{subequations}
The first line \eqref{eq:staticcondeq1} refers to equation \eqref{eq:systemAsMatrices:eq1}, and the second line \eqref{eq:staticcondeq2} to the hybrid part \eqref{eq:systemAsMatrices:eq2}.
Now, we want to demonstrate the solution procedure in detail: 
Unhybridized DG methods would require a~matrix of the same size as $\matL \in \IR^{KN \times KN}$ to be inverted implying rapid growth of the matrix size for high polynomial orders since $N = \frac{(p+1)(p+2)}{2} =\mathcal{O}(p^2)$ in 2D.
This is especially pronounced in comparison to other methods such as continuous finite elements, where the continuity requirements in the discrete space definition reduce the number of degrees of freedom. 
Using Schur complement reduction procedure, we obtain by substitution
\begin{align}
\label{eq:localsolve}
&\Ci \;=\; \matL^{-1} (\vecQ - \matM \Li)\\
\label{eq:hybridsolve}
&\left(- \matN \matL^{-1} \matM + \matP\right) \Li \;=\; \vecBmu^{(i)} - \matN \matL^{-1} \vecQ.
\end{align}
This leads to a memory (and time) efficient procedure to solve this system.
First, we invert $\matL$ to compute $\matL^{-1} \vecQ$ and $\matL^{-1} \matM$, which can be done in a local fashion because $\matL$ is a block diagonal matrix.
Therefore, $\matL^{-1} \vecQ$ and $\matL^{-1} \matM$ are often referred to as local solves. 
Then, we construct \eqref{eq:hybridsolve} using $\matL^{-1} \vecQ$ and $\matL^{-1} \matM$ and solve for $\Li$. 
Once $\Li$ is known, we can update $\Ci$ by substituting the updated value of $\Li$ into \eqref{eq:localsolve}.

The local solves can be implemented in several different ways:
In our \MatOct{} implementation, it turned out to be most efficient to explicitly compute~$\matL^{-1}$ in a~block-wise fashion and then to apply the inverse matrix to~$\vecQ$ and~$\matM$.
Thus we select a~number of elements and invert the corresponding blocks at once instead of inverting all element-blocks separately or inverting the entire matrix at once.
This block-wise inversion is implemented in the routine~\code{blkinv}.
The optimal block size depends on the utilized hardware especially on the cache sizes of the employed CPU.
Heuristically, we determined {$32 \cdot 2^{-p} \cdot N$} to be a~good choice in our case (for hardware details see {Table~\ref{tab:hdg_vs_ldg}}).
Optionally, one could parallelize the local solves since they do not depend on each other.

%% file: sections/implementationNew.tex
A description of data structures and algorithms related to meshing can be found in the first paper in series~\cite{FrankRAK2015}. 
For the sake of completeness, we briefly introduce transformation rules for element and edge integrals that are detailed in previous works~\cite{FrankRAK2015,ReuterAWFK2016} and emphasize differences due to edge unknown~$\lambda_h$.
Finally, we discuss the assembly of the block matrices in Sec.~\ref{sec:assembly}. 

\subsection{Backtransformation to reference element and reference interval}
\label{sec:transformtoThatEhat}

We use an~affine mapping from the reference triangle~$\hat{T} = \{\transpose{[0,0]},\transpose{[1,0]},\transpose{[0,1]}\}$ to any~$T_k = \{\vec{x}_{k1}, \vec{x}_{k2}, \vec{x}_{k3}\} \in \setT_h$,
\begin{equation}\label{eq:affinemappings}
\vec{F}_k :\quad  \hat{T} \ni \hat{\vec{x}}\mapsto \vecc{B}_k \hat{\vec{x}} + \vec{x}_{k1} = \vec{x}\in T_k\,,
\qquad\text{with}\quad
\IR^{2\times2} \ni \vecc{B}_k \coloneqq 
\left[ \vec{x}_{k2}-\vec{x}_{k1} \,\big|\, \vec{x}_{k3}-\vec{x}_{k1} \right]\,.
\end{equation}
It holds $0 < \det \vecc{B}_k = 2|T_k|$, and thus the component-wise definition of the mapping and its inverse read as
\begin{equation*}
\vec{F}_k(\hat{\vec{x}}) \;=\;
\begin{bmatrix}
B_k^{11}\,\hat{x}^1 + B_k^{12}\,\hat{x}^2 + {x}_{k1}^1\\
B_k^{21}\,\hat{x}^1 + B_k^{22}\,\hat{x}^2 + {x}_{k1}^2
\end{bmatrix}
\qquad\text{and}\qquad
\vec{F}_k^{-1}(\vec{x}) \;=\;
\frac{1}{2\,\abs{T_k}}
\begin{bmatrix}
B_k^{22}\,(x^1 - \cjs{x}_{k1}^1) - B_k^{12}\,(x^2 - {x}_{k1}^2)\\
B_k^{11}\,(x^2 - \cjs{x}_{k1}^2) - B_k^{21}\,(x^1 - {x}_{k1}^1)
\end{bmatrix}\;.
\end{equation*}
For functions~$\,w:T_k\rightarrow \IR\,$ and $\hat{w}:\hat{T}\rightarrow \IR\,$, we imply $\,\hat{w}=w\circ \vec{F}_k\,$, i.\,e., $\,w(\vec{x}) = \hat{w}(\hat{\vec{x}})\,$.
The gradient is transformed using the chain rule:
\begin{equation}\label{eq:rule:gradient}
\grad \;=\; \invtrans{\big(\hat{\grad}\vec{F}_k \big)}\,\hat{\grad}
\;=\;
\frac{1}{2\,\abs{T_k}}
\begin{bmatrix}
B_k^{22}\,\partial_{\hat{x}^1} - B_k^{21}\,\partial_{\hat{x}^2}\\
B_k^{11}\,\partial_{\hat{x}^2} - B_k^{12}\,\partial_{\hat{x}^1}
\end{bmatrix},
\end{equation}
where we abbreviated~$\hat{\grad} = \transpose{[\partial_{\hat{x}^1},\partial_{\hat{x}^2} ]}$.
This results in transformation formulas for integrals over an~element~$T_k$ or an~edge~$E_{kn} \subset T_k$ for a~function $w: \Omega \rightarrow \IR$
\begin{subequations}
  \begin{align}
  \label{eq:trafoRule:T}
  \int_{T_k}w(\vec{x})\,\dd\vec{x}
  &\;=\;\frac{\abs{T_k}}{\abs{\hat{T}}} \int_{\hat{T}}w\circ\vec{F}_k(\hat{\vec{x}})\,\dd\hat{\vec{x}}
  \;=\;2\abs{T_k} \int_{\hat{T}}w\circ\vec{F}_k(\hat{\vec{x}})\,\dd\hat{\vec{x}}
  \;=\;2\abs{T_k} \int_{\hat{T}}\hat{w}(\hat{\vec{x}})\,\dd\hat{\vec{x}}
  \;,\\
  \label{eq:trafoRule:E}
  \int_{E_{kn}} w(\vec{x})\,\dd\vec{x} 
  &\;=\; \frac{\abs{E_{kn}}}{\abs{\hat{E}_n}} \int_{\hat{E}_n} w\circ\vec{F}_k(\hat{\vec{x}})\,\dd\hat{\vec{x}}
  \;=\;\frac{\abs{E_{kn}}}{\abs{\hat{E}_n}} \int_{\hat{E}_n} \hat{w}(\hat{\vec{x}})\,\dd\hat{\vec{x}}
  \;.
  \end{align}
\end{subequations}
We introduce the mapping $\ghatn:[0,1] \rightarrow \refEn$ from the reference interval~$[0,1]$ to the $n$th edge of the reference element, which is given by
\begin{equation}
\ghat_1(s) = \begin{bmatrix} 1-s \\ s   \end{bmatrix}\,, \quad
\ghat_2(s) = \begin{bmatrix} 0   \\ 1-s \end{bmatrix}\,, \quad
\ghat_3(s) = \begin{bmatrix} s   \\ 0   \end{bmatrix}\,,
\end{equation}
and use this to transform Eq.~\eqref{eq:trafoRule:E} further
\begin{equation}
  \frac{\abs{E_{kn}}}{\abs{\hat{E}_n}} \int_{\hat{E}_n} \hat{w}(\hat{\vec{x}})\,\dd\hat{\vec{x}}
  =
  \frac{\abs{E_{kn}}}{\abs{\hat{E}_n}} \int_{0}^{1} \hat{w} \circ \ghatn(s)\, \abs{ \ghatn'(s) } \,\ds
  =
 \abs{E_{kn}} \int_{0}^{1} \hat{w} \circ \ghatn(s) \, \ds\,,
 \label{eq:mappingGamma}
\end{equation}
where we use the fact that $\abs{ \ghatn'(s) } = \abs{\Ehatn}$.

{A~difference compared to previous publications in series are edge integrals with basis functions from an~adjoining element and basis functions defined on the edge, e.g., $\sintEkn \pki\,\mkj\,\dx$ with $\kEdge = \rho( k, n )$ as defined in~\eqref{eq:mapkntokEdge}.
They}
are transformed according to transformation rules~\eqref{eq:trafoRule:E},~\eqref{eq:mappingGamma}
%
\begin{equation}
\label{eq:trafoRule:Emu}
\sintEkn  \pki \,\mujElementView\, \dx =
\frac{\left| \Ekn \right|}{\left| \refEn \right|} 
\sintErefk \vphat{i} (\vxhat) \, \mujElementView \circ \Fk (\vxhat) \,\dxhat
=
\abs{E_{kn}}
\sintUnit \vphat{i} \circ \ghatn(s) \,\mhat{j} \circ \betaMap(s) \, \ds\,,
\end{equation}
where we introduced an additional mapping~$\betaMap: [0,1] \rightarrow [0,1]$ that adapts the edge orientation to match the definition of~$\mujElementView = \mkj$ and is defined as
\begin{equation}
  \betaMap(s) =
  \begin{cases}
    s   &\text{if } \kappa(\rho(k,n),1) = k\,, \\
    1-s &\text{if } \kappa(\rho(k,n),2) = k
  \end{cases}
  \label{eq:betaMap}
\end{equation}
with~$\kappa(\rho(k,n),l) = \kappa(\kEdge,l)$ given in~\eqref{eq:mapKappa}.
This does not introduce any further terms into the equation since $\abs{\betaMap'(s)} = 1$.

{
Similarly to other publications in this series~\cite{FrankRAK2015,ReuterAWFK2016}, triangle and edge integrals are approximated using numerical quadrature rules after transformation to reference element~$\hat{T}$ or reference interval~$[0,1]$, respectively.
We abstain from reproducing the details here and refer the interested reader to previous works.
In all numerical experiments, we use quadrature rules of order~$2p+1$ on both elements and edges.
}

\subsection{Program structure}
\label{sec:structure}

\begin{figure}
\begin{center}
\includegraphics[width=\textwidth]{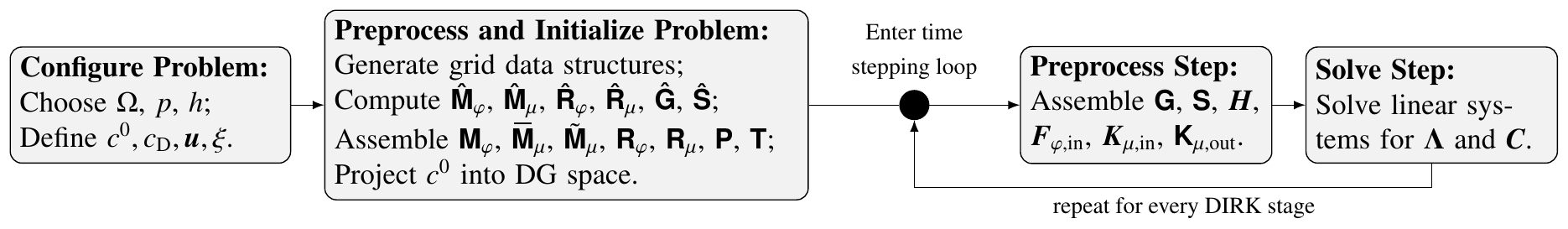}
\end{center}
\caption{\cjs{Structure of the solution algorithm. The last two steps are executed repeatedly as part of the time stepping loop.}}
\label{fig:structure}
\end{figure}

Following the structure outlined in Fig.~\ref{fig:structure}, the implemented solver starts with a~pre-processing phase that includes the definition of solver parameters, initial data, boundary conditions, and the right hand side function followed by the assembly of time-independent matrix blocks and global matrices.
In the time stepping loop, the time-dependent global matrices and right hand side vectors are assembled before solving the resulting linear system.
Note that this has to be done for every DIRK-stage (cf. Sec.~\ref{sec:timeDiscretization}).
Output files in VTK or Tecplot file formats are written after a~user-defined number of time steps.

\subsection{Assembly}
\label{sec:assembly}

In this section, the vectorized assembly of the block matrices in \eqref{eq:systemAsMatrices} is outlined.
For that purpose, the required terms are transformed to reference triangle~$\hat{T}$ or reference interval~$[0,1]$ and then evaluated by numerical quadrature.

As in previous papers in series~\cite{FrankRAK2015,ReuterAWFK2016}, we make extensive use of the Kronecker product $\vecc{A} \otimes \vecc{B}$ of two matrices $\vecc{A} \in \IR^{m_a\times n_a}$, $\vecc{B} \in \IR^{m_b\times n_b}$ defined as
\begin{equation}\label{eq:kron}
\vecc{A} \otimes \vecc{B} \coloneqq \left[ [\vecc{A}]_{i,j} \, \vecc{B}\right] \in \IR^{m_a m_b\times n_a n_b}\,.
\end{equation}
Additionally, we employ operation~$\vecc{A} \otimes_\mathrm{V} \vecc{B}$ with $m_b = r\,m_a, r \in \IN$ introduced in our last publication~\cite{ReuterAWFK2016} as
\begin{equation}\label{eq:kronVec}
\vecc{A} \otimes_\mathrm{V} \vecc{B} \coloneqq \left[ [\vecc{A}]_{i,j} \, [\vecc{B}]_{(i-1)r:ir,:}  \right] \in \IR^{m_b\times n_a n_b} \,,
\end{equation}
which can be interpreted as a~Kronecker product that takes a different right-hand-side for every row of the left-hand side.
This operation is implemented in the routine~\code{kronVec}.
In many cases, we must select edges matching a~certain criterion, e.g., edges in the interior~$\Ekn\in\inE$.
We denote this using the Kronecker delta symbol with a~matching subscript that indicates the criterion to be met, for example
\begin{equation*}
\delta_{\Ekn\in\inE} \coloneqq \left\{ \begin{matrix} 
1 & \text{if } \Ekn\in\inE \,,\\
0 & \text{otherwise} \,.
\end{matrix} \right.
\end{equation*}

Some of the block-matrices in system~\eqref{eq:systemAsMatrices} appeared in identical {or only slightly different} form in previous publications~\cite{FrankRAK2015,ReuterAWFK2016}.
For brevity, we abstain from reproducing the corresponding assembly steps here and refer to the existing descriptions.
This is the case for matrix~$\matRphi$ (identical to~$\matS^\mathrm{diag}$ in~\cite{FrankRAK2015}), source terms~$\vecH$ (cf. any of the papers in series), and mass matrix~$\matMphi$.
Vector~$\vecFphiD$ in equation~\eqref{eq:systemAsMatrices:eq1} contains Dirichlet data provided on inflow boundary edges and is very similar to vector~$\vecK_\mathrm{D}$ in~\cite{ReuterAWFK2016} only lacking the upwinding per quadrature point.
However, here we do not evaluate~$\vec{u}$ in all quadrature points during the assembly of the vector but once per time level and use it then for other terms as well (see Sec.~\ref{sec:assembly:matSmAndmatSout}).
We do the same for Dirichlet data~$c_\mathrm{D}$ and pass them together as per-quadrature-point values to the assembly routine \code{assembleVecEdgePhiIntVal}, which is the same as \code{assembleVecEdgePhiIntFuncContVal} in~\cite{ReuterAWFK2016} without the evaluation of the coefficient function.

Differences to the standard DG implementation in~\cite{ReuterAWFK2016} arise in all matrices that stem from edge integrals involving the new unknown~$\lambda$ (see Sec.~\ref{sec:assembly:matSmAndmatSout}--\ref{sec:assembly:vecKmuD}) and in matrices~$\matG^m$ (see Sec.~\ref{sec:assembly:matG}) that are assembled without projecting~$\vec{u}(t,\vx)$ into the broken polynomial space.


\subsubsection[Assembly of Gm]{Assembly of $\matGm$}
\label{sec:assembly:matG}
For the assembly of matrices~$\matGm$ from~\eqref{eq:matG}, we make use of the transformation rule for the gradient~\eqref{eq:rule:gradient}.
Due to the time-dependent function~$u^m(t,\vec{x})$ in the integrand, we cannot reduce the assembly to Kronecker products of reference matrices as we did for the mass matrix.
We apply transformation rules~\eqref{eq:rule:gradient},~\eqref{eq:trafoRule:T} and obtain
\begin{align*}
\intTk u^1(t, \vec{x}) \, \pkj \, \parX{1} \pki \, \dx
&\approx
\sumR \left(\matB^{22}_{k} \, [\vecU^1]_{k,r} \, [\matGhat]_{1,r,i,j} - \matB^{21}_{k} \,[\vecU^1]_{k,r}\, [\matGhat]_{2,i,j,r}\right) \,, \\
\intTk u^2(t, \vec{x}) \, \pkj \, \parX{2} \pki \, \dx
&\approx
\sumR \left( -\matB^{12}_{k}\, [\vecU^2]_{k,r}\, [\matGhat]_{1,r,i,j} + \matB^{11}_{k}\, [\vecU^2]_{k,r}\, [\matGhat]_{2,i,j,r} \right) 
\end{align*}
with multidimensional array~$\matGhat \in \R^{2\times N\times N\times R}$ that represents a~part of the contribution of the quadrature rule in every integration point~$\hat{\vec{q}}_r$ on the reference element~$\hat{T}$ and arrays~$\vec{U}^m\in\IR^{K\times R}$ that hold the velocity components evaluated in each quadrature point of each element
\begin{equation*}
	[\matGhat]_{m,i,j,r} \coloneqq \omr \, \parXmHat{} \, \phat_{ki} \, \phat_{kj} \,,
  \qquad
  [\vecU^m]_{k,r} \coloneqq u^m(t, \vec{F}_k(\refIntPtTri))\,.
\end{equation*}
The element-local matrix~$\matG^1_\Tk\in\IR^{N\times N}$ is then given as
\begin{align*}
\matG^{1}_{\Tk} 
&= 
\sumR 
\left(
B^{22}_{k}
\begin{bmatrix}
	\omr \, u^1(t,\vec{F}_k(\refIntPtTri)) \, \parX{1} \vphat{1}\, \vphat{1}  & \cdots & \omr \,\, u^1(t,\vec{F}_k(\refIntPtTri)) \, \parX{1} \vphat{1}\, \vphat{N} 
	\\
	\vdots  & \ddots & \vdots 
	\\
	\omr \, u^1(t,\vec{F}_k(\refIntPtTri)) \,   \parX{1} \vphat{N}\, \vphat{1}  & \cdots & \omr \, u^1(t,\vec{F}_k(\refIntPtTri)) \,  \parX{1} \vphat{N}\, \vphat{N} 
\end{bmatrix}\right.
\\
&\phantom{= \sumR }
\left.
- B^{21}_{k}
\begin{bmatrix}
\omr \, u^1(t,\vec{F}_k(\refIntPtTri)) \, \parX{2} \vphat{1} \, \vphat{1}  & \cdots & \omr \, u^1(t,\vec{F}_k(\refIntPtTri)) \,  \parX{2} \vphat{1} \,\vphat{N} 
\\
\vdots  & \ddots & \vdots 
\\
\omr \, u^1(t,\vec{F}_k(\refIntPtTri)) \, \parX{2} \vphat{N}\,\vphat{1}  & \cdots & \omr \, u^1(t,\vec{F}_k(\refIntPtTri)) \,\parX{2} \vphat{N} \,\vphat{N} 
\end{bmatrix}
\right)
\\
&= 
\sumR \left( B^{22}_{k} \,[\vecU^1]_{k,r} \,[\matGhat]_{1,:,:,r} - B^{21}_{k}\, [\vecU^1]_{k,r} \,[\matGhat]_{2,:,:,r} \right)\,.
\end{align*}
On the surface, this procedure appears to closely follow the assembly of~$\matG^m$ in~\cite{ReuterAWFK2016}; however, conceptually, there is a~big difference: here a~quadrature rule is applied to all elements at the same time using blocks of basis functions evaluated in each quadrature point, whereas the procedure to assemble~$\matG^m$ in~\cite{ReuterAWFK2016} builds the full matrix from the contributions of each degree of freedom of the projected DG representation of the velocity using already integrated reference blocks.
To speed up the implementation, we do not assemble a~global sparse matrix in each iteration of the \code{for}-loop over quadrature points and adding to the sparse matrix from the previous iteration.
Instead, we use the standard Kronecker product~\eqref{eq:kron} to build a dense $KN\times N$ vector of blocks, from which the global sparse matrix is constructed using the vectorial Kronecker operator~\eqref{eq:kronVec}
\begin{equation*}
\vecc{G}^1 = 
\sumR \vecc{I}_{K\times K} \otimes_\mathrm{V} \left( 
  \begin{bmatrix} B_1^{22} \, [\vec{U}^1]_{1,r} \\ \vdots \\ B_K^{22} [\vec{U}^1]_{K,r} \end{bmatrix}  \otimes [\matGhat]_{1,:,:,r} -
  \begin{bmatrix} B_1^{21} \, [\vec{U}^1]_{1,r} \\ \vdots \\ B_K^{21} [\vec{U}^1]_{K,r} \end{bmatrix}  \otimes [\matGhat]_{2,:,:,r}
\right) \,,
\end{equation*}
where~$\vecc{I}_{K\times K}$ is the ${K\times K}$~identity matrix.
$\matG^{2}$ is assembled analogously.
The procedure for both matrices is implemented in the routine~\code{assembleMatElemDphiPhiFuncContVec}.

\subsubsection[Assembly of S and Sout]{Assembly of $\matS$ and $\matSout$}
\label{sec:assembly:matSmAndmatSout}
Matrices~$\matS$ and~$\matSout$ are assembled together, the set of relevant edges is expanded to~$\inE\cup\bcEout$.
On a~relevant edge~$\Ekn$, we transform terms of the form given in Eq.~\eqref{eq:SbarEntry} using transformation rule~\eqref{eq:trafoRule:Emu} and approximate the integral by a~one-dimensional numerical quadrature rule
\begin{align*}
\int_{\Ekn} (\vecu \cdot \nkn) \, \pki \, \mujElementView \, \ds 
&= \abs{\Ekn} \int_0^1 \left((\vec{u}(t) \circ \Fk \circ \ghatns) \cdot \nkn \right) \, \vphat{i} \circ \ghatns \, \mhat{j} \circ \betaMap(s) \, \ds \\
&\approx \sum_{r=1}^R \underbrace{\left((\vecu(t) \circ \Fk \circ \ghatn(\hat{q}_r)) \cdot \nkn \right)}_{\eqqcolon [\vecU_{\vNormal}]_{k,n,r}} \, \underbrace{\omr \,  \vphat{i} \circ \ghatn(\hat{q}_r) \, \mhat{j} \circ \betaMap(\hat{q}_r)}_{\eqqcolon [\matShat]_{i,j,n,r,l}} \,,
\end{align*}
where~$\vecU_{\vNormal}\in\IR^{K\times 3\times R}$ holds the normal velocity evaluated in each quadrature point, and the subscript~$l$ in~$\matShat \in \IR^{N\times \Nedge \times 3 \times R \times 2}$ covers the two cases of~$\betaMap$ in~\eqref{eq:betaMap}.
This allows us to assemble the global matrix as
\begin{equation}
\matS = \sum_{n=1}^3 \underbrace{\begin{bmatrix} 
  \delta_{E_{1n}=E_{1}} & \dots & \delta_{E_{1n}=E_{\Kedge}} \\
  \vdots & \ddots & \vdots \\
  \delta_{E_{Kn}=E_{1}} & \dots & \delta_{E_{Kn}=E_{\Kedge}}
\end{bmatrix}}_{\eqqcolon \permutMat} \otimes_\mathrm{V} \left( 
\sum_{r=1}^R \sum_{l=1}^2 \begin{bmatrix}
  \delta_{E_{1n}\in\inE\cup\bcEout} \, \abs{E_{1n}} \, 
    \delta_{\kappa(\rho(1,n),l) = 1} \, [\vecU_{\vNormal}]_{1,n,r} \\
  \vdots \\
  \delta_{E_{Kn}\in\inE\cup\bcEout} \, \abs{E_{Kn}} \, 
    \delta_{\kappa(\rho(K,n),l) = K} \, [\vecU_{\vNormal}]_{1,n,r} 
\end{bmatrix} \otimes [\matShat]_{:,:,n,r,l} \right)\,.
\label{eq:Delta_n}
\end{equation}
Here we introduce the permutation matrix~$\permutMat\in\IR^{K\times\Kedge}$, $n\in\{1,2,3\}$ that has a single entry per row indicating the correspondence $\Ekn = \Ekbar$ for all elements and edges.
It takes care of the necessary permutation from the element-based view of the assembly to the edge-based view of the edge degrees of freedom.
The assembly of~$\matS$ is implemented in the routine~\code{assembleMatEdgePhiIntMuVal}.

\subsubsection[Assembly of Rmu, T, and Kmuout]{Assembly of $\matRmu$, $\matT$, and $\matKmuOut$}
\label{sec:assembly:matRmuAndmatT}

Matrices $\matRmu$, $\matT$, and $\matKmuOut$ are all constructed using similar terms with the only difference being the set of edges considered or the roles of~$\vphi$ and~$\mu$ interchanged.
The terms in Eqs.~\eqref{eq:RmuEntry},~\eqref{eq:matTentry}, or~\eqref{eq:matKmuOutEntry} are transformed using transformation rule~\eqref{eq:trafoRule:Emu} yielding
\begin{equation*}
	\sintEkn  \pki \, \mujElementView \, \ds	=
	\abs{\Ekn}
	\underbrace{\sintUnit \vphat{i} \circ \ghatns \, \mhat{j} \circ \betaMap(s) \, \ds}_{\eqqcolon [\matRmuHat]_{i,j,n,l}}\,,
\end{equation*}
where the index~$l$ in~$\matRmuHat\in\IR^{N\times\Nedge\times3\times2}$ plays the same role as it did for~$\matShat$ before covering the two cases of~$\betaMap$.
Thus we can assemble~$\matRmu$ and~$\matT$ as
\begin{equation*}
\matRmu 
= \sum_{n=1}^3 \sum_{l=1}^2 \left( \begin{bmatrix}
  \abs{E_{1n}}\, \delta_{E_{1n}\in\inE} & & \\
  & \ddots & \\
  & & \abs{E_{Kn}}\, \delta_{E_{Kn}\in\inE}
\end{bmatrix} \, \permutMat \right) \otimes [\matRmuHat]_{:,:,n,l} 
= \transpose{\matT} 
\end{equation*}
with~$\permutMat$ from Eq.~\eqref{eq:Delta_n}.
These matrices are time-independent and thus assembled only once in~\code{preprocessProblem}, using routine~\code{assembleMatEdgePhiIntMu}.

Matrix~$\matKmuOut$ only differs in the set of edges considered, hence it can be assembled almost identically as
\begin{equation*}
\transpose{(\matKmuOut)}
= \sum_{n=1}^3 \sum_{l=1}^2 \left( \begin{bmatrix}
  \abs{E_{1n}}\, \delta_{E_{1n}\in\bcEout} & & \\
  & \ddots & \\
  & & \abs{E_{Kn}}\, \delta_{E_{Kn}\in\bcEout}
\end{bmatrix} \, \permutMat  \right) \otimes [\matRmuHat]_{:,:,n,l} \,.
\end{equation*}
Note that, in fact, we assemble the transpose of $\matKmuOut$ and thus can reuse the same function.
However, due to the time-dependent velocity field, the set of outflow edges~$\bcEout$ can change over time, and we have to do this in every stage of the time stepping method.


\subsubsection[Assembly of MmuBar and MmuTilde]{Assembly of $\matMmuBar$ and $\matMmuTilde$}
\label{sec:assembly:matMmuBarAndMatMuTilde}
For the hybrid mass matrices~$\matMmuBar$ and~$\matMmuTilde$ (cf. Eq.~\eqref{eq:matMmuBar1} and~\eqref{eq:matMmuTildeEntry}), we apply transformation rule~\eqref{eq:trafoRule:Emu} to obtain
\begin{equation*}
\int_{\Ekn} \mujElementView \, \muiElementView \, \ds 
= \abs{\Ekn} \int_0^1 \mhat{j} \circ \betaMap(s) \, \mhat{i} \circ \betaMap(s) \, \ds
= \abs{\Ekn} \underbrace{\int_0^1 \mhat{j}(s) \, \mhat{i}(s) \, \ds}_{\eqqcolon [\hat{\matM}_\mu]_{i,j}}\,.
\end{equation*}
With the help of permutation matrices~$\permutMat$ (cf. Eq.~\eqref{eq:Delta_n}), we obtain the global matrices
\begin{align*}
\matMmuBar &= 
\sum_{n=1}^3 \left( \transpose{(\permutMat)} \, \begin{bmatrix}
\abs{E_{1n}}\, \delta_{E_{1n}\in\inE} & & \\
 & \ddots &  \\
  &  & \abs{E_{Kn}} \,\delta_{E_{Kn}\in\inE}
\end{bmatrix} \, \permutMat \right) \otimes  \hat{\matM}_\mu\,,\quad
\\
\matMmuTilde &=
\sum_{n=1}^3 \left( \transpose{(\permutMat)} \, \begin{bmatrix}
\abs{E_{1n}}\, \delta_{E_{1n}\in\bcE} & & \\
 & \ddots &  \\
  &  & \abs{E_{Kn}} \,\delta_{E_{Kn}\in\bcE}
\end{bmatrix} \, \permutMat \right) \otimes  \hat{\matM}_\mu\,,
\end{align*}
which are implemented in a~common assembly routine~\code{assembleMatEdgeMuMu}.

\subsubsection[Assembly of KmuIn]{Assembly of $\vecKmuD$}
\label{sec:assembly:vecKmuD}
Last, we consider the Dirichlet boundary contributions on inflow boundary edges in the hybrid equation.
We transform the term in Eq.~\eqref{sec:assembly:vecKmuD} as before and approximate it using numerical quadrature
\begin{equation*}
\int_{\Ekn} c_\mathrm{D}\,\muiElementView \, \ds 
= \abs{\Ekn} \int_0^1 c_\mathrm{D}(t) \circ \Fk \circ \ghatns \, \mhat{i}(s) \, \ds
\approx \abs{\Ekn} \sum_{r=1}^R \omr \, c_\mathrm{D}(t) \circ \Fk \circ \ghatn(\hat{q}_r) \, \mhat{i}(\hat{q}_r) \,.
\end{equation*}
We can omit the mapping~$\betaMap$ here, since we only consider boundary edges, and our numbering of the mesh entities ensures that only the first case of the definition of~$\betaMap$ (cf. Eq.~\eqref{eq:betaMap}) is relevant here.
Thus the global vector is assembled as
\begin{equation*}
\vecKmuD = 
\sum_{n=1}^3 \abs{\Ekn} \sum_{r=1}^R \omr \, \mhat{i}(\hat{q}_r) \, \transpose{(\permutMat)} \, 
\begin{bmatrix} \delta_{E_{1n}\in\bcEin} & & \\ & \ddots & \\ & & \delta_{E_{Kn}\in\bcEin} \end{bmatrix} [\vecC_\mathrm{D}]_{:,n,r} \,,
\end{equation*}
which is implemented in~\code{assembleVecEdgeMuFuncCont}.

%% file: sections/results.tex
In this section, we verify our implementation by means of convergence experiments followed by some performance analysis of the code. In addition, a~runtime comparison of the presented HDG discretization to a~time-implicit version of the DG discretization from the previous publication~\cite{ReuterAWFK2016} is given to illustrate some points mentioned in this work.

\subsection{Analytical convergence tests}

Our implementation is verified by comparing the experimental orders of convergence to the analytically predicted ones for smooth solutions.
For that, we choose an~exact solution~$c(t,\vec{x})$ and a~velocity field~$\vecu(t,\vec{x})$, with which we derive boundary data~$c_\mathrm{D}$ and source term~$\source$ analytically by substituting~$c$ and~$\vecu$ into~\eqref{eq:model}--\eqref{eq:flux}.
The discretization error~$\|c_h - c\|_{L^2(\Omega)}$ is computed as the~$L^2$-norm of the difference between the numerical and the analytical solutions at the end time (cf.~\cite{FrankRAK2015}).
From that, the experimental orders of convergence~$\text{EOC}$ is given by 
\begin{equation*}
\text{EOC} \coloneqq \ln\left(\frac{\|c_{h_{j-1}} - c\|_{L^2(\Omega)}}{\|c_{h_{j}} - c\|_{L^2(\Omega)}}\right) \Bigg/ \ln\left(\frac{h_{j-1}}{h_j} \right)\,.
\end{equation*}

\subsubsection{Steady problem}
\label{sec:steadyproblem}

\begin{table}[!ht]
\small
\begin{tabularx}{\linewidth}{@{}LLllLllLllLllLl@{}}\toprule
$p$ & \multicolumn{2}{l}{0} && \multicolumn{2}{l}{1} && \multicolumn{2}{l}{2} && \multicolumn{2}{l}{3} && \multicolumn{2}{l}{4} \\
\cmidrule{2-3} \cmidrule{5-6} \cmidrule{8-9} \cmidrule{11-12} \cmidrule{14-15}
$j$ & $\|c_h-c\|$ & $\text{EOC}$ && $\|c_h-c\|$ & $\text{EOC}$ && $\|c_h-c\|$ & $\text{EOC}$ && $\|c_h-c\|$ & $\text{EOC}$ && $\|c_h-c\|$ & $\text{EOC}$\\
\midrule
1 & 2.69e-01 &  --- && 6.42e-02 &  --- && 7.35e-03 &  --- && 1.50e-03 &  --- && 1.87e-04 &  --- \\
2 & 1.86e-01 & 0.53 && 1.75e-02 & 1.88 && 7.41e-04 & 3.31 && 9.79e-05 & 3.94 && 6.16e-06 & 4.92 \\
3 & 1.18e-01 & 0.66 && 4.32e-03 & 2.02 && 8.55e-05 & 3.12 && 6.26e-06 & 3.97 && 1.95e-07 & 4.98 \\
4 & 6.87e-02 & 0.78 && 1.07e-03 & 2.01 && 1.04e-05 & 3.03 && 3.95e-07 & 3.99 && 6.11e-09 & 4.99 \\
5 & 3.78e-02 & 0.86 && 2.68e-04 & 2.00 && 1.30e-06 & 3.01 && 2.48e-08 & 3.99 && 1.92e-10 & 5.00 \\
6 & 2.00e-02 & 0.92 && 6.71e-05 & 2.00 && 1.62e-07 & 3.00 && 1.55e-09 & 4.00 && 6.00e-12 & 5.00 \\
\bottomrule
\end{tabularx}
\caption{$L^2(\Omega)$ discretization errors for the steady problem in Sec.~\ref{sec:steadyproblem} and experimental orders of convergence for different polynomial degrees. We have~$h_j = \frac{1}{3\cdot2^j}$ and $K=18\cdot4^j$ triangles in the \mbox{$j$th} refinement level.}
\label{tab:conv:steady}
\end{table}

To verify our spatial discretization we choose the exact solution~$c(\vec{x}) \coloneqq \cos(7x^1)\,\cos(7x^2)$ and velocity field~$\vec{u}(\vec{x})\coloneqq \transpose{[\exp((x^1+x^2)/2), \exp(x^1-x^2)/2)]}$ on the unit square~$(0,1)^2$.
We omit the time discretization and solve the problem directly for~$c_{h_j}$ on a~sequence of increasingly finer meshes with element sizes $h_j = \frac{1}{3\cdot 2^j}$ yielding the expected orders of convergence~$\text{EOC} = p + 1$ for $p > 0$ as listed in Table~\ref{tab:conv:steady}.

\subsubsection{Unsteady problem (ODE)}
\label{sec:unsteadyproblem:ode}

\begin{table}[!ht]
\small
\begin{tabularx}{\linewidth}{@{}LLllLllLllLllLl@{}}\toprule
$p$ & \multicolumn{2}{l}{0} && \multicolumn{2}{l}{1} && \multicolumn{2}{l}{2} && \multicolumn{2}{l}{3} && \multicolumn{2}{l}{4} \\
\cmidrule{2-3} \cmidrule{5-6} \cmidrule{8-9} \cmidrule{11-12} \cmidrule{14-15}
$j$ & $\|c_h-c\|$ & $\text{EOC}$ && $\|c_h-c\|$ & $\text{EOC}$ && $\|c_h-c\|$ & $\text{EOC}$ && $\|c_h-c\|$ & $\text{EOC}$ && $\|c_h-c\|$ & $\text{EOC}$\\
\midrule
1 & 4.25e-02 &  --- && 8.30e-05 &  --- && 6.79e-06 &  --- && 1.13e-08 &  --- && 1.13e-08 &  --- \\
2 & 2.14e-02 & 0.99 && 2.13e-05 & 1.96 && 8.53e-07 & 3.00 && 7.20e-10 & 3.97 && 7.20e-10 & 3.97 \\
3 & 1.08e-02 & 0.99 && 5.40e-06 & 1.98 && 1.07e-07 & 3.00 && 4.54e-11 & 3.99 && 4.54e-11 & 3.99 \\
4 & 5.39e-03 & 1.00 && 1.36e-06 & 1.99 && 1.34e-08 & 3.00 && 2.85e-12 & 3.99 && 2.86e-12 & 3.99 \\
5 & 2.70e-03 & 1.00 && 3.40e-07 & 2.00 && 1.67e-09 & 3.00 && 2.10e-13 & 3.76 && 2.43e-13 & 3.56 \\
\bottomrule
\end{tabularx}
\caption{$L^2(\Omega)$ discretization errors for the unsteady problem in Sec.~\ref{sec:unsteadyproblem:ode} measured at $t_{\text{end}} = 2$ and experimental orders of convergence for different polynomial degrees. We have~$h_j = \frac{1}{3\cdot2^4}$,~$K=4608$, and~$\Delta t_j = \frac{1}{5\cdot2^j}$ in the \mbox{$j$th} refinement level.}
\label{tab:conv:unsteady:ode}
\end{table}

We test our implementation of the DIRK schemes (cf. Sec.~\ref{sec:timeDiscretization}) using exact solution~$c(t,\vec{x}) = \exp(-t)$ and velocity field~$\vecu(t,\vec{x}) = \vec{0}$ on~$\Omega = (0,1)^2$ in the time interval~$J=[0,2]$.
The mesh size is fixed as~$h_j = \frac{1}{3\cdot2^4}$, and the time step size is~$\Delta t_j = \frac{1}{5\cdot2^j}$.
The order of the DIRK scheme is chosen to be~$\min(p+1,4)$. Table~\ref{tab:conv:unsteady:ode} shows the expected orders of convergence~$\text{EOC} = p + 1$ for $p < 4$.

\subsubsection{Unsteady problem (PDE)}
\label{sec:unsteadyproblem:pde}

\begin{table}[!ht]
\small
\begin{tabularx}{\linewidth}{@{}LLllLllLllLllLl@{}}\toprule
$p$ & \multicolumn{2}{l}{0} && \multicolumn{2}{l}{1} && \multicolumn{2}{l}{2} && \multicolumn{2}{l}{3} && \multicolumn{2}{l}{4} \\
\cmidrule{2-3} \cmidrule{5-6} \cmidrule{8-9} \cmidrule{11-12} \cmidrule{14-15}
$j$ & $\|c_h-c\|$ & $\text{EOC}$ && $\|c_h-c\|$ & $\text{EOC}$ && $\|c_h-c\|$ & $\text{EOC}$ && $\|c_h-c\|$ & $\text{EOC}$ && $\|c_h-c\|$ & $\text{EOC}$\\
\midrule
1 & 2.69e-01 &  --- && 6.42e-02 &  --- && 7.35e-03 &  --- && 1.50e-03 &  --- && {1.87e-04} & { ---} \\ 
2 & 1.86e-01 & 0.53 && 1.75e-02 & 1.88 && 7.41e-04 & 3.31 && 9.79e-05 & 3.94 && {6.16e-06} & {4.92} \\ 
3 & 1.18e-01 & 0.66 && 4.32e-03 & 2.02 && 8.55e-05 & 3.12 && 6.26e-06 & 3.97 && {1.96e-07} & {4.97} \\ 
4 & 6.86e-02 & 0.78 && 1.07e-03 & 2.01 && 1.04e-05 & 3.03 && 3.96e-07 & 3.98 && {7.44e-09} & {4.72} \\ 
5 & 3.78e-02 & 0.86 && 2.68e-04 & 2.00 && 1.30e-06 & 3.01 && 2.52e-08 & 3.97 && {7.74e-10} & {3.27} \\ 
\bottomrule
\end{tabularx}
\caption{$L^2(\Omega)$ discretization errors for the unsteady problem in Sec.~\ref{sec:unsteadyproblem:pde} measured at $t_{\text{end}} = 2$ and experimental orders of convergence for different polynomial degrees. We have~$h_j = \frac{1}{3\cdot2^j}$,~$K=18\cdot4^j$ triangles, and~$\Delta t_j = \frac{1}{5\cdot2^j}$ in the \mbox{$j$th} refinement level. For $p = 4$, we reduced the time step to~$\Delta t_j = \frac{1}{5\cdot2^{j+2}}$.}
\label{tab:conv:unsteady:pde}
\end{table}

To produce a~combined verification of the time and space discretization, we use exact solution given by $c(t,\vec{x}) = \cos(7x^1)\,\cos(7x^2)\,+\,\exp(-t)$ and velocity field~$\vec{u}(\vec{x})\coloneqq \transpose{[\exp((x^1+x^2)/2), \exp(x^1-x^2)/2)]}$ on~$\Omega = (0,1)^2$ in the time interval~$J=[0,2]$.
Mesh size and time step are refined as~$h_j = \frac{1}{3\cdot2^j}$ and~$\Delta t_j = \frac{1}{5\cdot2^j}$, respectively.
The order of the DIRK scheme is chosen as~$\min(p+1,4)$, and Table~\ref{tab:conv:unsteady:pde} shows the expected orders of convergence~$\text{EOC} = p + 1$ for $p > 0$.

\subsection{Comparison to the unhybridized DG implementation}
\label{sec:solid_body}


\begin{figure}[!ht]
\centering
\hfill
\begin{subfigure}[t]{.378\textwidth}
\includegraphics[width=\textwidth]{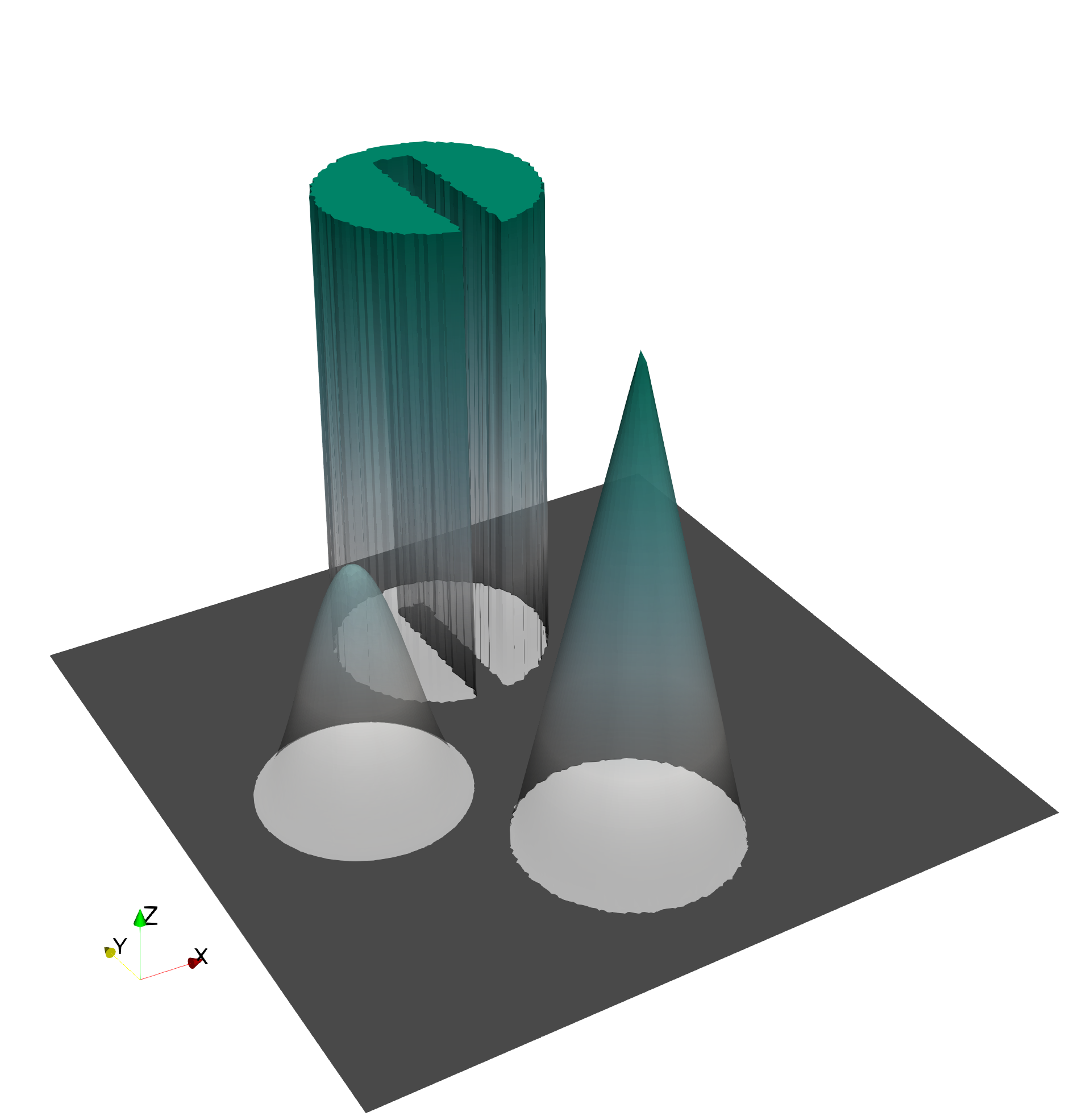}
\caption{Exact solution.}
\label{fig:solid-body:3d:hdg-initial}
\end{subfigure}%
\hfill
\begin{subfigure}[t]{.055\textwidth}
\end{subfigure}
\hfill%
\begin{subfigure}[t]{.378\textwidth}
\includegraphics[width=\textwidth]{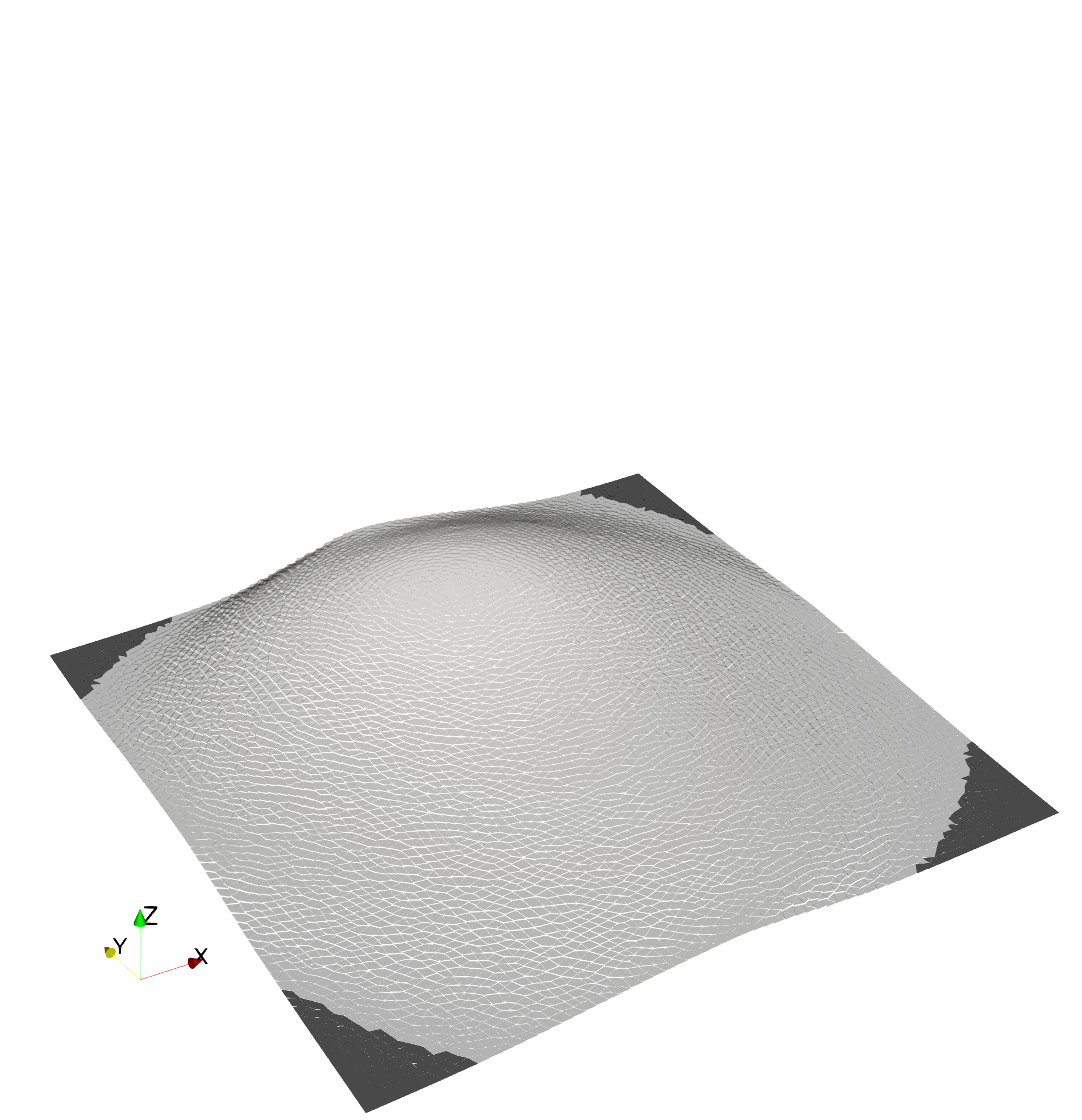}
\caption{$p=0$.}
\label{fig:solid-body:3d:hdg_0}
\end{subfigure}%
\hfill
\\
\hfill
\begin{subfigure}[t]{.378\textwidth}
\includegraphics[width=\textwidth]{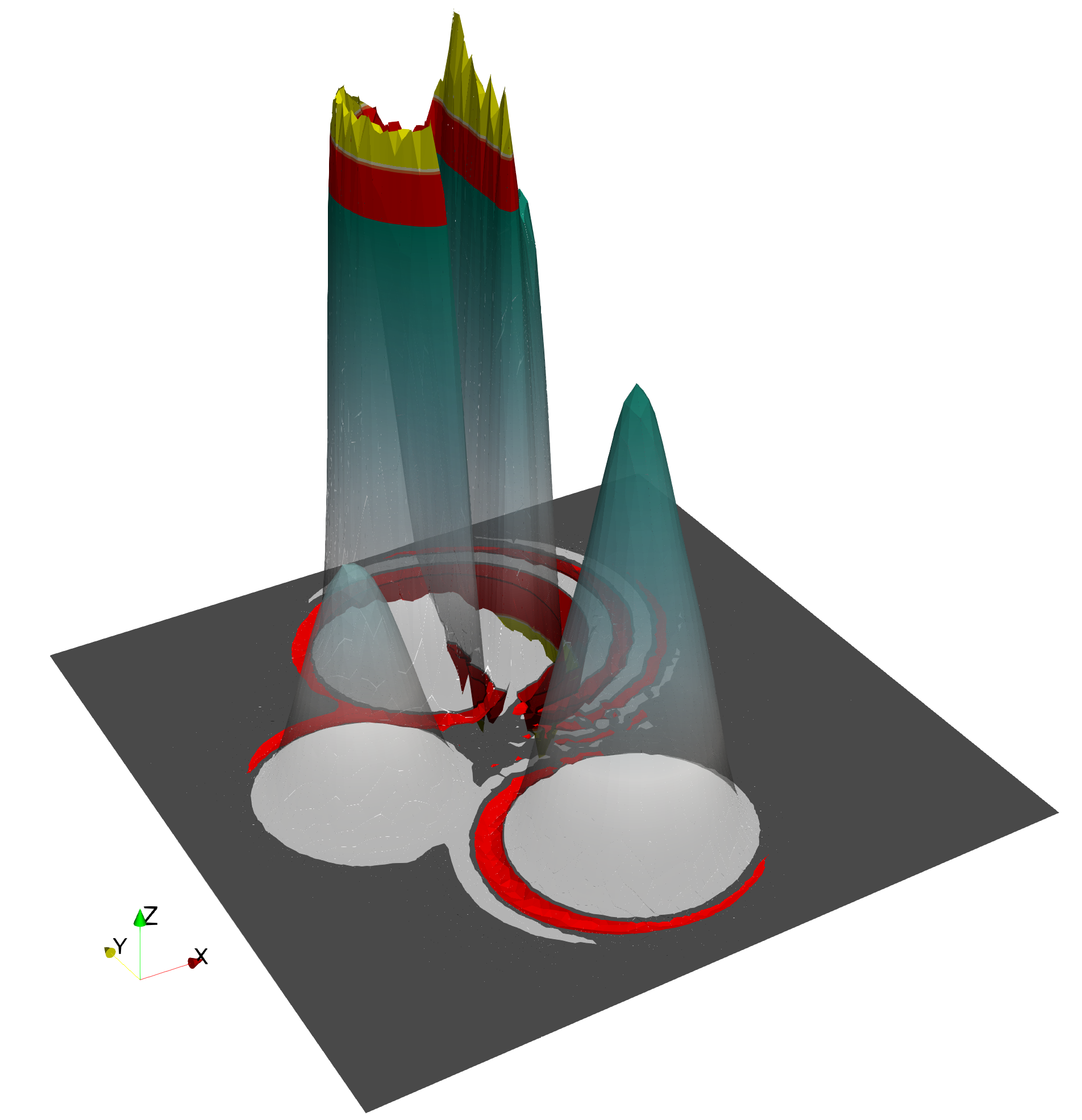}
\caption{$p=1$.}
\label{fig:solid-body:3d:hdg_1}
\end{subfigure}%
\hfill%
\begin{subfigure}[t]{.055\textwidth}
\includegraphics[width=\textwidth]{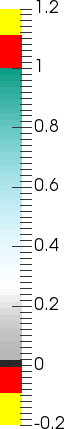}
\end{subfigure}
\hfill%
\begin{subfigure}[t]{.378\textwidth}
\includegraphics[width=\textwidth]{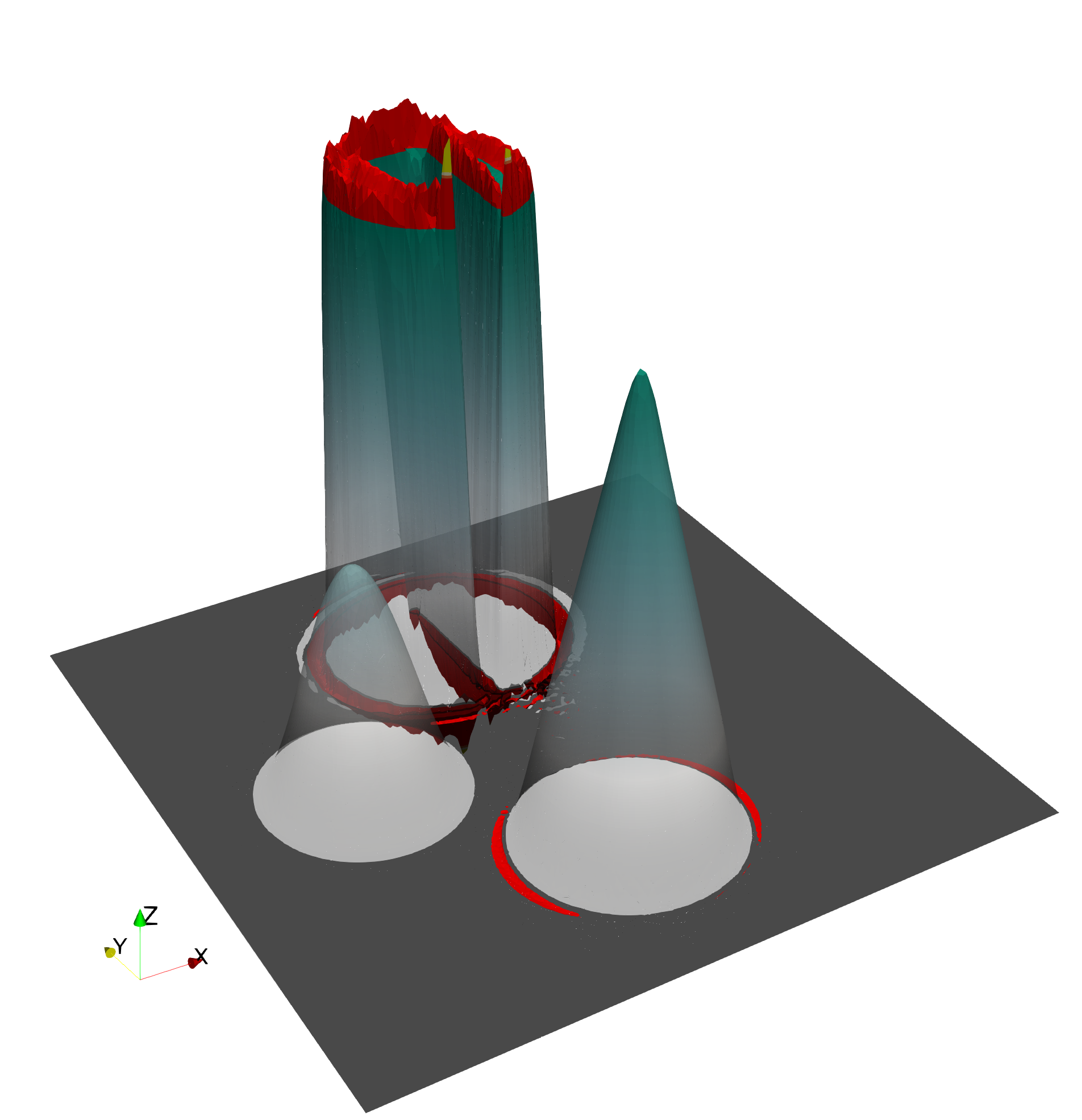}
\caption{$p=2$.}
\label{fig:solid-body:3d:hdg_2}
\end{subfigure}%
\hfill
\\
\hfill
\begin{subfigure}[t]{.378\textwidth}
\includegraphics[width=\textwidth]{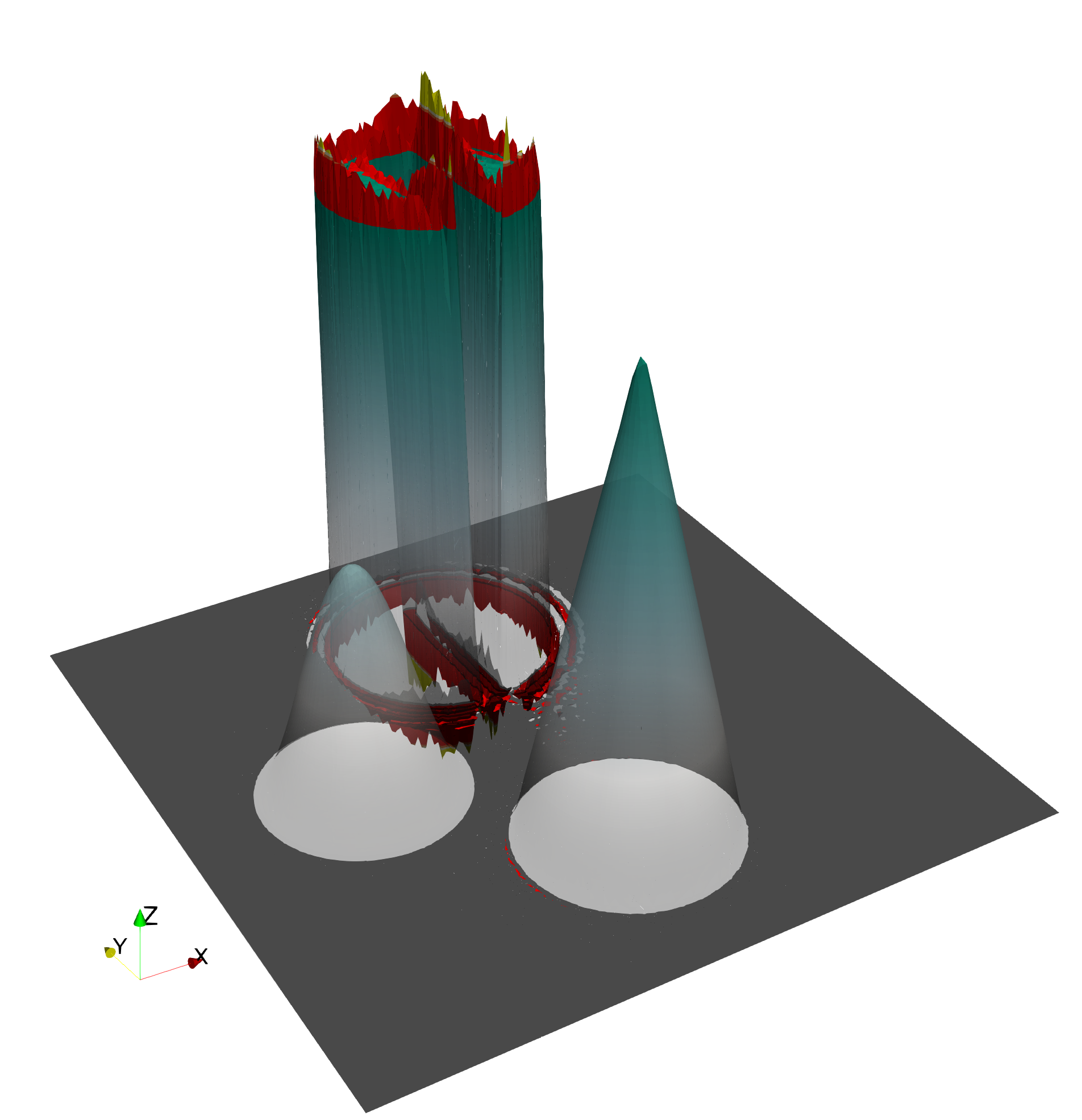}
\caption{$p=3$.}
\label{fig:solid-body:3d:hdg_3}
\end{subfigure}%
\hfill
\begin{subfigure}[t]{.055\textwidth}
\end{subfigure}
\hfill%
\begin{subfigure}[t]{.378\textwidth}
\includegraphics[width=\textwidth]{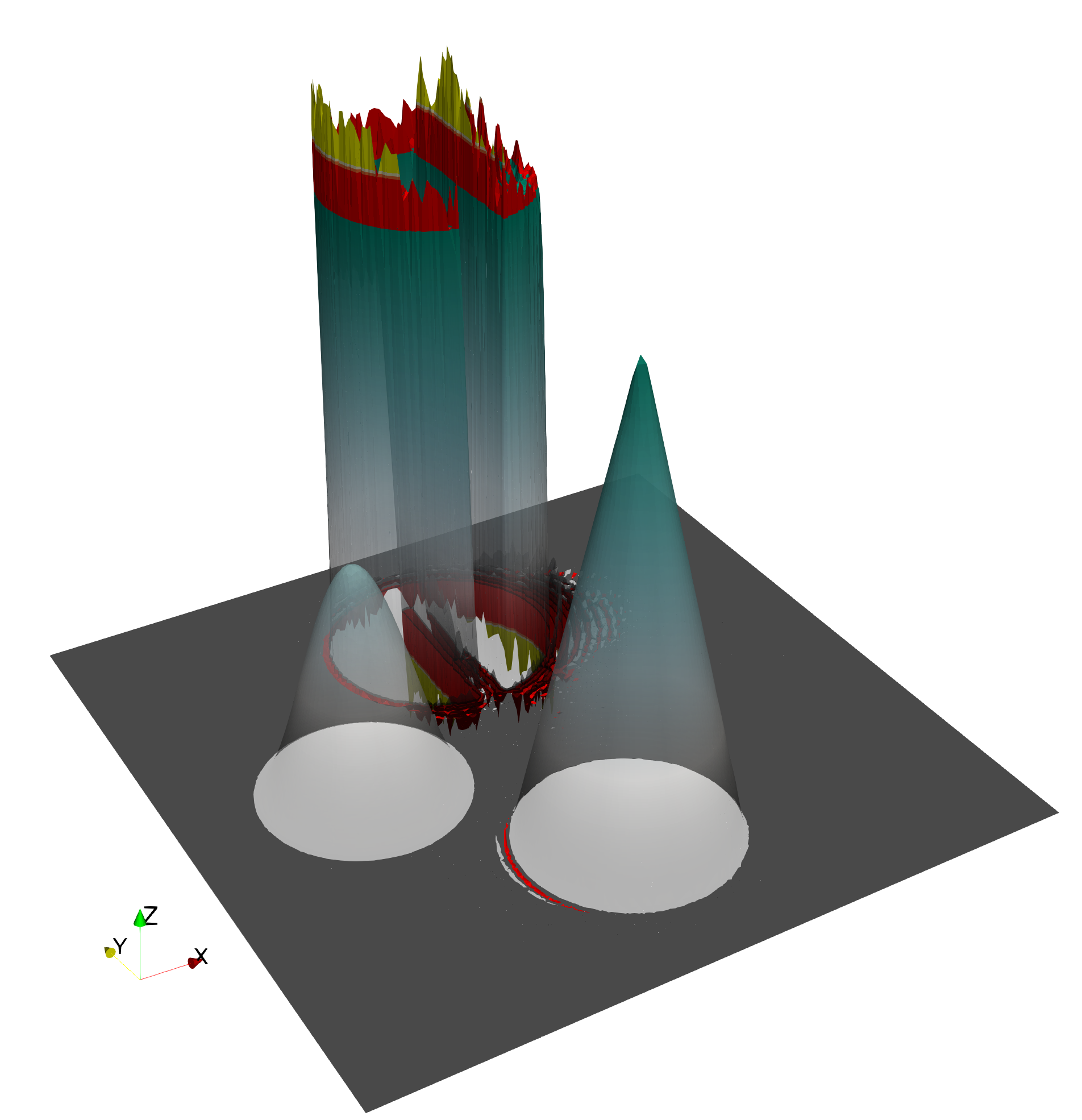}
\caption{$p=4$.}
\label{fig:solid-body:3d:hdg_4}
\end{subfigure}%
\hfill
\caption{HDG~solutions for the solid body rotation benchmark for different polynomial orders at end time~$t_\mathrm{end}=2\pi$.}
\label{fig:solid-body:3d:hdg}
\end{figure}

\begin{figure}[!ht]
\centering
\hfill
\begin{subfigure}[t]{.378\textwidth}
\includegraphics[width=\textwidth]{fig/initial_data}
\caption{Exact solution.}
\label{fig:solid-body:3d:ldg-initial}
\end{subfigure}%
\hfill
\begin{subfigure}[t]{.055\textwidth}
\end{subfigure}
\hfill
\begin{subfigure}[t]{.378\textwidth}
\includegraphics[width=\textwidth]{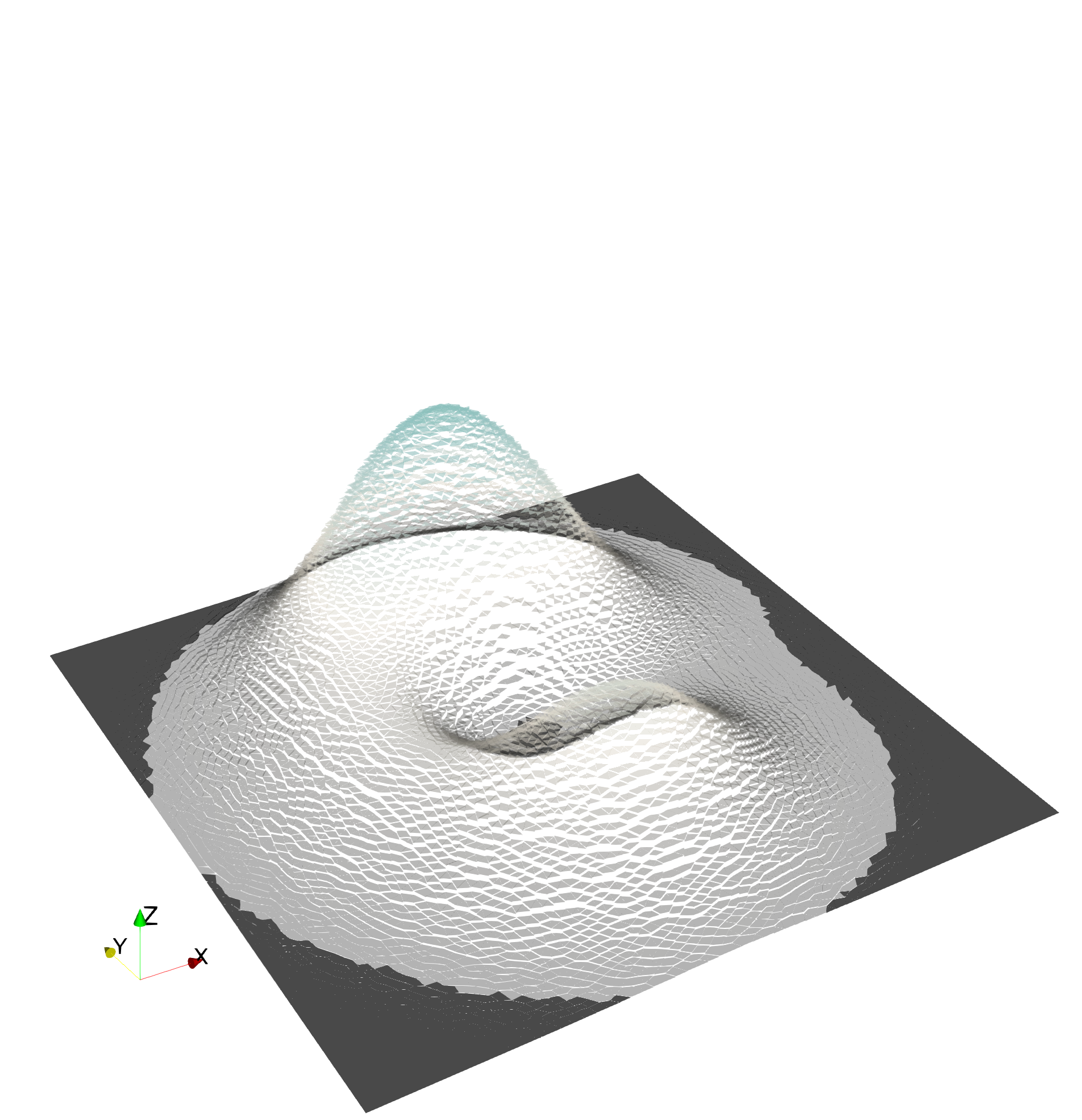}
\caption{$p=0$.}
\label{fig:solid-body:3d:ldg_0}
\end{subfigure}%
\hfill
\\
\hfill
\begin{subfigure}[t]{.378\textwidth}
\includegraphics[width=\textwidth]{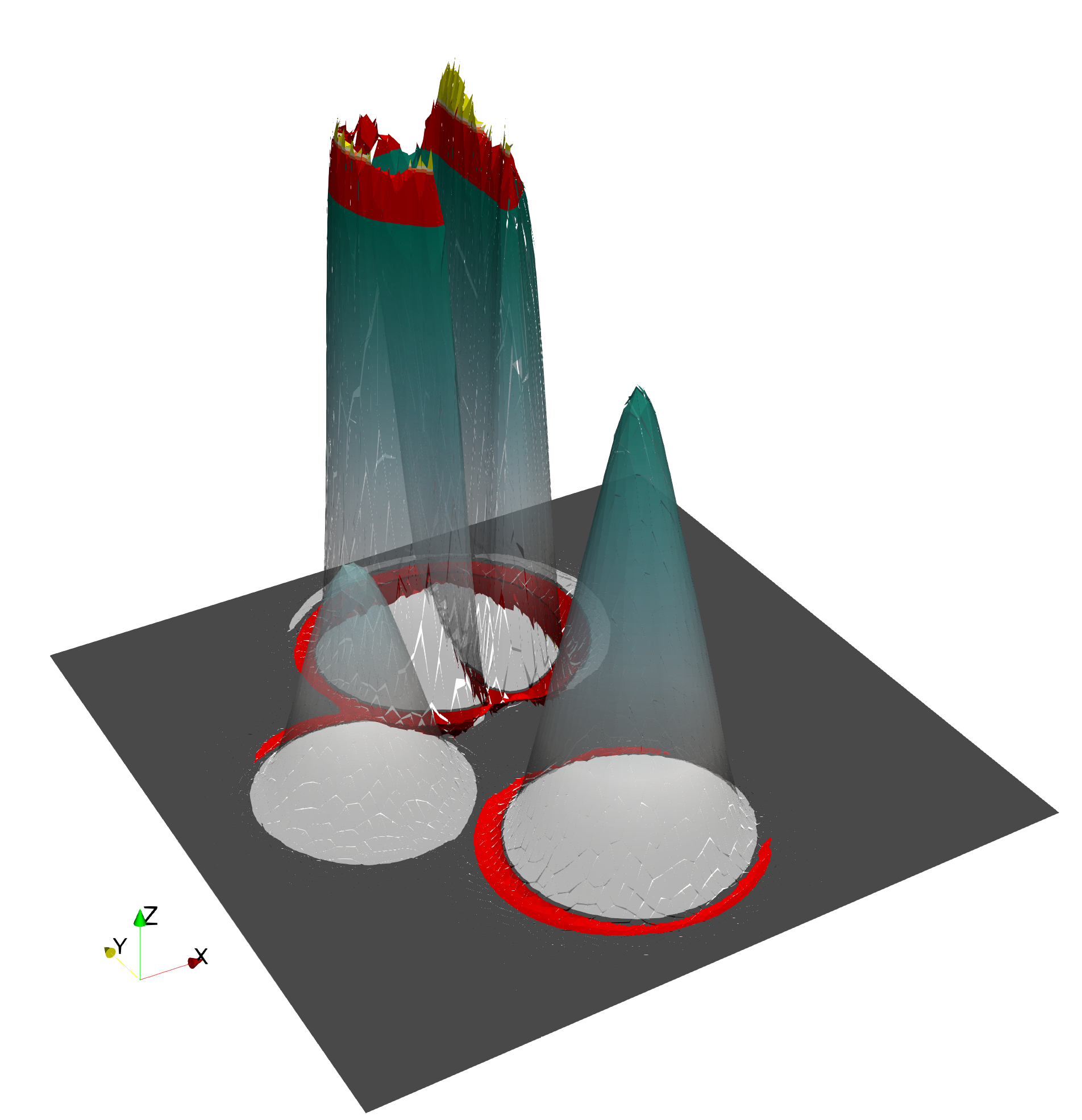}
\caption{$p=1$.}
\label{fig:solid-body:3d:ldg_1}
\end{subfigure}%
\hfill
\begin{subfigure}[t]{.055\textwidth}
\includegraphics[width=\textwidth]{fig/legend_vert}
\end{subfigure}
\hfill
\begin{subfigure}[t]{.378\textwidth}
\includegraphics[width=\textwidth]{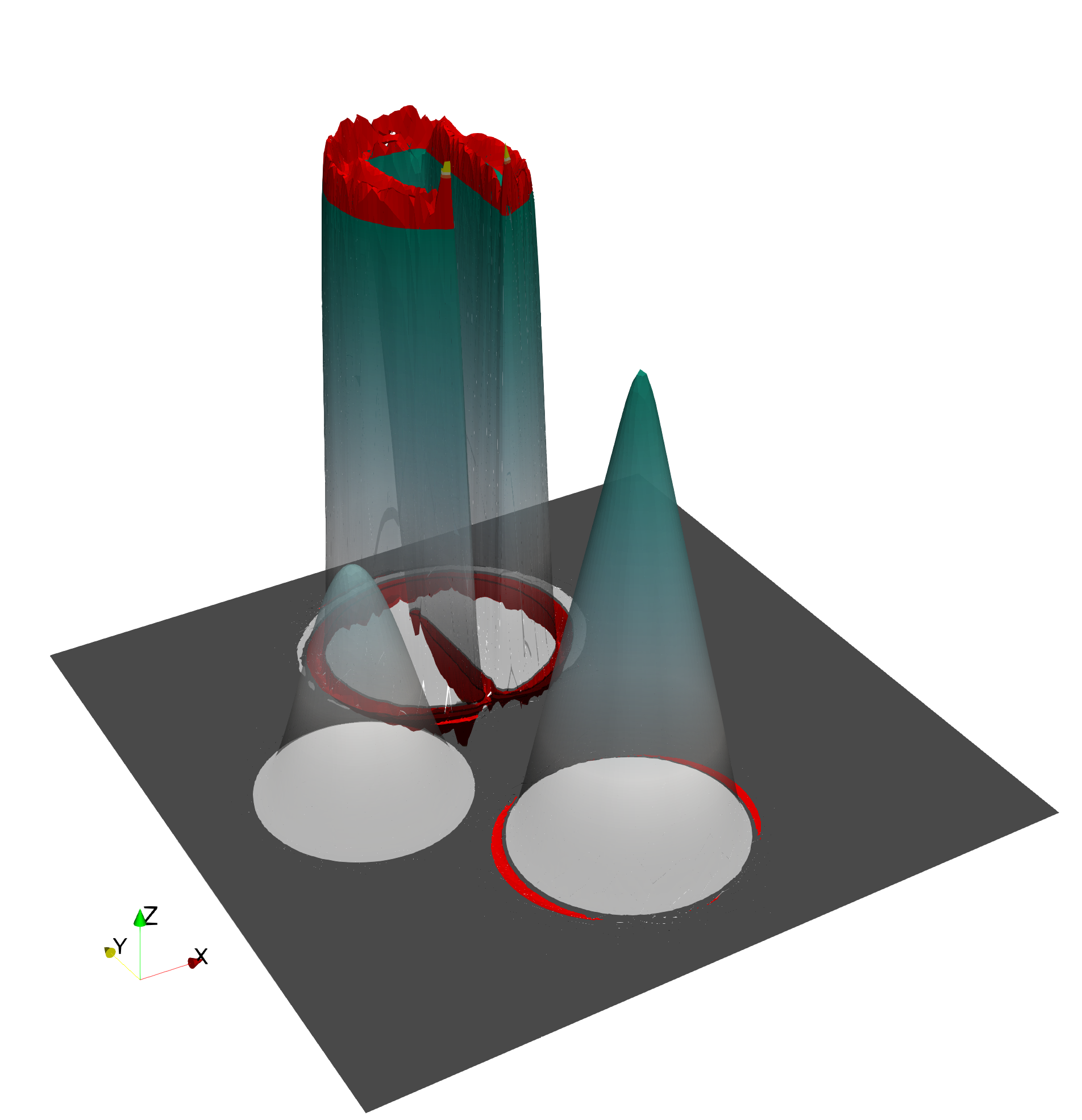}
\caption{$p=2$.}
\label{fig:solid-body:3d:ldg_2}
\end{subfigure}%
\hfill
\\
\hfill
\begin{subfigure}[t]{.378\textwidth}
\includegraphics[width=\textwidth]{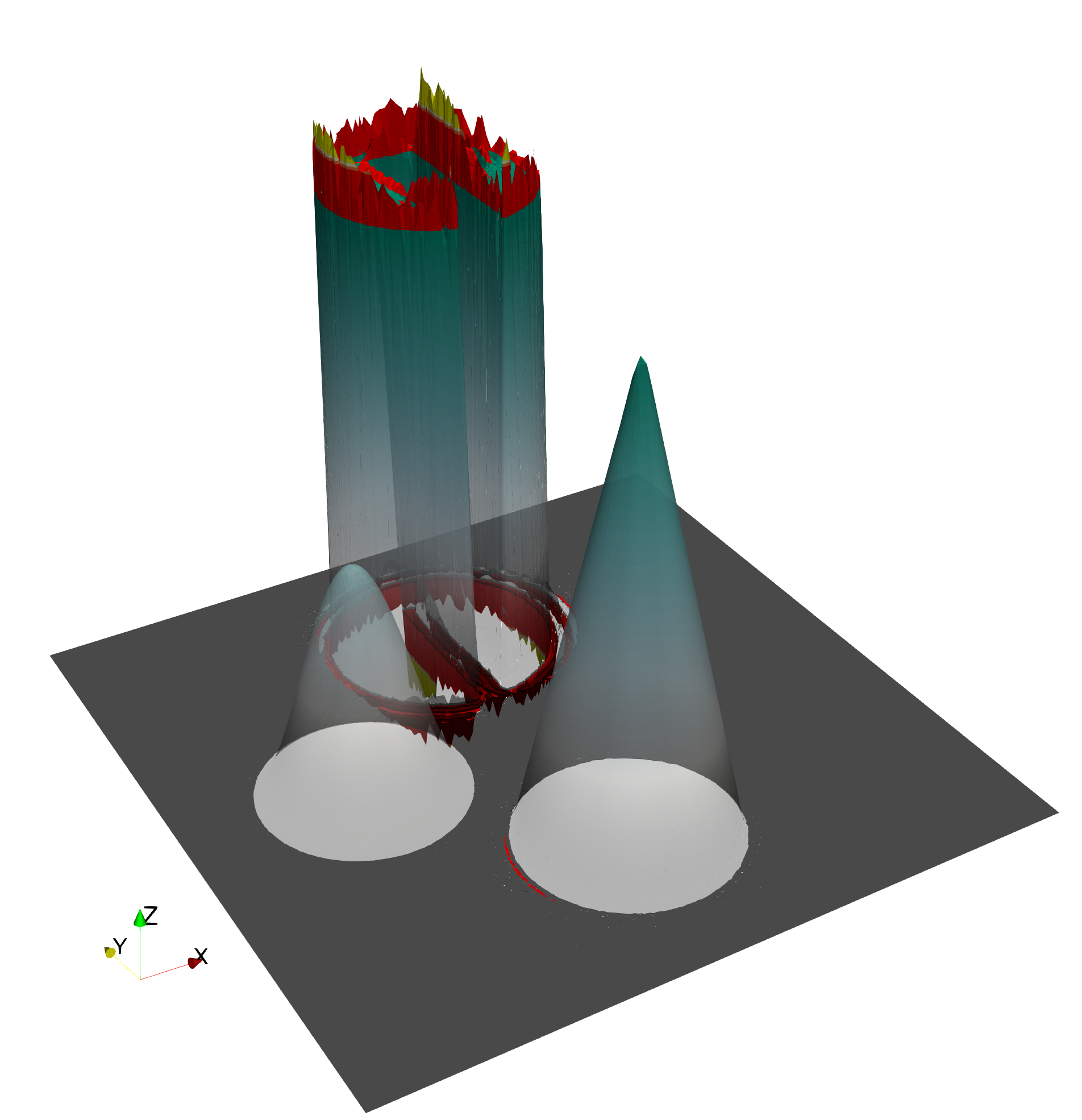}
\caption{$p=3$.}
\label{fig:solid-body:3d:ldg_3}
\end{subfigure}%
\hfill
\begin{subfigure}[t]{.055\textwidth}
\end{subfigure}
\hfill
\begin{subfigure}[t]{.378\textwidth}
\includegraphics[width=\textwidth]{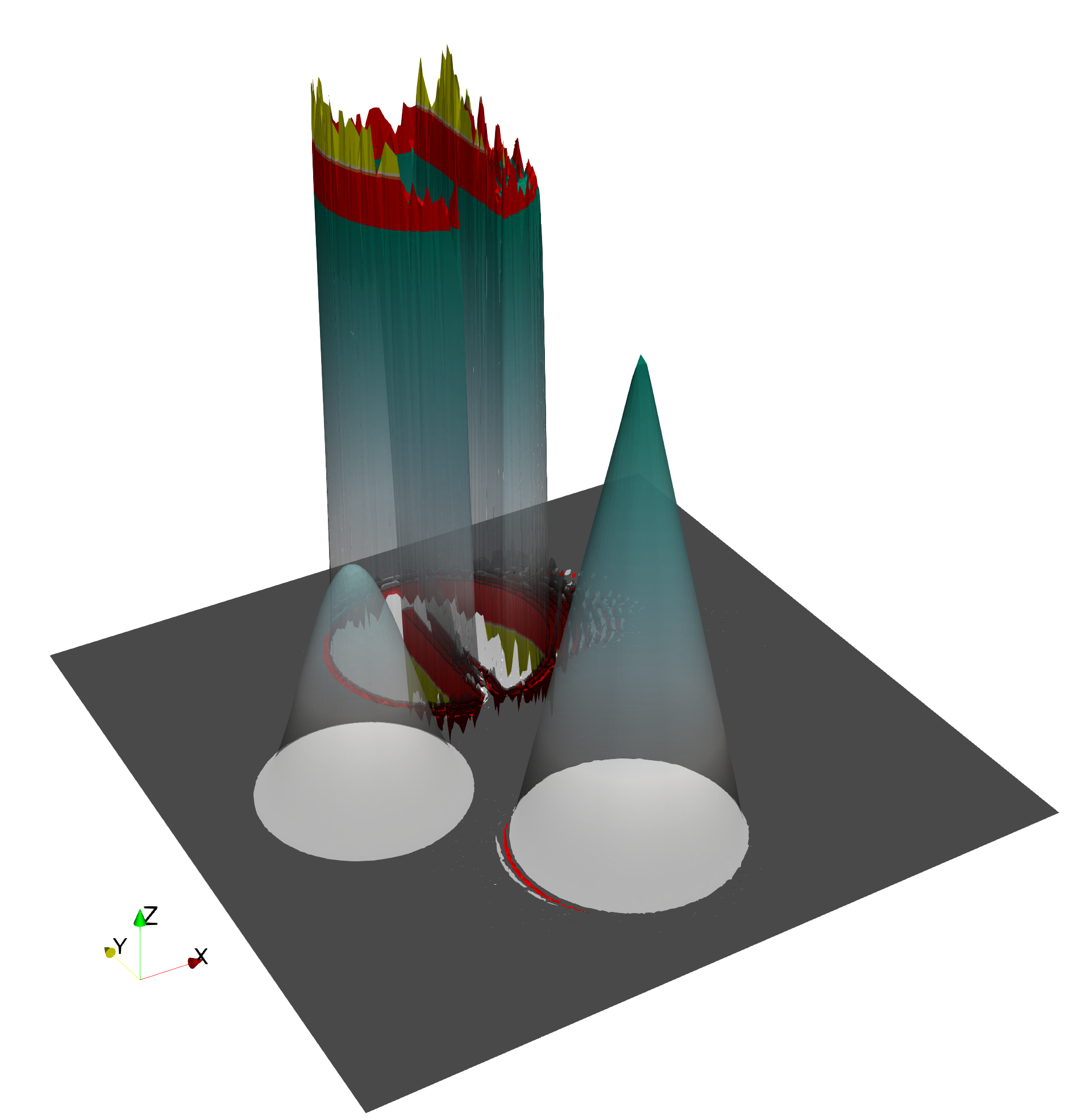}
\caption{$p=4$.}
\label{fig:solid-body:3d:ldg_4}
\end{subfigure}%
\hfill
\caption{DG~solutions for the solid body rotation benchmark for different polynomial orders at end time~$t_\mathrm{end}=2\pi$.}
\label{fig:solid-body:3d:ldg}
\end{figure}

\begin{figure}[!ht]
\centering
\includegraphics[width=.7\textwidth]{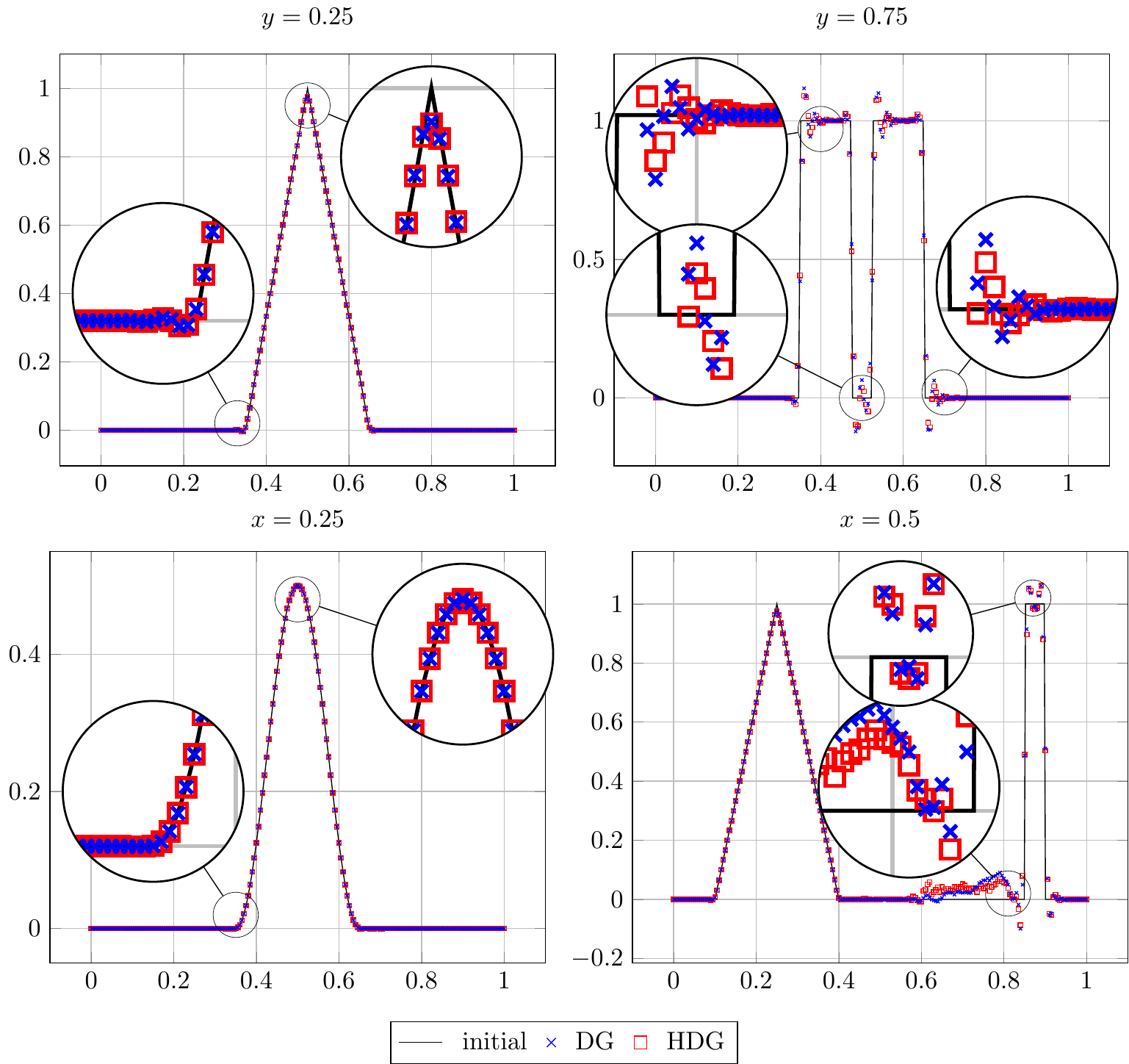}
\caption{Cross-sections of the DG solutions for the solid body rotation benchmark with~$p=3$ at end time $t_\mathrm{end}=2\pi$.}
\label{fig:solid-body:line}
\end{figure}

Our previous work in series was concerned with the same model problem~$\eqref{eq:model}$ and the same linear flux~$\eqref{eq:flux}$; however, an~\emph{unhybridized discontinuous Galerkin} discretization with explicit strong stability preserving Runge--Kutta methods was used---in contrast to a~HDG discretization and diagonally \emph{implicit} Runge-Kutta schemes utilized in the current study.
To make both implementations comparable, we implemented a~variant of the DG solver that incorporates the DIRK schemes and computes the solid body rotation benchmark proposed by LeVeque~\cite{LeVeque1996} and used in~\cite{ReuterAWFK2016} to investigate the performance of slope limiters.
The initial solution consists of a~slotted cylinder, a sharp cone, and a smooth hump (see Fig.~\ref{fig:solid-body:3d:hdg-initial}) placed in a~square domain~$\Omega=(0,1)^2$ with velocity field~$\vecu(\vx) = \transpose{[0.5 - x^2, x^1 - 0.5]}$ producing a~full counterclockwise rotation of the initial scene over time interval~$J=(0,2\pi)$.
With $r = 0.0225$ and
$G(\vec{x},\vec{x}_0) \coloneqq \frac{1}{0.15} \|\vec{x}-\vec{x}_0\|_2$,
we choose homogeneous Dirichlet boundary conditions $c_\mathrm{D} = 0$ and right-hand-side function $\source = 0$ with initial data satisfying
\begin{equation*}
c^0(\vec{x}) = \left\{
\begin{array}{lll}
  1                                               & \quad\text{if}\quad \parbox{.36\textwidth}{$(x^1 - 0.5)^2 + (x^2 - 0.75)^2 \le r$\\$\land\;(x^1\le0.475 \lor x^1\ge0.525 \lor x^2\ge0.85)$} & \mbox{(slotted cylinder)} \\
  1-G(\vec{x},\transpose{[0.5,0.25]})                       & \quad\text{if}\quad (x^1 - 0.5)^2 + (x^2 - 0.25)^2 \le r & \mbox{(sharp cone)}       \\
  \frac{1}{4}(1+\cos(\pi G(\vec{x},\transpose{[0.25,0.5]})) & \quad\text{if}\quad (x^1 - 0.25)^2 + (x^2 - 0.5)^2 \le r & \mbox{(smooth hump)}      \\
  0                                               & \quad\text{otherwise} & 
\end{array}
\right\}
\end{equation*}
An~unstructured mesh generated by MATLAB's \texttt{initmesh} with maximum element size~$h=2^{-6}$ is used resulting in 14006 elements and 320 time steps.

\begin{table}[!ht]
\small
\parbox[b][][t]{.7\linewidth}{%
\begin{tabularx}{\linewidth}{@{}LLLlLL@{}}\toprule
& \multicolumn{2}{l}{$\|c_h(t_\mathrm{end}) - c^0\|_{L^2(\Omega)}$} & & \multicolumn{2}{l}{runtime $[s]$} \\
\cmidrule{2-3}\cmidrule{5-6}
{$p$} & DG & HDG && DG & HDG \\
\midrule
0 & 1.87e-01 & 2.35e-01 && 27.1  & 39.8 \\
1 & 7.25e-02 & 7.78e-02 && 157   & 218 \\
2 & 5.53e-02 & 5.44e-02 && 797   & 717 \\
3 & 4.02e-02 & 4.13e-02 && 3980  & 3166 \\
4 & 4.16e-02 & 4.18e-02 && 10996 & 7199 \\
\bottomrule
\end{tabularx}
}
\hfill
\parbox[b][][t]{.27\linewidth}{%
\begin{tabularx}{\linewidth}{@{}lL@{}}\toprule
\multicolumn{2}{l}{Hardware\,/\,software details} \\
\midrule
CPU & Intel Core-i7 4790 (Haswell) \\
RAM & 32 GBytes \\
MATLAB & R2017a \\
\bottomrule
\end{tabularx}
\vfill
}
\caption{Comparison of $L^2$-errors for the solid body rotation benchmark (cf. Sec.~\ref{sec:solid_body}) at end time~$t_\mathrm{end}=2\pi$ (using initial data as the exact solution) and runtimes for DG and HDG solvers (left).
Details of the employed hardware and software (right).
}
\label{tab:hdg_vs_ldg}
\end{table}

Figures~\ref{fig:solid-body:3d:hdg} and~\ref{fig:solid-body:3d:ldg} show the computed solution at end time~$t_\mathrm{end}=2\pi$ for different polynomial degrees.
We chose a~color map (inspired by~\cite{Badia2017}) that emphasizes violations of the discrete maximum principle: 
in the range~$[0,1]$ we have a~color gradient from black via white to green;
values in the range~$[-0.1,0)\cup(1,1.1]$ are colored red, and values in~$(-\infty,-0.1)\cup(1.1,+\infty)$ are colored yellow.
To make oscillations in the bottom range better visible, we gradually reduce the opacity from above to zero.
Figure~\ref{fig:solid-body:line} presents intersection lines for~$p=3$.

Clearly, the lowest-order approximation is unusable for this kind of problem in both implementations with numerical diffusion killing off most (DG) or all (HDG) of the variability of the solution.
For~$p>0$, solutions from HDG and DG are in good agreement for all approximation orders.
This finding is substantiated by the intersection lines in Fig.~\ref{fig:solid-body:line} and the $L^2$-errors shown in Table~\ref{tab:hdg_vs_ldg} exhibiting only minor differences between both discretizations.
However, also clearly visible are severe violations of the discrete maximum principles and oscillations in the wake of cylinder and cone, which do not become less pronounced with increasing approximation order.
This type of behavior can be alleviated using slope limiters as shown in our previous work~\cite{ReuterAWFK2016}, where a~post-processing step in each time level restricted the updated solution at the vertices to the bounds provided by the mean values of the adjacent elements.
Unfortunately, designing slope limiters for implicit time stepping methods is not a~trivial task and lies out of scope of this work.

When comparing the runtimes for both discretizations in Table~\ref{tab:hdg_vs_ldg}, it becomes clear that the static condensation outlined in Sec.~\ref{sec:staticcondensation} becomes advantageous especially for higher approximation orders making HDG a~superior approach for time-implicit high-order discretizations.

\subsection{Performance analysis}

\begin{table}[!ht]
\small
\begin{tabularx}{\linewidth}{@{}lLLLLLLLLLL@{}}\toprule
 & \multicolumn{5}{l}{$h_\mathrm{max}=2^{-4}$ $(K=872)$} & \multicolumn{5}{l}{$h_\mathrm{max}=2^{-6}$ $(K=14006)$} \\
\cmidrule(lr){2-6} \cmidrule(l){7-11}
 & \mbox{runtime} & \multicolumn{2}{l}{assembly step} & \multicolumn{2}{l}{solver step} & runtime & \multicolumn{2}{l}{assembly step} & \multicolumn{2}{l}{solver step} \\
\cmidrule(lr){3-4} \cmidrule(lr){5-6} \cmidrule(lr){8-9} \cmidrule(l){10-11}
$p$ & & $\matG$ & $\matS$ & $\matL^{-1}$ & $\lambda$ & & $\matG$ & $\matS$ & $\matL^{-1}$ & $\lambda$ \\
\midrule
0 &  0.641 & \multicolumn{2}{c}{0.181 (28.2$\%$)} & \multicolumn{2}{c}{0.047  (7.3$\%$)} &   1.156 & \multicolumn{2}{c}{  0.367 (31.7$\%$)} & \multicolumn{2}{c}{ 0.617 (53.4$\%$)} \\
&& 18.2$\%$ &  8.3$\%$ & 25.5$\%$ & 61.7$\%$ && 13.4$\%$ & 19.6$\%$ &  4.5$\%$ & 91.9$\%$ \\
1 &  0.919 & \multicolumn{2}{c}{0.351 (38.2$\%$)} & \multicolumn{2}{c}{0.278 (30.3$\%$)} &   6.698 & \multicolumn{2}{c}{  1.820 (27.2$\%$)} & \multicolumn{2}{c}{ 4.674 (69.8$\%$)} \\
&& 29.3$\%$ & 17.9$\%$ & 33.1$\%$ & 59.0$\%$ && 42.6$\%$ & 29.3$\%$ & 23.8$\%$ & 71.2$\%$ \\
2 &  1.818 & \multicolumn{2}{c}{0.655 (36.0$\%$)} & \multicolumn{2}{c}{0.852 (46.9$\%$)} &  21.945 & \multicolumn{2}{c}{  6.271 (28.6$\%$)} & \multicolumn{2}{c}{15.491 (70.6$\%$)} \\
&& 41.5$\%$ & 21.4$\%$ & 34.5$\%$ & 56.0$\%$ && 55.8$\%$ & 31.0$\%$ & 25.3$\%$ & 66.1$\%$ \\
3 &  5.627 & \multicolumn{2}{c}{3.079 (54.7$\%$)} & \multicolumn{2}{c}{2.353 (41.8$\%$)} & 100.089 & \multicolumn{2}{c}{ 53.218 (53.2$\%$)} & \multicolumn{2}{c}{46.688 (46.6$\%$)} \\
&& 72.2$\%$ & 15.6$\%$ & 33.1$\%$ & 53.7$\%$ && 79.2$\%$ & 17.0$\%$ & 26.5$\%$ & 60.2$\%$ \\
4 & 10.159 & \multicolumn{2}{c}{6.273 (61.7$\%$)} & \multicolumn{2}{c}{3.685 (36.3$\%$)} & 230.148 & \multicolumn{2}{c}{151.809 (66.0$\%$)} & \multicolumn{2}{c}{78.151 (34.0$\%$)} \\
&& 80.7$\%$ & 12.3$\%$ & 27.3$\%$ & 53.9$\%$ && 88.2$\%$ &  9.9$\%$ & 26.8$\%$ & 56.6$\%$ \\
\bottomrule
\end{tabularx}
\caption{Runtime distribution for 10 time steps of the solid body rotation benchmark (cf.~Sec.~\ref{sec:solid_body}) without initialization tasks (i.e., only time stepping loop) measured using MATLAB's profiler.
We compare different mesh sizes and approximation orders with runtimes given in seconds.
Percentages for assembly and solvers are relative to the runtime of the time stepping loop, percentages for the assembly of~$\matG$ and~$\matS$ (cf.~secs.~\ref{sec:assembly:matSmAndmatSout} and~\ref{sec:assembly:matG}) are relative to assembly runtime, and percentages for inversion of~$\matL$ (cf. sec.~\ref{sec:staticcondensation}) and solving~\eqref{eq:hybridsolve} for~$\lambda$ are relative to solver step runtimes.
Details on the employed hardware and software are given in Table~\ref{tab:hdg_vs_ldg} (right).
}
\label{tab:perf}
\end{table}

A~major advantage of the hybridized DG method is the fact that the globally coupled linear equation system resulting from the discretization is relatively compact and easy to solve (see Sec.~\ref{sec:staticcondensation} for details).
To illustrate this, we present some performance results that show the runtime distribution among the different steps of the code.
We disregard pre-processing and initialization tasks as these are only performed once and consider ten time steps of the solid body rotation benchmark presented in Sec.~\ref{sec:solid_body}. Using MATLAB's profiler to determine the runtime share of each instruction we see --- just as expected --- that the linear solvers together with the assembly of the time-dependent block matrices, particularly~$\matG$ and~$\matS$, are responsible for the majority of the computation time.
Table~\ref{tab:perf} shows the runtime shares for the most expensive parts of the code.

First of all, the overall runtime clearly increases with mesh size and polynomial approximation order simply due to the increasing number of degrees of freedom.
More intriguing is the fact that the assembly step becomes more dominant than the linear system solves with the increasing polynomial degree.
The primary reason for this is the assembly of matrices~$\matG^m$, which do not only grow in size due to the increasing number of degrees of freedom but also require a~loop over all quadrature points (see Sec.~\ref{sec:assembly:matG} for details).
With the order of the quadrature rule (and thus the number of quadrature points) increasing with the polynomial degree, the computational complexity of this operation grows quickly.
Nevertheless, the total runtime increase for higher polynomial approximation orders is not as pronounced as for the unhybridized DG solver as shown in~Sec.~\ref{sec:solid_body}.

In contrast to the shift in the assembly step, the runtime distribution between the block-wise inversion of~$\matL$ and solving~\eqref{eq:hybridsolve} is very similar throughout all approximation orders and different mesh sizes.
We would like to point out that the local solves with~$\matL^{-1}$ are element-local and thus could be easily parallelized with virtually perfect scaling.

%% file: sections/registerofroutines.tex
In this section, the routines added since the first two papers in series \cite{FrankRAK2015,ReuterAWFK2016} are presented in form of a~two-part list:
first, the scripts implementing the solution algorithm in the order they are executed followed by the alphabetically ordered list of assembly and integration routines used in the solution algorithm.
In the code available on GitHub~\cite{FESTUNGGithub}, all routines check for correctly provided function arguments using \MatOct's~\code{validateattributes} excluded here for brevity.

All data that is needed throughout the entire algorithm (e.\,g., mesh data structures, pre-computed reference blocks, etc.) are passed between steps in a~struct always called~\code{problemData}.
In all routines, the argument \code{g} is a struct containing information about the triangulation.
Input parameter \code{N} is the number of degrees of freedom for the 2D polynomials $\phih$, and \code{Nmu} is the number of degrees of freedom for the 1D polynomials $\muh$.
Parameter~\code{qOrd} is the order of the quadrature rule, and~\code{basesOnQuad} is a~struct that contains basis functions~$\phih$ and~$\muh$ evaluated in quadrature points of the reference element~$\hat{T}$, edges~$\hat{E}_n$ of the reference element, and reference interval~$[0,1]$.
In the time stepping routines, \code{nStep} always refers to the current time step index~$n$ and \code{nSubStep} to the current Runge-Kutta stage~$i$ (cf.~Sec.~\ref{sec:timeDiscretization}).

\subsection{Solution algorithm}

\code{problemData = configureProblem(problemData)} is executed first and fills the \code{problemData}-struct with all basic configuration options. 
Problem parameters are to be modified inside this routine.
\begin{lstlisting}
function problemData = configureProblem(problemData)
%% Parameters.
% Choose testcase
problemData = setdefault(problemData, 'testcase', 'solid_body'); 

% Mark run as convergence test (enforces Friedrichs-Keller triangulation)
problemData = setdefault(problemData, 'isConvergence', false);

% Maximum edge length of triangle
problemData = setdefault(problemData, 'hmax', 2^-6);

% Local polynomial approximation order (0 to 4)
problemData = setdefault(problemData, 'p', 2);

% Order of Runge-Kutta method
problemData = setdefault(problemData, 'ordRK', min(problemData.p + 1, 4));

% Order of quadrature rule
problemData = setdefault(problemData, 'qOrd', 2*problemData.p + 1);

% Load testcase
problemData = getTestcase(problemData, problemData.testcase);

% Visualization settings
problemData = setdefault(problemData, 'isVisGrid', false);  % visualization of grid
problemData = setdefault(problemData, 'isVisSol', true);  % visualization of solution
problemData = setdefault(problemData, 'outputFrequency', 50); % no visualization of every timestep
problemData = setdefault(problemData, 'outputBasename', ['output' filesep 'hdg_advection']); 
problemData = setdefault(problemData, 'outputTypes', { 'vtk', 'tec' });  % Type of visualization files

% HDG stabilization parameter
problemData = setdefault(problemData, 'stab', 1.0); 

% Enable blockwise local solves and specify block size
problemData = setdefault(problemData, 'isBlockSolve', true);
problemData = setdefault(problemData, 'blockSolveSize', 2^(5 - problemData.p));

%% Domain and triangulation configuration.
% Triangulate unit square using pdetool (if available or Friedrichs-Keller otherwise).
if ~problemData.isConvergence && license('checkout','PDE_Toolbox')
  problemData.generateGridData = @(hmax) domainPolygon([0 1 1 0], [0 0 1 1], hmax);
else
  fprintf('PDE_Toolbox not available. Using Friedrichs-Keller triangulation.\n');
  problemData.generateGridData = @domainSquare;
end % if
end % function
\end{lstlisting}

\code{problemData = preprocessProblem(problemData)} performs all pre-processing tasks and one-time computations prior to the projection of the initial data and the time stepping loop.
This includes mesh generation, evaluation of basis functions in quadrature points, computation of blocks for the reference elements and edges, and assembly of time-independent block matrices.
\begin{lstlisting}
function problemData = preprocessProblem(problemData)
%% Triangulation.
problemData.g = problemData.generateGridData(problemData.hmax);
if problemData.isVisGrid,  visualizeGrid(problemData.g);  end

%% Globally constant parameters.
problemData.K    = problemData.g.numT;  % number of triangles
problemData.N    = nchoosek(problemData.p + 2, problemData.p); % number of local DOFs
problemData.Nmu  = problemData.p + 1; % number of local DOFs on Faces
problemData.dt   = problemData.tEnd / problemData.numSteps;

% [K x 3] arrays that mark local edges (E0T) or vertices (V0T) that are
% interior or have a certain boundary type.
problemData.g.markE0Tint  = problemData.generateMarkE0Tint(problemData.g);
problemData.g.markE0Tbdr  = problemData.generateMarkE0Tbdr(problemData.g);

% Precompute some repeatedly evaluated fields
problemData.g = computeDerivedGridData(problemData.g);

%% Configuration output.
fprintf('Computing with polynomial order %d (%d local DOFs) on %d triangles and %d edges.\n', ...
        problemData.p, problemData.N, problemData.K, problemData.g.numE)

qOrd = problemData.qOrd;
N = problemData.N;
Nmu = problemData.Nmu;

%% Lookup table for basis function.
problemData.basesOnQuad = computeBasesOnQuad(N, struct, [qOrd, qOrd + 1]);
problemData.basesOnQuad = computeBasesOnQuadEdge(Nmu, problemData.basesOnQuad, [qOrd, qOrd + 1]);

%% Computation of matrices on the reference triangle.
problemData.hatM = integrateRefElemPhiPhi(N, problemData.basesOnQuad, qOrd);
hatMmu = integrateRefEdgeMuMu(Nmu, problemData.basesOnQuad, qOrd);
problemData.hatRmu = integrateRefEdgePhiIntMu([N, Nmu], problemData.basesOnQuad, qOrd);
hatRphi = integrateRefEdgePhiIntPhiInt(N, problemData.basesOnQuad, qOrd);
problemData.hatG = integrateRefElemDphiPhiPerQuad(N, problemData.basesOnQuad, qOrd);
problemData.hatS = integrateRefEdgePhiIntMuPerQuad([N, Nmu], problemData.basesOnQuad, qOrd);

%% Assembly of time-independent global matrices.
problemData.globMphi = assembleMatElemPhiPhi(problemData.g, problemData.hatM);
problemData.globRmu = assembleMatEdgePhiIntMu(problemData.g, problemData.g.markE0Tint, ...
                                              problemData.hatRmu);
problemData.globRphi = assembleMatEdgePhiIntPhiInt(problemData.g, problemData.g.markE0Tint, hatRphi);
problemData.globMmuBar = assembleMatEdgeMuMu(problemData.g, problemData.g.markE0Tint, hatMmu);
problemData.globMmuTilde = assembleMatEdgeMuMu(problemData.g, problemData.g.markE0Tbdr, hatMmu);
problemData.globP = problemData.stab .* problemData.globMmuBar + problemData.globMmuTilde;
problemData.globT = problemData.globRmu';
end % function
\end{lstlisting}

\code{problemData = initializeProblem(problemData)} projects initial data~$c^0$ and visualizes the initial state.
\begin{lstlisting}
function problemData = initializeProblem(problemData)
problemData.isFinished = false;
%% Initial data on elements and faces
if ~problemData.isStationary
  problemData.cDisc = projectFuncCont2DataDisc(problemData.g, problemData.c0Cont, problemData.qOrd, ...
      problemData.hatM, problemData.basesOnQuad);
  fprintf('L2 error w.r.t. the initial condition: %g\n', computeL2Error(problemData.g, ...
      problemData.cDisc, problemData.c0Cont, problemData.qOrd + 1, problemData.basesOnQuad));
  %% visualization of inital condition.
  if problemData.isVisSol
    cLgr = projectDataDisc2DataLagr(problemData.cDisc);
    visualizeDataLagr(problemData.g, cLgr, 'u_h', problemData.outputBasename,0,problemData.outputTypes)
  end
  fprintf('Starting time integration from 0 to %g using time step size %g (%d steps).\n', ...
      problemData.tEnd, problemData.dt, problemData.numSteps)
end % if
end % function
\end{lstlisting}

\code{problemData = preprocessStep(problemData, nStep)} is the first step in each iteration of the time stepping loop and applies mass matrix~$\matMphi$ to solution vector~$\vecC^n$ required in each stage of the DIRK schemes.
\begin{lstlisting}
function problemData = preprocessStep(problemData, nStep)
if ~problemData.isStationary
  problemData.globMcDisc = problemData.globMphi * reshape(problemData.cDisc', [], 1);
end % if
end % function
\end{lstlisting}

\code{problemData = solveStep(problemData, nStep)} is the main routine of the time stepping loop and determines Runge-Kutta coefficients before initiating sub-stepping to carry out the Runge-Kutta stages.
\begin{lstlisting}
function problemData = solveStep(problemData, nStep)
if problemData.isStationary
  problemData.t = 0;
else
  % Obtain Runge-Kutta rule and initialize solution vectors
  [problemData.t, problemData.A, problemData.b] = rungeKuttaImplicit(problemData.ordRK, ...
                                                      problemData.dt, (nStep - 1) * problemData.dt);
  problemData.cDiscRK = cell(length(problemData.t), 1);
end % if

% Carry out RK steps
problemData.isSubSteppingFinished = false;
problemData = iterateSubSteps(problemData, nStep);
end % function
\end{lstlisting}

\code{problemData = preprocessSubStep(problemData, nStep, nSubStep)} is the main assembly routine of the current Runge-Kutta stage, where velocity function~$\vecu$ is evaluated at quadrature points of all edges, inflow and outflow boundaries are determined, and all time-dependent block matrices and right-hand-side vectors are built.
\begin{lstlisting}
function problemData = preprocessSubStep(problemData, nStep, nSubStep)
% Select time level for current Runge-Kutta stage
t = problemData.t(nSubStep);
cDCont = @(x1, x2) problemData.cDCont(t, x1 ,x2);
u1Cont = @(x1, x2) problemData.u1Cont(t, x1, x2);
u2Cont = @(x1, x2) problemData.u2Cont(t, x1, x2);
fCont = @(x1,x2) problemData.fCont(t, x1, x2);

% Evaluate normal advection velocity on every edge.
uNormalQ0E0T = computeFuncContNuOnQuadEdge(problemData.g, u1Cont, u2Cont, problemData.qOrd);

% Determine inflow- and outflow edges
[~, W] = quadRule1D(problemData.qOrd);
markE0TbdrOut = problemData.g.markE0Tbdr & sum(bsxfun(@times, reshape(W,1,1,[]), uNormalQ0E0T), 3) > 0;
markE0TbdrIn = problemData.g.markE0Tbdr & ~markE0TbdrOut;

% Assemble source term.
problemData.globH = problemData.globMphi * reshape(projectFuncCont2DataDisc(problemData.g, fCont, ...
                        problemData.qOrd, problemData.hatM, problemData.basesOnQuad)', [], 1);

% Assemble element integral contributions.
problemData.globG = assembleMatElemDphiPhiFuncContVec(problemData.g, problemData.hatG, ...
                        u1Cont, u2Cont, problemData.qOrd);

% Assemble flux on interior and outflow edges
problemData.globS = assembleMatEdgePhiIntMuVal(problemData.g, problemData.g.markE0Tint|markE0TbdrOut,...
                        problemData.hatS, uNormalQ0E0T);
problemData.globKmuOut = assembleMatEdgePhiIntMu(problemData.g, markE0TbdrOut, problemData.hatRmu)';
                      
% Assemble inflow boundary conditions.
cQ0E0T = computeFuncContOnQuadEdge(problemData.g, cDCont, problemData.qOrd);
problemData.globFphiIn = assembleVecEdgePhiIntVal(problemData.g, markE0TbdrIn, ...
                            cQ0E0T .* uNormalQ0E0T, problemData.N, problemData.basesOnQuad);
problemData.globKmuIn = assembleVecEdgeMuFuncCont(problemData.g, markE0TbdrIn, ...
                            cDCont, problemData.basesOnQuad, problemData.qOrd);
end % function
\end{lstlisting}

\code{problemData = solveSubStep(problemData, nStep, nSubStep)} builds and solves the linear systems described in Sec.~\ref{sec:staticcondensation} for the current Runge-Kutta stage.
\begin{lstlisting}
function problemData = solveSubStep(problemData, nStep, nSubStep)
K = problemData.K;
N = problemData.N;

if problemData.isStationary
  matL = - problemData.globG{1} - problemData.globG{2} + problemData.stab * problemData.globRphi;
  vecQ = problemData.globH - problemData.globFphiIn;
  matM = problemData.globS - problemData.stab * problemData.globRmu;
else
  dtA = problemData.dt * problemData.A;
  
  cDiscRkRHS = zeros(K * N, 1);
  for i = 1 : nSubStep - 1
    cDiscRkRHS = cDiscRkRHS + dtA(nSubStep, i) .* problemData.cDiscRK{i};
  end

  matLbar = - problemData.globG{1} - problemData.globG{2} + problemData.stab * problemData.globRphi;
  matL = problemData.globMphi + dtA(nSubStep, nSubStep) .* matLbar;
  vecBphi = problemData.globH - problemData.globFphiIn;
  vecQ =  problemData.globMcDisc + dtA(nSubStep, nSubStep) * vecBphi + cDiscRkRHS;
  matMbar = problemData.globS - problemData.stab * problemData.globRmu;
  matM = dtA(nSubStep, nSubStep) .* matMbar;
end % if

%% Computing local solves
if problemData.isBlockSolve
  matLinv = blkinv(matL, problemData.blockSolveSize * N);
  LinvQ = matLinv * vecQ;
  LinvM = matLinv * matM;
else
  LinvQ = matL \ vecQ;
  LinvM = matL \ matM;
end % if

%% Solving global system for lambda
matN = - problemData.stab * problemData.globT - problemData.globKmuOut ;
matP = problemData.globP;

lambdaDisc = (-matN * LinvM + matP) \ (problemData.globKmuIn - matN * LinvQ);

%% Reconstructing local solutions from updated lambda
problemData.cDisc = LinvQ - LinvM * lambdaDisc;

if ~problemData.isStationary
  problemData.cDiscRK{nSubStep} = vecBphi - matLbar * problemData.cDisc - matMbar * lambdaDisc;
end % if

problemData.cDisc = reshape(problemData.cDisc, N, K)';
end % function
\end{lstlisting}

\code{problemData = postprocessSubStep(problemData, nStep, nSubStep)} determines whether \code{...SubStep}-routines must be executed once more to carry out the remaining Runge-Kutta stages or if control has to be returned to function~\code{solveStep}.
\begin{lstlisting}
function problemData = postprocessSubStep(problemData, nStep, nSubStep)
problemData.isSubSteppingFinished = problemData.isStationary || nSubStep >= length(problemData.t);
end % function
\end{lstlisting}

\code{problemData = postprocessStep(problemData, nStep)} checks if all time integration steps have been carried out or if another iteration of the time stepping loop is required.
\begin{lstlisting}
function problemData = postprocessStep(problemData, nStep)
problemData.isFinished = problemData.isStationary || nStep >= problemData.numSteps;
end % function

\end{lstlisting}

\code{problemData = outputStep(problemData, outputStep)} writes visualization output files for the updated solution.
\begin{lstlisting}
function problemData = outputStep(problemData, nStep)
%% visualization
if problemData.isVisSol && mod(nStep, problemData.outputFrequency) == 0
  cLagrange = projectDataDisc2DataLagr(problemData.cDisc);
  visualizeDataLagr(problemData.g, cLagrange, 'u_h', problemData.outputBasename, ...
                    nStep, problemData.outputTypes);
end % if
end % function

\end{lstlisting}

\code{problemData = postprocessProblem(problemData)} is executed after termination of the time stepping loop and computes the $L^2$-error of the solution if the analytical solution is provided.
\begin{lstlisting}
function problemData = postprocessProblem(problemData)
%% Visualization
if problemData.isVisSol
  cLagrange = projectDataDisc2DataLagr(problemData.cDisc);
  visualizeDataLagr(problemData.g, cLagrange, 'u_h', problemData.outputBasename, ...
                    problemData.numSteps, problemData.outputTypes);
end % if

fprintf('Finished simulation at t_end = %g\n', problemData.tEnd);
%% Error evaluation
if isfield(problemData, 'cCont')
  problemData.error = computeL2Error(problemData.g, problemData.cDisc, ...
    @(x1, x2) problemData.cCont(problemData.tEnd, x1, x2), problemData.qOrd, problemData.basesOnQuad);
  fprintf('L2 error w.r.t. the analytical solution: %g\n', problemData.error);
else
  problemData.error = computeL2Error(problemData.g, problemData.cDisc, ...
    problemData.c0Cont, problemData.qOrd, problemData.basesOnQuad);
  fprintf('L2 error w.r.t. the initial condition: %g\n', problemData.error);
end % if
end % function

\end{lstlisting}

\subsection{Helper routines}

\code{ret = assembleMatEdgeMuMu(g, markE0T, refEdgeMuMu)} assembles mass-matrices~$\matMmuBar$ and~$\matMmuTilde$ for the edge-based basis functions~$\muh$ as described in Sec.~\ref{sec:assembly:matMmuBarAndMatMuTilde}.
Array \code{markE0T} plays the role of the Kronecker delta in the matrix definition by selecting the relevant edges for either matrix. 
\code{refEdgeMuMu} is the mass matrix~$\matMhat_\mu$ on the reference interval~$[0,1]$ computed by~\code{integrateRefEdgeMuMu}.
\begin{lstlisting}
function ret = assembleMatEdgeMuMu(g, markE0T, refEdgeMuMu)
Nmu = size(refEdgeMuMu, 1); Kedge = g.numE;
ret = sparse(Kedge * Nmu, Kedge * Nmu);
for n = 1 : 3
  Kkn = g.areaE0T(:, n) .*  markE0T(:, n) ;
  ret = ret + kron(sparse(g.E0T(:, n), g.E0T(:, n), Kkn, Kedge, Kedge), refEdgeMuMu);
end % for
end % function
\end{lstlisting}

\code{ret = assembleMatEdgePhiIntMu(g, markE0T, refEdgePhiIntMu)} assembles matrices~$\matRmu$,~$\matT$, and~$\matKmuOut$ as detailed in Sec.~\ref{sec:assembly:matRmuAndmatT}.
As before, \code{markE0T} is used to select relevant edges, and \code{refEdgePhiIntMu} stores the pre-computed reference blocks~$\matRmuHat$ given by~\code{integrateRefEdgePhiIntMu}.
\begin{lstlisting}
function ret = assembleMatEdgePhiIntMu(g, markE0T, refEdgePhiIntMu)
K = g.numT;  Kedge = g.numE; 
[N, Nmu, ~, ~] = size(refEdgePhiIntMu);
ret = sparse(K * N, Kedge * Nmu);
for n = 1 : 3
  for l = 1 : 2
    Rkn = markE0T(:,n) .*  g.markSideE0T(:, n, l) .* g.areaE0T(:, n);
    ret = ret + kron(sparse(1 : K, g.E0T(:, n), Rkn, K, Kedge), refEdgePhiIntMu(:, :, n, l));
  end % for
end % for
end % function
\end{lstlisting}

\code{ret = assembleMatEdgePhiIntMuVal(g, markE0T, refEdgePhiIntMuPerQuad, valOnQuad)} builds the global matrices~$\matS$ and~$\matSout$ given in Sec.~\ref{sec:assembly:matSmAndmatSout}.
Both matrices are assembled at once by specifying all interior and outflow edges as~\code{markE0T}.
The reference blocks per quadrature point~$\matShat$ are computed by~$\code{integrateRefEdgePhiIntMuPerQuad}$ and specified in the parameter~\code{refEdgePhiIntMuPerQuad}.
Parameter~\code{valOnQuad} stores the normal velocity~$\vecU_{\vNormal}$ evaluated at the quadrature points of each edge.
\begin{lstlisting}
function ret = assembleMatEdgePhiIntMuVal(g, markE0T, refEdgePhiIntMuPerQuad, valOnQuad)
K = g.numT; Kbar = g.numE;
[N, Nmu, ~, R] = size(refEdgePhiIntMuPerQuad{1});

ret = sparse(K*N, Kbar*Nmu);
for n = 1 : 3
  RknTimesVal = sparse(K*N, Nmu);
  markAreaE0T = markE0T(:, n) .* g.areaE0T(:, n);
  for l = 1 : 2
    markAreaSideE0T = markAreaE0T .*  g.markSideE0T(:, n, l);
    for r = 1 : R
      RknTimesVal = RknTimesVal + kron(markAreaSideE0T .* valOnQuad(:, n, r), ...
                      refEdgePhiIntMuPerQuad{l}(:, :, n, r));
    end % for r
  end % for l
  ret = ret + kronVec(sparse(1:K, g.E0T(:, n), ones(K, 1), K, Kbar), RknTimesVal);
end % for n
end % function
\end{lstlisting}

\code{ret = assembleMatElemDphiPhiFuncContVec(g, refElemDphiPhiPerQuad, funcCont1, funcCont2, qOrd)} assembles matrices~$\matG^m$ described in Sec.~\ref{sec:assembly:matG}.
Reference blocks~$\matGhat$ are provided in parameter~\code{refElemDphiPhiPerQuad} and computed by~\code{integrateRefElemDphiPhiPerQuad}.
The two components of the continuous function in the integrand -- here the velocity components -- are given as function handles in~\code{funcCont1} and~\code{funcCont2}, respectively.
\begin{lstlisting}
function ret = assembleMatElemDphiPhiFuncContVec(g, refElemDphiPhiPerQuad, funcCont1, funcCont2, qOrd)
K = g.numT; [N, ~, R] = size(refElemDphiPhiPerQuad{1});

if nargin < 5, p = (sqrt(8*N+1)-3)/2;  qOrd = max(2*p, 1); end
[Q1, Q2, ~] = quadRule2D(qOrd);

ret = { zeros(K*N, N), zeros(K*N,N) };
for r = 1 : R
  valOnQuad1 = funcCont1(g.mapRef2Phy(1, Q1(r), Q2(r)), g.mapRef2Phy(2, Q1(r), Q2(r)));
  ret{1} = ret{1} + kron(g.B(:,2,2) .* valOnQuad1, refElemDphiPhiPerQuad{1}(:, :, r)) ...
                  - kron(g.B(:,2,1) .* valOnQuad1, refElemDphiPhiPerQuad{2}(:, :, r));
        
  valOnQuad2 = funcCont2(g.mapRef2Phy(1, Q1(r), Q2(r)), g.mapRef2Phy(2, Q1(r), Q2(r)));
  ret{2} = ret{2} - kron(g.B(:,1,2) .* valOnQuad2, refElemDphiPhiPerQuad{1}(:, :, r)) ...
                  + kron(g.B(:,1,1) .* valOnQuad2, refElemDphiPhiPerQuad{2}(:, :, r));
end % for
ret{1} = kronVec(speye(K,K), ret{1});
ret{2} = kronVec(speye(K,K), ret{2});
end % function
\end{lstlisting}

\code{ret = assembleVecEdgeMuFuncCont(g, markE0T, funcCont, basesOnQuad, qOrd)} builds the vector with Dirichlet data~$\vecKmuD$ (cf. Sec.~\ref{sec:assembly:vecKmuD}).
Inflow edges are specified in~\code{markE0T}, and Dirichlet data is specified as a~function handle in~\code{funcCont}.
\begin{lstlisting}
function ret = assembleVecEdgeMuFuncCont(g, markE0T, funcCont, basesOnQuad, qOrd)
[Q, W] = quadRule1D(qOrd);
[R, Nmu] = size(basesOnQuad.mu{qOrd});

ret = zeros(g.numE, Nmu);
for n = 1 : 3
  [Q1, Q2] = gammaMap(n, Q);
  funcOnQuad = funcCont(g.mapRef2Phy(1, Q1, Q2), g.mapRef2Phy(2, Q1, Q2));
  Kkn = markE0T(:, n) .* g.areaE0T(:,n);
  ret(g.E0T(:, n), :) = ret(g.E0T(:, n), :) + (repmat(Kkn, 1, R) .* funcOnQuad) ...
                          * ( repmat(W', 1, Nmu) .* basesOnQuad.mu{qOrd} );
end % for

ret = reshape(ret', [], 1);
end % function
\end{lstlisting}


\code{invA = blkinv(A, blockSize)} inverts a~square block-diagonal matrix~\code{A} built up from blocks of size~$n\times n$, which is used to compute~$\matL^{-1}$ as described in Sec.~\ref{sec:staticcondensation}.
To improve the computational performance of this method, the inversion is performed using a~blocking-technique where a~pre-specified number of diagonal blocks are inverted together.
Parameter~\code{blockSize} specifies the blocking size and must be a~multiple of~$n$.
\begin{lstlisting}
function invA = blkinv(A, blockSize)
n = size(A,1);
numBlocks = ceil(n / blockSize);
idxEnd = (numBlocks - 1) * blockSize;
idxs = mat2cell(1 : idxEnd, 1, ones(1, numBlocks - 1) * blockSize);
cellA = cellfun(@(X) A(X, X), idxs, 'UniformOutput', false);
cellInvA = cellfun(@(X) inv(X), cellA, 'UniformOutput', false);
invA = blkdiag(cellInvA{:}, inv(A(idxEnd + 1 : end, idxEnd + 1 : end)));
end % function
\end{lstlisting}

\code{basesOnQuadEdge = computeBasesOnQuadEdge(N, basesOnQuadEdge, requiredOrders)} evaluates edge basis function~$\muh$ at the quadrature points of the reference interval~$[0,1]$.
Parameter \code{basesOnQuadEdge} is a~(possibly empty) struct to which the computed fields~\code{mu} and~\code{thetaMu} (the latter corresponding to~$\mhat{} \circ \betaMap$) are added.
The quadrature orders can be optionally specified in~\code{requiredOrders} and with defaults equal to~$2p$ and~$2p+1$.
\begin{lstlisting}
function basesOnQuadEdge = computeBasesOnQuadEdge(N, basesOnQuadEdge, requiredOrders)
if nargin < 3
    p = N - 1;
    if p > 0
        requiredOrders = [2*p, 2*p+1];
    else
        requiredOrders = 1;
    end % if
end % if

basesOnQuadEdge.mu = cell(max(requiredOrders),1);
basesOnQuadEdge.thetaMu = cell(max(requiredOrders),1);
for it = 1 : length(requiredOrders)
    ord = requiredOrders(it);
    [Q, ~] = quadRule1D(ord);
    R = length(Q);
    
    basesOnQuadEdge.mu{ord} = zeros(R, N);
    for i = 1 : N
        basesOnQuadEdge.mu{ord}(:, i) = phi1D(i, Q);
    end
    
    basesOnQuadEdge.thetaMu{ord} = zeros(R, N, 2);
    basesOnQuadEdge.thetaMu{ord}(:, :, 1) = basesOnQuadEdge.mu{ord};
    basesOnQuadEdge.thetaMu{ord}(:, :, 2) = flipud(basesOnQuadEdge.mu{ord});
end % for
end % function
\end{lstlisting}

\code{g = computeDerivedGridData(g)} enriches grid data structure~\code{g} by an~additional~$K\times 3\times 2$ field~\code{markSideE0T} that marks for each element-local edge, whether the element has local index~1 or~2 at the edge thus providing a~way to determine~$l$ in the mapping~$\kappa$ (cf. eq.~\eqref{eq:mapKappa}) for given element~$T_k$ and edge~$\Ekn = E_{\kEdge}$.
\begin{lstlisting}
function g = computeDerivedGridData(g)
g.markSideE0T = true(g.numT, 3, 2);
for n = 1 : 3
  % Mark element to have local id 1 or 2 at the given edge
  g.markSideE0T(:, n, 1) = g.T0E(g.E0T(:, n), 2) ~= (1:g.numT)';
  g.markSideE0T(:, n, 2) = ~g.markSideE0T(:, n, 1);
end % for
end % function
\end{lstlisting}

\code{ret = computeFuncContOnQuadEdge(g, funcCont, qOrd)} evaluates a~given function handle~\code{funcCont} in all quadrature points of all edges.
This is used to determine the Dirichlet boundary data~$c_\mathrm{D}$ before multiplying it with the normal velocity and assembling~$\vecFphiD$ in routine \code{assembleVecEdgePhiIntVal}.
\begin{lstlisting}
function ret = computeFuncContOnQuadEdge(g, funcCont, qOrd)
K = g.numT;  [Q, W] = quadRule1D(qOrd);
ret = zeros(K, 3, length(W));
for n = 1 : 3
  [Q1, Q2] = gammaMap(n, Q);
  ret(:, n, :) = funcCont(g.mapRef2Phy(1, Q1, Q2), g.mapRef2Phy(2, Q1, Q2));
end % for
end % function
\end{lstlisting}

\code{problemData = getTestcase(problemData, testcase)} defines initial and boundary data as well as the analytical solution, if available, for test cases given in Sec.~\ref{sec:results}.
To shorten the presentation, the actual definitions are omitted here but can be found in GitHub~\cite{FESTUNGGithub}.
\begin{lstlisting}
function problemData = getTestcase(problemData, testcase)  
switch testcase
  case 'solid_body' 
    % LeVeque's solid body rotation 
  case 'stationary' 
    % Stationary analytical example    
  case 'transient_ode'
    % Transient ODE example    
  case 'transient' 
    % Transient analytical example   
  otherwise
    error('Invalid testcase "%s".', testcase);
end % switch

fprintf('Loaded testcase "%s".\n', testcase);

% Time stepping parameters
problemData = setdefault(problemData, 'numSteps', numSteps);
problemData = setdefault(problemData, 'tEnd', tEnd);
problemData = setdefault(problemData, 'isStationary', isStationary);

% Store defined functions in struct
problemData = setdefault(problemData, 'c0Cont', c0Cont);
problemData = setdefault(problemData, 'fCont', fCont);
problemData = setdefault(problemData, 'u1Cont', u1Cont);
problemData = setdefault(problemData, 'u2Cont', u2Cont);
problemData = setdefault(problemData, 'cDCont', cDCont);
if isAnalytical, problemData = setdefault(problemData, 'cCont', cCont); end

% Specify boundary conditions
problemData = setdefault(problemData, 'generateMarkE0Tint', @(g) g.idE0T == 0);
problemData = setdefault(problemData, 'generateMarkE0Tbdr', @(g) g.idE0T ~= 0);

end % function
\end{lstlisting}

\code{ret = integrateRefEdgeMuMu(N, basesOnQuadEdge, qOrd)} computes reference blocks~$\matMhat_\mu$ (see Sec.~\ref{sec:assembly:matMmuBarAndMatMuTilde}).
\begin{lstlisting}
function ret = integrateRefEdgeMuMu(N, basesOnQuad, qOrd)
if nargin < 3, p = N-1;  qOrd = 2*p+1;  end
[~, W] = quadRule1D(qOrd);
ret = zeros(N, N); % [N x N]
for i = 1 : N
  for j = 1 : N
    ret(i, j) = sum( W' .* basesOnQuad.mu{qOrd}(:,i) .* basesOnQuad.mu{qOrd}(:,j) );
  end % for
end % for
end % function
\end{lstlisting}

\code{ret = integrateRefEdgePhiIntMu(N, basesOnQuad, qOrd)} computes reference blocks~$\matRmuHat$ (see Sec.~\ref{sec:assembly:matRmuAndmatT}).
\begin{lstlisting}
function ret = integrateRefEdgePhiIntMu(N, basesOnQuad, qOrd)
if nargin < 3, p = (sqrt(8*N(1)+1)-3)/2;  qOrd = 2*p+1;  end
[~, W] = quadRule1D(qOrd);

ret = zeros(N(1), N(2), 3, 2); % [N(1) x N(2) x 3 x 2]
for n = 1 : 3
    for i = 1 : N(1)
        for j = 1 : N(2)
            ret(i, j, n, 1) =  W * (basesOnQuad.phi1D{qOrd}(:, i, n) .* basesOnQuad.mu{qOrd}(:, j) );
            ret(i, j, n, 2) =  W * (basesOnQuad.phi1D{qOrd}(:, i, n) .* basesOnQuad.thetaMu{qOrd}(:, j, 2) );
        end % for
    end % for
end % for
end % function
\end{lstlisting}

\code{ret = integrateRefEdgePhiIntMuPerQuad(N, basesOnQuad, qOrd)} computes reference blocks~$\matShat$ (see Sec.~\ref{sec:assembly:matSmAndmatSout}).
\begin{lstlisting}
function ret = integrateRefEdgePhiIntMuPerQuad(N, basesOnQuad, qOrd)
if nargin < 3, p = (sqrt(8*N+1)-3)/2;  qOrd = 2*p+1; end
[~, W] = quadRule1D(qOrd); R = length(W);
ret = { zeros(N(1), N(2), 3, R); zeros(N(1), N(2), 3, R) };
for n = 1 : 3 
  for i = 1 : N(1)
    for j = 1 : N(2)
      ret{1}(i, j, n, :) = W(:) .* basesOnQuad.phi1D{qOrd}(:, i, n) ...
                              .* basesOnQuad.mu{qOrd}(:, j);
      ret{2}(i, j, n, :) = W(:) .* basesOnQuad.phi1D{qOrd}(:, i, n) ...
                              .* basesOnQuad.thetaMu{qOrd}(:, j, 2);
    end % for j
  end % for i
end % for n
end % function
\end{lstlisting}

\code{ret = integrateRefElemDphiPhiPerQuad(N, basesOnQuad, qOrd)} computes reference blocks~$\matGhat$ (see Sec.~\ref{sec:assembly:matG}).
\begin{lstlisting}
function ret = integrateRefElemDphiPhiPerQuad(N, basesOnQuad, qOrd)
if nargin < 3, p = (sqrt(8*N+1)-3)/2;  qOrd = max(2*p, 1);  end
[~, ~, W] = quadRule2D(qOrd); R = length(W);
ret = { zeros(N, N, R); zeros(N, N, R) };
for i = 1 : N
  for j = 1 : N
    for m = 1 : 2
      ret{m}(i, j, :) =  W(:) .* basesOnQuad.phi2D{qOrd}(:, j) .* basesOnQuad.gradPhi2D{qOrd}(:, i, m);
    end % for
  end % for
end % for
end % function
\end{lstlisting}

\code{[t, A, b, c] = rungeKuttaImplicit(ord, tau, t0)} provides the Butcher tableau for DIRK schemes of orders one to four (see Sec.~\ref{sec:timeDiscretization}) with the order given in~\code{ord}.
The current time level is provided in~\code{t0}, and the time step size is given in~\code{tau}.
\begin{lstlisting}
function [t, A, b, c] = rungeKuttaImplicit(ord, tau, t0)
switch ord
  case 1
    A = 1;
    
  case 2
    lambda = 0.5 * (2 - sqrt(2));
    A = [     lambda,      0 ; 
          (1-lambda), lambda ];
    
  case 3
    alph = 0.4358665215084589994160194511935568425292;
    bet  = 0.5 * (1.+ alph);
    b1    = -0.25 * (6 * alph * alph - 16 * alph + 1);
    b2    =  0.25 * (6 * alph * alph - 20 * alph + 5);

    A = [       alph,     0,     0 ; 
          bet - alph,  alph,     0 ; 
                  b1,    b2,  alph ];
    
  case 4
    gamm = 1./4.; 
    A = [      gamm,           0,        0,        0,     0;
              1./2.,        gamm,        0,        0,     0;
            17./50.,     -1./25.,     gamm,        0,     0;
         371./1360., -137./2720., 15./544.,     gamm,     0;
            25./24.,    -49./48., 125./16., -85./12.,  gamm ];
    
  otherwise
    error('Invalid order for implicit Runge-Kutta')
end % switch
    
b = A(end, :);
c = sum(A, 2);
t = t0 + c * tau;
end % function
\end{lstlisting}

%% file: FESTUNG-Pt3.bbl
\begin{thebibliography}{10}
\expandafter\ifx\csname url\endcsname\relax
  \def\url#1{\texttt{#1}}\fi
\expandafter\ifx\csname urlprefix\endcsname\relax\def\urlprefix{URL }\fi
\expandafter\ifx\csname href\endcsname\relax
  \def\href#1#2{#2} \def\path#1{#1}\fi

\bibitem{ReedHill1973}
H.~Reed, T.~R. Hill, Triangular mesh methods for the neutron transport
  equation, Tech. Rep. LA-UR-73-479, Los Alamos Scientific Laboratory, NM
  (1973).

\bibitem{ABCM}
D.~N. Arnold, F.~Brezzi, B.~Cockburn, L.~D. Marini, Unified analysis of
  {discontinuous {{Galerkin}}} methods for elliptic problems, SIAM Journal on
  Numerical Analysis 39 (2002) 1749--1779.

\bibitem{ShuOverviewDG}
C.-W. Shu, Discontinuous {G}alerkin method for time-dependent problems: survey
  and recent developments, in: Recent developments in discontinuous {G}alerkin
  finite element methods for partial differential equations, Vol. 157 of IMA
  Vol. Math. Appl., Springer, Cham, 2014, pp. 25--62.
\newblock \href {http://dx.doi.org/10.1007/978-3-319-01818-8_2}
  {\path{doi:10.1007/978-3-319-01818-8_2}}.

\bibitem{BaBo2011}
F.~Bassi, L.~Botti, A.~Colombo, D.~D. Pietro, P.~Tesini, On the flexibility of
  agglomeration based physical space discontinuous {Galerkin} discretizations,
  Journal of Computational Physics 231 (2011) 45--65.

\bibitem{Wang201553}
L.~Wang, P.-O. Persson, A high-order discontinuous {G}alerkin method with
  unstructured space–time meshes for two-dimensional compressible flows on
  domains with large deformations, Computers \& Fluids 118 (2015) 53 -- 68.
\newblock \href {http://dx.doi.org/10.1016/j.compfluid.2015.05.026}
  {\path{doi:10.1016/j.compfluid.2015.05.026}}.

\bibitem{ortwein201401}
P.~Ortwein, T.~Binder, S.~Copplestone, A.~Mirza, P.~Nizenkov, M.~Pfeiffer,
  T.~Stindl, S.~Fasoulas, C.-D. Munz, Parallel performance of a discontinuous
  {G}alerkin spectral element method based {PIC}-{DSMC} solver, in: W.~E.
  Nagel, D.~H. Kröner, M.~M. Resch (Eds.), High Performance Computing in
  Science and Engineering ‘14, Springer, 2015, pp. 671--681.
\newblock \href {http://dx.doi.org/10.1007/978-3-319-10810-0_44}
  {\path{doi:10.1007/978-3-319-10810-0_44}}.

\bibitem{Fidkowski2005}
K.~Fidkowski, T.~Oliver, J.~Lu, D.~L. Darmofal, p--{Multigrid} solution of
  high-order {discontinuous {Galerkin}} discretizations of the compressible
  {Navier--Stokes} equations, Journal of Computational Physics 207 (2005)
  92--113.

\bibitem{Luo2006}
H.~Luo, J.~D. Baum, R.~L\"{o}hner,
  \href{http://dx.doi.org/10.1016/j.jcp.2005.06.019}{A p-multigrid
  discontinuous {G}alerkin method for the {E}uler equations on unstructured
  grids}, Journal of Computational Physics 211~(2) (2006) 767--783.
\newblock \href {http://dx.doi.org/10.1016/j.jcp.2005.06.019}
  {\path{doi:10.1016/j.jcp.2005.06.019}}.
\newline\urlprefix\url{http://dx.doi.org/10.1016/j.jcp.2005.06.019}

\bibitem{Bassi2009}
F.~Bassi, A.~Ghidoni, S.~Rebay, P.~Tesini,
  \href{http://dx.doi.org/10.1002/fld.1917}{High-order accurate p-multigrid
  discontinuous {G}alerkin solution of the euler equations}, International
  Journal for Numerical Methods in Fluids 60~(8) (2009) 847--865.
\newblock \href {http://dx.doi.org/10.1002/fld.1917}
  {\path{doi:10.1002/fld.1917}}.
\newline\urlprefix\url{http://dx.doi.org/10.1002/fld.1917}

\bibitem{AizingerKK2015}
V.~Aizinger, D.~Kuzmin, L.~Korous,
  \href{http://www.sciencedirect.com/science/article/pii/S0096300315006645}{Scale
  separation in fast hierarchical solvers for discontinuous {G}alerkin
  methods}, Applied Mathematics and Computation 266 (2015) 838--849.
\newblock \href {http://dx.doi.org/10.1016/j.amc.2015.05.047}
  {\path{doi:10.1016/j.amc.2015.05.047}}.
\newline\urlprefix\url{http://www.sciencedirect.com/science/article/pii/S0096300315006645}

\bibitem{JaustSA2016}
A.~Jaust, J.~Sch\"{u}tz, V.~Aizinger,
  \href{http://dx.doi.org/10.1002/pamm.201610411}{An efficient linear solver
  for the hybridized discontinuous {G}alerkin method}, PAMM 16~(1) (2016)
  845--846.
\newblock \href {http://dx.doi.org/10.1002/pamm.201610411}
  {\path{doi:10.1002/pamm.201610411}}.
\newline\urlprefix\url{http://dx.doi.org/10.1002/pamm.201610411}

\bibitem{SchuetzAizinger2017}
J.~Sch\"{u}tz, V.~Aizinger,
  \href{http://www.sciencedirect.com/science/article/pii/S0377042716306288}{A
  hierarchical scale separation approach for the hybridized discontinuous
  {G}alerkin method}, Journal of Computational and Applied Mathematics 317
  (2017) 500--509.
\newblock \href {http://dx.doi.org/10.1016/j.cam.2016.12.018}
  {\path{doi:10.1016/j.cam.2016.12.018}}.
\newline\urlprefix\url{http://www.sciencedirect.com/science/article/pii/S0377042716306288}

\bibitem{Thiele2017}
C.~Thiele, M.~Araya-Polo, F.~O. Alpak, B.~Riviere, F.~Frank,
  \href{http://www.sciencedirect.com/science/article/pii/S0898122117303735}{Inexact
  hierarchical scale separation: A two-scale approach for linear systems from
  discontinuous {G}alerkin discretizations}, Computers and Mathematics with
  Applications (2017) in press\href
  {http://dx.doi.org/http://dx.doi.org/10.1016/j.camwa.2017.06.025}
  {\path{doi:http://dx.doi.org/10.1016/j.camwa.2017.06.025}}.
\newline\urlprefix\url{http://www.sciencedirect.com/science/article/pii/S0898122117303735}

\bibitem{CGorHDG}
R.~M. Kirby, S.~J. Sherwin, B.~Cockburn, To {CG} or {HDG}: {A} comparative
  study, Journal of Scientific Computing 51~(1) (2012) 183--212.
\newblock \href {http://dx.doi.org/10.1007/s10915-011-9501-7}
  {\path{doi:10.1007/s10915-011-9501-7}}.

\bibitem{Yakovlev2016}
S.~Yakovlev, D.~Moxey, R.~M. Kirby, S.~J. Sherwin, To {CG} or to {HDG}: {A}
  comparative study in 3d, Journal of Scientific Computing 67~(1) (2016)
  192--220.
\newblock \href {http://dx.doi.org/10.1007/s10915-015-0076-6}
  {\path{doi:10.1007/s10915-015-0076-6}}.

\bibitem{WBMS13}
M.~Woopen, A.~Balan, G.~May, J.~Sch\"utz, A comparison of hybridized and
  standard {DG} methods for target-based $hp$-adaptive simulation of
  compressible flow, Computers and Fluids 98 (2014) 3--16.

\bibitem{Veubeke}
B.~F. de~Veubeke, Displacement and equilibrium models in the finite element
  method., in: O.C.Zienkiewicz, G.S.Holister (Eds.), Stress Analysis, Wiley New
  York, 1965, Ch.~9, pp. 145--197.

\bibitem{AB85}
D.~Arnold, F.~Brezzi, Mixed and nonconforming {Finite Element} methods:
  Implementation, postprocessing and error estimates, Mathematical Modelling
  and Numerical Analysis 19 (1985) 7--32.

\bibitem{BDM85}
F.~Brezzi, J.~Douglas, L.~D. Marini, Two families of mixed finite elements for
  second order elliptic problems, Numerische Mathematik 47 (1985) 217--235.

\bibitem{CoGo04}
B.~Cockburn, J.~Gopalakrishnan, A characterization of hybridized mixed methods
  for second order elliptic problems, SIAM Journal on Numerical Analysis 42
  (2004) 283--301.

\bibitem{COGOLA}
B.~Cockburn, J.~Gopalakrishnan, R.~Lazarov, Unified hybridization of
  {discontinuous {{Galerkin}}}, mixed, and continuous {{{Galerkin}}} methods
  for second order elliptic problems, SIAM Journal on Numerical Analysis 47
  (2009) 1319--1365.

\bibitem{CockburnHDGStokes}
B.~Cockburn, J.~Gopalakrishnan, N.~C. Nguyen, J.~Peraire, F.-J. Sayas, Analysis
  of {HDG} methods for {S}tokes flow, Mathematics of Computation 80~(274)
  (2011) 723--760.
\newblock \href {http://dx.doi.org/10.1090/S0025-5718-2010-02410-X}
  {\path{doi:10.1090/S0025-5718-2010-02410-X}}.

\bibitem{EgWal13}
H.~Egger, C.~Waluga, hp-analysis of a hybrid {DG} method for {Stokes} flow, IMA
  Journal of Numerical Analysis 33~(2) (2013) 687--721.

\bibitem{Egger2013}
H.~Egger, C.~Waluga, A hybrid discontinuous {G}alerkin method for
  {Darcy}-{Stokes} problems, in: R.~Bank, M.~Holst, O.~Widlund, J.~Xu (Eds.),
  Domain Decomposition Methods in Science and Engineering XX, Springer Berlin
  Heidelberg, Berlin, Heidelberg, 2013, pp. 663--670.
\newblock \href {http://dx.doi.org/10.1007/978-3-642-35275-1_79}
  {\path{doi:10.1007/978-3-642-35275-1_79}}.

\bibitem{NguyenNavierStokes}
N.~Nguyen, J.~Peraire, B.~Cockburn, A hybridizable discontinuous {G}alerkin
  method for the incompressible {N}avier-{S}tokes equations, AIAA Paper
  2010-362.

\bibitem{NgPe12}
N.~C. Nguyen, J.~Peraire, Hybridizable discontinuous {{Galerkin}} methods for
  partial differential equations in continuum mechanics, Journal of
  Computational Physics 231 (2012) 5955--5988.

\bibitem{NgPeCo11}
N.~C. Nguyen, J.~Peraire, B.~Cockburn, High-order implicit hybridizable
  discontinuous {{Galerkin}} methods for acoustics and elastodynamics, Journal
  of Computational Physics 230 (2011) 3695--3718.

\bibitem{SchMa11}
J.~Sch\"utz, G.~May, A hybrid mixed method for the compressible {Navier-Stokes}
  equations, Journal of Computational Physics 240 (2013) 58--75.

\bibitem{NguyenMaxwell}
N.~C. Nguyen, J.~Peraire, B.~Cockburn, Hybridizable discontinuous {G}alerkin
  methods for the time-harmonic {M}axwell's equations, Journal of Computational
  Physics 230~(19) (2011) 7151--7175.
\newblock \href {http://dx.doi.org/10.1016/j.jcp.2011.05.018}
  {\path{doi:10.1016/j.jcp.2011.05.018}}.

\bibitem{NPC09}
N.~C. Nguyen, J.~Peraire, B.~Cockburn, An implicit high-order hybridizable
  {discontinuous {{Galerkin}}} method for nonlinear convection-diffusion
  equations, Journal of Computational Physics 228 (2009) 8841--8855.

\bibitem{NPC09L}
N.~C. Nguyen, J.~Peraire, B.~Cockburn, An implicit high-order hybridizable
  discontinuous {{Galerkin}} method for linear convection-diffusion equations,
  Journal of Computational Physics 228 (2009) 3232--3254.

\bibitem{EgSch09}
H.~Egger, J.~Sch\"oberl, A hybrid mixed {discontinuous {{Galerkin}}} finite
  element method for convection-diffusion problems, IMA Journal of Numerical
  Analysis 30 (2010) 1206--1234.

\bibitem{BuiThanh2015}
T.~Bui-Thanh, From {G}odunov to a unified hybridized discontinuous {G}alerkin
  framework for partial differential equations, Journal of Computational
  Physics 295 (2015) 114--146.
\newblock \href {http://dx.doi.org/10.1016/j.jcp.2015.04.009}
  {\path{doi:10.1016/j.jcp.2015.04.009}}.

\bibitem{FrankRAK2015}
F.~Frank, B.~Reuter, V.~Aizinger, P.~Knabner, {FESTUNG}: A {MATLAB}/{GNU}
  {O}ctave toolbox for the discontinuous {G}alerkin method. {P}art {I}:
  {D}iffusion operator, Computers and Mathematics with Applications 70~(1)
  (2015) 11 -- 46.
\newblock \href {http://dx.doi.org/10.1016/j.camwa.2015.04.013}
  {\path{doi:10.1016/j.camwa.2015.04.013}}.

\bibitem{ReuterAWFK2016}
B.~Reuter, V.~Aizinger, M.~Wieland, F.~Frank, P.~Knabner, {FESTUNG}: A
  {MATLAB}/{GNU} {O}ctave toolbox for the discontinuous {G}alerkin method.
  {P}art {II}: {A}dvection operator and slope limiting, Computers and
  Mathematics with Applications 72~(7) (2016) 1896--1925.
\newblock \href {http://dx.doi.org/10.1016/j.camwa.2016.08.006}
  {\path{doi:10.1016/j.camwa.2016.08.006}}.

\bibitem{FESTUNG}
B.~Reuter, F.~Frank, V.~Aizinger,
  \href{https://www.math.fau.de/FESTUNG}{{FESTUNG} --- {T}he {F}inite {E}lement
  {S}imulation {T}oolbox for {UN}structured {G}rids} (2018).
\newline\urlprefix\url{https://www.math.fau.de/FESTUNG}

\bibitem{FESTUNGGithub}
B.~Reuter, F.~Frank, \href{https://github.com/FESTUNG}{{FESTUNG}: {T}he
  {F}inite {E}lement {S}imulation {T}oolbox for {UN}structured {G}rids} (2018).
\newblock \href {http://dx.doi.org/10.5281/zenodo.1215561}
  {\path{doi:10.5281/zenodo.1215561}}.
\newline\urlprefix\url{https://github.com/FESTUNG}

\bibitem{Kuzmin2010}
D.~Kuzmin, A vertex-based hierarchical slope limiter for adaptive discontinuous
  {G}alerkin methods, Journal of Computational and Applied Mathematics 233~(12)
  (2010) 3077--3085, {F}inite Element Methods in Engineering and Science
  (FEMTEC 2009).
\newblock \href {http://dx.doi.org/10.1016/j.cam.2009.05.028}
  {\path{doi:10.1016/j.cam.2009.05.028}}.

\bibitem{Aizinger2011}
V.~Aizinger, A geometry independent slope limiter for the discontinuous
  {G}alerkin method, in: E.~Krause, Y.~Shokin, M.~Resch, D.~Kr{\"o}ner,
  N.~Shokina (Eds.), Computational Science and High Performance Computing IV,
  Vol. 115 of Notes on Numerical Fluid Mechanics and Multidisciplinary Design,
  Springer Berlin Heidelberg, 2011, pp. 207--217.
\newblock \href {http://dx.doi.org/10.1007/978-3-642-17770-5_16}
  {\path{doi:10.1007/978-3-642-17770-5_16}}.

\bibitem{ReuterRAK2017}
B.~Reuter, A.~Rupp, V.~Aizinger, P.~Knabner, {FESTUNG}: A {MATLAB}/{GNU}
  {O}ctave toolbox for the discontinuous {G}alerkin method. {P}art {IV}:
  {G}eneric problem framework and model coupling interface, in preparation.

\bibitem{GR1}
E.~Godlewski, P.-A. Raviart, Hyperbolic Systems of {Conservation Laws},
  Ellipses Paris, 1991.

\bibitem{GR2}
E.~Godlewski, P.-A. Raviart, Numerical Approximation of Hyperbolic Systems of
  {Conservation Laws}, Springer New York, 1996.

\bibitem{Kroener}
D.~Kr\"oner, Numerical Schemes for Conservation Laws, Wiley Teubner, 1997.

\bibitem{LF}
P.~Lax, Weak solutions of nonlinear hyperbolic equations and their numerical
  computation, Communications on Pure and Applied Mathematics 7 (1954)
  159--193.

\bibitem{Rusanov}
V.~V. Rusanov, The calculation of the interaction of non-stationary shock waves
  with barriers, Journal of Computational Mathematics and Mathematical Physics
  (in Russian) 1 (1961) 267--279.

\bibitem{Alex77}
R.~Alexander, Diagonally implicit {{Runge-Kutta}} methods for stiff
  {O}.{D}.{E}.'s., SIAM Journal of Numerical Analysis 14 (1977) 1006--1021.

\bibitem{HaiWanII}
E.~Hairer, G.~Wanner, Solving Ordinary Differential Equations II, Springer
  Series in Computational Mathematics, 1991.

\bibitem{JS13}
A.~Jaust, J.~Sch\"utz, A temporally adaptive hybridized discontinuous
  {{Galerkin}} method for time-dependent compressible flows, Computers and
  Fluids 98 (2014) 177--185.
\newblock \href {http://dx.doi.org/10.1016/j.compfluid.2014.01.019}
  {\path{doi:10.1016/j.compfluid.2014.01.019}}.

\bibitem{LeVeque1996}
R.~J. Le{V}eque, High-resolution conservative algorithms for advection in
  incompressible flow, SIAM Journal on Numerical Analysis 33~(2) (1996)
  627--665.
\newblock \href {http://dx.doi.org/10.2307/2158391}
  {\path{doi:10.2307/2158391}}.

\bibitem{Badia2017}
S.~Badia, J.~Bonilla, A.~Hierro, {D}ifferentiable monotonicity-preserving
  schemes for discontinuous {G}alerkin methods on arbitrary meshes, Computer
  Methods in Applied Mechanics and Engineering 320 (2017) 582 -- 605.
\newblock \href {http://dx.doi.org/10.1016/j.cma.2017.03.032}
  {\path{doi:10.1016/j.cma.2017.03.032}}.

\end{thebibliography}
